\documentclass[12pt]{article}
\usepackage[utf8]{inputenc}
\usepackage{cite}

\usepackage{url}

\usepackage{caption}

\usepackage{amssymb, amsmath, amsthm}
\usepackage{tikz}

\usepackage{thm-restate}

\usepackage{cellspace}
\usepackage[pdftex,hypertexnames=false,linktocpage=true,hidelinks]{hyperref}
\usepackage{cleveref}

\usepackage{fullpage}

\usepackage{amsfonts}

\usepackage{enumitem}

\newtheorem{theorem}{Theorem}[section]

\newtheorem{lemma}[theorem]{Lemma}

\newtheorem{cor}[theorem]{Corollary}

\newtheorem{prop}[theorem]{Proposition}
\newtheorem{claim}[theorem]{Claim}
\newtheorem{observation}[theorem]{Observation}
\theoremstyle{definition}
\newtheorem{definition}[theorem]{Definition}

\newtheorem{remark}[theorem]{Remark}

\newtheorem{ques}{Question}

\numberwithin{equation}{section}

\newcommand{\acopyof}{}

\tikzset{vtx/.style={inner sep=2.7pt, outer sep=0pt, circle, fill=red,draw}}

\newcommand{\vc}[1]{\ensuremath{\vcenter{\hbox{#1}}}}
\tikzset{flag_pic/.style={scale=1}}  
\tikzset{unlabeled_vertex/.style={inner sep=1.7pt, outer sep=0pt, circle, fill}} 
\tikzset{labeled_vertex/.style={inner sep=2.2pt, outer sep=0pt, rectangle, fill=gray, draw=black}} 
\tikzset{edge_color0/.style={color=black,line width=1.2pt,opacity=0.5,dashed}} 
\tikzset{edge_color1/.style={color=red,  line width=1.2pt,opacity=0}} 
\tikzset{edge_color2/.style={color=blue, line width=1.2pt,opacity=1}} 
\tikzset{edge_color3/.style={color=green!80!black,line width=1.2pt}} 
\tikzset{edge_color4/.style={color=orange, line width=1.2pt}} 
\tikzset{edge_color5/.style={color=red,  line width=1.2pt,dotted}} 
\tikzset{edge_color6/.style={color=blue, line width=1.2pt,dotted}} 
\tikzset{edge_color7/.style={color=green, line width=1.2pt,dotted}} 
\tikzset{edge_color8/.style={color=gray, line width=1.2pt}} 
\tikzset{edge_colorroot/.style={color=red, line width=1.7pt}} 
\tikzset{edge_thin/.style={color=black}} 
\tikzset{edge_hidden/.style={color=black,dotted,opacity=0}} 
\tikzset{vertex_color1/.style={inner sep=1.7pt, outer sep=0pt, draw, circle, fill=red}} 
\tikzset{vertex_color2/.style={inner sep=1.7pt, outer sep=0pt, draw, circle, fill=blue}} 
\tikzset{vertex_color3/.style={inner sep=1.7pt, outer sep=0pt, draw, circle, fill=green}} 
\tikzset{vertex_color4/.style={inner sep=1.7pt, outer sep=0pt, draw, circle, fill=pink}} 
\tikzset{vertex_color5/.style={inner sep=1.7pt, outer sep=0pt, draw, circle, fill=pink,label=below:{$5$}}} 
\tikzset{vertex_color6/.style={inner sep=1.7pt, outer sep=0pt, draw, circle, fill=pink,label=below:{$6$}}} 
\tikzset{vertex_color7/.style={inner sep=1.7pt, outer sep=0pt, draw, circle, fill=pink,label=below:{$7$}}} 
\tikzset{vertex_color8/.style={inner sep=1.7pt, outer sep=0pt, draw, circle, fill=pink,label=below:{$8$}}} 
\tikzset{vertex_color9/.style={inner sep=1.7pt, outer sep=0pt, draw, circle, fill=pink,label=below:{$9$}}} 
\tikzset{vertex_color10/.style={inner sep=1.7pt, outer sep=0pt, draw, circle, fill=pink,label=below:{$10$}}} 
\tikzset{vertex_color11/.style={inner sep=1.7pt, outer sep=0pt, draw, circle, fill=pink,label=below:{$11$}}} 
\tikzset{vertex_color12/.style={inner sep=1.7pt, outer sep=0pt, draw, circle, fill=pink,label=below:{$12$}}} 
\tikzset{vertex_color13/.style={inner sep=1.7pt, outer sep=0pt, draw, circle, fill=pink,label=below:{$13$}}} 
\tikzset{vertex_color14/.style={inner sep=1.7pt, outer sep=0pt, draw, circle, fill=pink,label=below:{$14$}}}\tikzset{labeled_vertex_color1/.style={inner sep=2.2pt, outer sep=0pt, draw, rectangle, fill=red}} 
\tikzset{labeled_vertex_color2/.style={inner sep=2.2pt, outer sep=0pt, draw, rectangle, fill=blue}} 
\tikzset{labeled_vertex_color3/.style={inner sep=2.2pt, outer sep=0pt, draw, rectangle, fill=green}} 
\tikzset{labeled_vertex_color4/.style={inner sep=2.2pt, outer sep=0pt, draw, rectangle, fill=pink}} 

\def\outercycle#1#2{ 
\pgfmathtruncatemacro{\plusone}{#1+1} 
\pgfmathtruncatemacro{\minusone}{#1-1} 
\draw  \foreach \x in {0,...,\minusone}{(0.5*\x,0) coordinate(x\x)};} 

\def\drawhypervertex#1#2{ \pgfmathtruncatemacro{\plusone}{#1+1}  \draw[edge_color2] (x#1)++(0,-0.2-0.2*#2)+(-0.2,0) -- +(0.2,0) +(0,0) node[fill=white,outer sep=0,inner sep=0]{\tiny \plusone};}

\def\drawhyperedge#1#2{\pgfmathtruncatemacro{\plusone}{#1+1} \draw[dotted] (x0)++(0,-0.2-0.2*#1)--++(0.5*#2-0.5,0);} 

\def\labelvertex#1{\pgfmathtruncatemacro{\vertexlabel}{#1+1 } \draw (x#1) node{\color{yellow}\tiny\vertexlabel}; }

\tikzset{
vtx/.style={inner sep=1.1pt, outer sep=0pt, circle, fill,draw}, 
vtxl/.style={inner sep=1.1pt, outer sep=0pt, rectangle, fill=yellow,draw=black}, 
hyperedge/.style={fill=pink,opacity=0.5,draw=black}, 
}

\definecolor{darkpastelgreen}{rgb}{0.01, 0.75, 0.24}

\title{Positive co-degree densities and jumps}
\author{J\'ozsef Balogh \footnote{Department of Mathematics, University of Illinois Urbana-Champaign, IL, USA
E-mail: \texttt{jobal@illinois.edu}.
Research supported in part by NSF grants  RTG DMS-1937241 and FRG DMS-2152488, the Arnold O. Beckman Research Award (UIUC Campus Research Board RB 24012).} \and 
Anastasia Halfpap\footnote{Department of Mathematics, Iowa State University, Ames, IA, USA.  E-mail: \texttt{ahalfpap@iastate.edu}. Research supported by NSF DSM-2152490.} 
\and 
Bernard Lidick\'y \footnote{Department of Mathematics, Iowa State University, Ames, IA, USA.  E-mail: \texttt{lidicky@iastate.edu}. Research supported by NSF DSM-2152490 and Scott Hanna Professorship.} 
\and
Cory Palmer \footnote{Department of Mathematical Sciences, University of Montana. Email: \texttt{cory.palmer@umontana.edu}. Research supported by a grant from the Simons Foundation \#712036.}
}
\date{}

\newcommand{\oururl}{\url{https://lidicky.name/pub/pco/}}

\begin{document}

\maketitle

\begin{abstract}
The \textit{minimum positive co-degree}  of a nonempty $r$-graph $H$, denoted by $\delta_{r-1}^+(H)$, is the largest integer $k$ such that for every $(r-1)$-set $S \subset V(H)$, if $S$  is contained in a hyperedge of $H$, then $S$ is contained in at least $k$ hyperedges of $H$. Given a family $\mathcal{F}$ of $r$-graphs, the \textit{positive co-degree Tur\'an function} $\mathrm{co^+ex}(n,\mathcal{F})$
is the maximum of $\delta_{r-1}^+(H)$ over all $n$-vertex $r$-graphs $H$ containing no member of $\mathcal{F}$. The \textit{positive co-degree density} of $\mathcal{F}$ is $\gamma^+(\mathcal{F}) = \underset{n \rightarrow \infty}{\lim} \frac{\mathrm{co^+ex}(n,\mathcal{F})}{n}.$ While the existence of $\gamma^+(\mathcal{F})$ is proved for all families $\mathcal{F}$, only few positive co-degree densities are known exactly. 

For a fixed $r \geq 2$, we call $\alpha \in [0,1]$ an \textit{achievable value} if there exists a family of $r$-graphs $\mathcal{F}$ with $\gamma^+(\mathcal{F}) = \alpha$, and call $\alpha$ a \textit{jump} if for some $\delta > 0$, there is no family $\mathcal{F}$ with $\gamma^+(\mathcal{F}) \in (\alpha, \alpha + \delta)$. Halfpap, Lemons, and Palmer~\cite{halfpap2024positive} showed that every $\alpha \in [0, \frac{1}{r})$ is a jump. We extend this result by showing that every $\alpha \in [0, \frac{2}{2r -1})$ is a jump. We also show that for $r = 3$, the set of achievable values is infinite, more precisely,  $\frac{k-2}{2k-3}$ for every $k \geq 4$ is achievable. Finally, we determine two additional achievable values for $r=3$ using flag algebra calculations.

\end{abstract}


\section{Introduction}

An $r$-\textit{graph} is a hypergraph in which all hyperedges have size $r$.  We often refer to the hyperedges of an $r$-graph as $r$-\textit{edges}. Given a family of $r$-graphs $\mathcal{F}$, the \textit{Tur\'an number} $\mathrm{ex}(n,\mathcal{F})$ is the maximum number of $r$-edges possible in an $n$-vertex $r$-graph that contains no member of $\mathcal{F}$ as a subhypergraph. When $r = 2$, the function $\mathrm{ex}(n,F)$ is well-studied and relatively well-understood. Given a set of $r$-graphs $\mathcal{F}$, we define the \textit{Tur\'an density} of $\mathcal{F}$ to be
\[ \pi(\mathcal{F}) := \underset{n \rightarrow \infty}{\lim} \frac{\mathrm{ex}(n,\mathcal{F})}{\binom{n}{r}}.\]
The Erd\H{o}s-Stone Theorem~\cite{ErSt} as pointed out by Erd\H{o}s-Simonovits~\cite{ErSi} determines the Tur\'an density of every $2$-graph as a function of its chromatic number.

\begin{theorem}[Erd\H{o}s-Stone]\label{ess}

Let $F$ be a $2$-graph with $\chi(F) = k$. Then
$\pi(F) = 1 - \frac{1}{k-1}. $

\end{theorem}
Note that an extension of Theorem~\ref{ess} also  claims that for every family $\mathcal{F}$ of $2$-graphs, $\pi(\mathcal{F})$ is equal to the minimum of $1-1/(\chi(F)-1)$ over $F \in \mathcal{F}$. While Theorem~\ref{ess} does not give us perfect information about Tur\'an numbers (in particular, for bipartite graphs $F$, it only demonstrates that $\mathrm{ex}(n,F) = o(n^2)$), it yields a good ``approximate'' understanding of Tur\'an numbers by fully describing Tur\'an densities.

For $r \geq 3$, we do not have an analogue to Theorem~\ref{ess}, and much less is known about Tur\'an densities.
Even $\pi(\mathcal{F})$ could be smaller than 
$\min\{\pi(F) : F \in \mathcal{F}\}$, as observed by Balogh~\cite{jozsituran}.
Not only do we lack a general theory, but the Tur\'an densities of many small $r$-graphs remain unknown, despite great effort. Famously, the Tur\'an density of the tetrahedron $K_4^3$ is still undetermined. The difficulty of determining Tur\'an densities for hypergraphs has motivated the study of various other hypergraph extremal functions, typically maximizing some variant of the minimum degree. In particular, given an $r$-graph $H$, the \textit{co-degree} of a set $S \in \binom{V(H)}{r-1}$ is the number of $r$-edges containing $S$, and the {\it minimum co-degree} of $H$, denoted by $\delta_{r-1}(H)$, is the smallest co-degree realized by an $(r-1)$-set contained in $V(H)$. The \textit{co-degree Tur\'an function} of a family of $r$-graphs $\mathcal{F}$, denoted by $\mathrm{coex}(n,\mathcal{F})$, is the largest possible minimum co-degree of an $n$-vertex $r$-graph containing no member of $\mathcal{F}$ as a subhypergraph.

Mubayi and Zhao \cite{Mubayicode} showed that the \textit{co-degree density}
\[\gamma(\mathcal{F}) := \underset{n \rightarrow \infty}{\lim} \frac{\mathrm{coex}(n,\mathcal{F})}{n}\]
exists for every family of $r$-graphs $\mathcal{F}$, and studied the general behavior of $\mathrm{coex}(n,\mathcal{F})$. Note that for every family of $2$-graphs $\mathcal{F}$ we have $\gamma(\mathcal{F}) = \pi(\mathcal{F})$; however, for $r \geq 3$, co-degree Tur\'an problems are not equivalent to Tur\'an problems, and co-degree density does not behave in the manner suggested by Theorem~\ref{ess}. We first define the notion of a \textit{jump} in density.

\begin{definition}
Fix $r \geq 2$. Suppose $\varphi$ is a function that maps families of $r$-graphs to $[0,1]$. We say that $\alpha \in [0, 1)$ is a $\varphi$-{\it jump} if there exists $\delta \in (0, 1 - \alpha)$ such that for no family $\mathcal{F}$ of $r$-graphs, $\varphi(\mathcal{F}) \in (\alpha, \alpha + \delta)$.
\end{definition}

While Theorem~\ref{ess} shows that co-degree density (and  Tur\'an density) jumps everywhere when $r = 2$, Mubayi and Zhao showed that co-degree density does not jump when $r \geq 3$.

\begin{theorem}[Mubayi-Zhao \cite{Mubayicode}]\label{mzjumps}

For $r \geq 3$, no $\alpha \in [0,1)$ is a $\gamma$-jump.

\end{theorem}

This co-degree phenomenon was further investigated in \cite{piga2023smallcodegree,ding2023vanishingcodegree}.

We remark that Theorem~\ref{mzjumps} suggests a substantial difference in behavior between Tur\'an and co-degree Tur\'an problems for hypergraphs.  We know for every $r \geq 3$ that every $\alpha \in [0, r!/r^r)$ is a $\pi$-jump. On the other hand, $\pi$ is also known to \textit{not} jump in infinitely many places for $r \geq 3$; see~\cite{hg-no-jump, do-jump}.

The \textit{minimum positive co-degree} of an $r$-graph $H$ is the largest integer $k$ such that, whenever $S \in \binom{V(H)}{r-1}$ is contained in some $r$-edge of $H$, then $S$ is contained in at least $k$ $r$-edges of $H$. 
The edgeless $r$-graph is defined to have positive co-degree zero. 
We denote the minimum positive co-degree of $H$ by $\delta_{r-1}^+(H)$. 
We define the \textit{positive co-degree Tur\'an number}, denoted by  $\mathrm{co^+ex}(n,\mathcal{F})$, to be the largest possible minimum positive co-degree of an $n$-vertex $r$-graph containing no member of $\mathcal{F}$ as a subhypergraph.

Balogh, Lemons, and Palmer~\cite{BLP} introduced the minimum positive co-degree as an alternative notion of minimum degree in $r$-graphs. Since then this parameter has already been studied from several angles. The concept of $\mathrm{co^+ex}(n,\mathcal{F})$ was recently introduced by Halfpap, Lemons, and Palmer~\cite{halfpap2024positive}. The investigation of minimum positive co-degree as an extremal parameter is partially motivated by the admissibility of constructions that mimic the extremal graphs for classical questions. For example, given an $r$-graph $H$, the \textit{t-blow-up} $H[t]$ of $H$ is the $r$-graph obtained by replacing each vertex $v_i \in V(H)$ with a class $V_i$ of $t$ vertices, where a set of $r$-vertices spans an $r$-edge if and only if they belong to $r$ distinct classes of $H[t]$ which correspond to an $r$-edge in $H$. In classical Tur\'an theory, graph blow-ups yield extremal or nearly extremal constructions for all non-bipartite forbidden graphs. Blow-ups also occur as extremal examples for other types of thresholds---for instance, one of the constructions demonstrating the tightness of Dirac's Theorem is a slightly unbalanced blow-up of an edge (i.e., a complete bipartite graph). 

For $r \geq 3$ and every $r$-graph $H$ a sufficiently large blow-up of $H$ has minimum co-degree $0$, which means that even after adding $o(n^3)$ hyperedges,  blow-ups will not provide extremal constructions for minimum co-degree density problems. However, $H[t]$ inherits the \textit{positive} co-degree properties of $H$. Thus, blow-ups (as well as other constructions with co-degree $0$ sets, such as $r$-graphs containing large strongly independent sets or multiple components) are potential extremal examples for positive co-degree problems. Previous results suggest that extremal constructions for positive co-degree problems in fact \textit{do} often look analogous to classical extremal constructions. See \cite{halfpap2024positive} for extremal constructions avoiding some small $3$-graphs, and \cite{halfpapmagnanspanning} on positive co-degree analogs of Dirac's Theorem, for which the extremal constructions also naturally generalize the graph extremal examples (and have minimum co-degree 0). Due to the expanded range of potential constructions, $\mathrm{co^+ex}(n,\mathcal{F})$ and $\mathrm{coex}(n,\mathcal{F})$ are generally not equal, and they appear to behave differently. 

Define the \textit{positive co-degree density} of a family of $r$-graphs $\mathcal{F}$ as the limit
\[\gamma^+(\mathcal{F}) := \underset{n \rightarrow \infty}{\lim} \frac{\mathrm{co^+ex}(n,\mathcal{F})}{n}. \]
The existence of $\gamma^+(F)$ was established by Halfpap, Lemons, and Palmer~\cite{halfpap2024positive} via a constructive argument, which can be generalized to finite families $\mathcal{F}$.
Pikhurko~\cite{pikhurko2023limitpositiveelldegreeturan} gave a probabilistic argument establishing that $\gamma^+(\mathcal{F})$ exists for all families $\mathcal{F}$. 

\begin{prop}[Halfpap-Lemons-Palmer~\cite{halfpap2024positive}]\label{first jump}

Fix $r \geq 2$ and let $\mathcal{F}$ be a family of $r$-graphs. Then
\[\gamma^+(\mathcal{F}) \in \{0\} \cup \left[\frac{1}{r}, 1\right].\]
\end{prop}

In other words, every $\alpha \in [0,1/r)$ is a $\gamma^+$-jump. 
Proposition~\ref{first jump} describes behavior similar to that of the classical Tur\'an density. Every $r$-graph $F$ that is contained in some blow-up of an $r$-edge can be shown to be ``degenerate'', having $\mathrm{co^+ex}(n,F) = o(n)$, so $\gamma^+(F) = 0$. On the other hand, if $F$ is not contained in any blow-up of an $r$-edge, then $\gamma^+(F) \geq \frac{1}{r}$, since the balanced $n$-vertex blow-up of an $r$-edge has minimum positive co-degree approximately $\frac{n}{r}$. 
An $r$-graph $F$ is \emph{$k$-partite} if there is a partition of $V(F)$ into $k$ classes such that each edge  intersects each part at most once.

Although Proposition~\ref{first jump} suggests that $\gamma$ and $\gamma^+$ exhibit fundamentally different beha\-viors, it is not clear what behavior to expect from $\gamma^+$. Currently, we know the exact value $\gamma^+(F)$  only for very few $3$-graphs $F$. Halfpap, Lemons, and Palmer~\cite{halfpap2024positive} determined the values of $\gamma^+(F)$ for many small $3$-graphs $F$, and bounded $\gamma^+(F)$ in some other instances. 
Various authors have reported improvements on several of these initial bounds, with the current best known values summarized in Table~\ref{previous bounds}. For comparison, the best-known bounds on $\pi$ and $\gamma$ for these $3$-graphs are also provided.
Graphs not defined in Table~\ref{fig-3graphs} are
defined by their edge sets as
\begin{align*}
K_4^{3-} &= \{123, 124, 134\}, \quad F_5 = \{123, 124, 345\},\quad
\mathbb{F} = \{ 123, 345, 156, 246, 147, 257, 367 \},\\
C_\ell &= \{123, 234, 345,\ldots, (\ell-2)(\ell-1)\ell, (\ell-1)\ell 1, \ell 12\},
\quad
C_\ell^- = C_\ell - \{\ell 12\}.
\end{align*}
See~\cite{halfpap2024positive,l2norm} for more details about $3$-graphs in Table~\ref{previous bounds}. Note that the $3$-graph denoted here as $K_4^{3-}$ is often called $K_4^{-}$ in the literature; we adopt the notation $K_4^{3-}$ to distinguish this $3$-graph from the 2-graph obtained by deleting an edge from $K_4$, which we denote by $K_4^-$.

\begin{table}
\begin{center}
\begin{tabular}{|l||c|c|c|c|c|c|}
\hline
 $F$ & $\leq \pi(F) $  & $\pi(F) \leq $ & $\leq \gamma(F) $ & $\gamma(F) \leq   $ & $\leq \gamma^+(F)  $ & $\gamma^+(F) \leq $ \\
 \hline
 \hline
{$K_4^{3-}$} & 2/7 \cite{K43-extremalfrankl}  & 0.28689 \cite{flagmatic} & $1/4$ \cite{codegreeconj} & $1/4$  \cite{FalgasK4-}   & $ 1/3 $ \cite{halfpap2024positive}& $  1/3 $ \cite{halfpap2024positive}\\
{$F_5$} & 2/9  \cite{Bollobascancellative} & 2/9 \cite{F5Frankl} & 0 \cite{l2norm} & 0 \cite{l2norm} & $1/3$ \cite{halfpap2024positive}& $1/3$ \cite{halfpap2024positive}\\
{$F_{3,2}$} & 4/9  \cite{F32Mubayi} & 4/9 \cite{F32furedi} & 1/3 \cite{codF32falgas} & 1/3 \cite{codF32falgas} & $1/2$ \cite{halfpap2024positive} & $1/2$ \cite{halfpap2024positive}\\
{$\mathbb{F}$} & 3/4 \cite{VeraSos} & 3/4 \cite{FanoFuredi} & 1/2 \cite{MubayiFano} & 1/2 \cite{MubayiFano} & 2/3 \cite{halfpap2024positive} & 2/3 \cite{halfpap2024positive}\\
{$K_4^3$} &5/9 \cite{MR177847} & 0.5615 \cite{BaberTuran} & 1/2 \cite{MR1829685} & 0.529 \cite{l2norm} & $1/2$ \cite{halfpap2024positive} &  0.543  \cite{volec} \\ 

{$F_{3,3}$} & 3/4 \cite{F32Mubayi} & 3/4 \cite{F32Mubayi} & 1/2 \cite{l2norm} & 0.604 \cite{l2norm} & 3/5 \cite{halfpap2024positive}& 0.616  \\
{$C_{5}$} & $2\sqrt{3} - 3$ \cite{F32Mubayi} & 0.46829 \cite{flagmatic} & 1/3 \cite{l2norm} & 0.3993 \cite{l2norm} & $1/2$ \cite{halfpap2024positive} & 1/2 \cite{Wu} \\
{$C_{7}$} &$2\sqrt{3} - 3$ \cite{F32Mubayi} &
0.464186 
& 1/3 \cite{l2norm} & 0.371 & 1/2 \cite{halfpap2025} & 1/2 \cite{halfpap2025} \\
{$C_{5}^-$} &1/4 \cite{F32Mubayi} & 1/4\cite{lidicky2024c5-} & 0 \cite{l2norm} & 0 \cite{piga} & 1/3 \cite{halfpap2024positive} & 1/3 \cite{Wu} \\

{$J_4$} & $1/2$ \cite{DaisyBollobas} & $0.50409$ \cite{flagmatic} & $1/4$ \cite{l2norm} & $0.473$ \cite{l2norm} & $4/7$ \cite{halfpap2024positive} & 4/7  \\

{$F_{4,2}$} & $4/9$ & $0.4933328$  & $1/3$ & $0.4185 $ & $3/5$  & $3/5$  \\
\hline
\end{tabular}
\captionof{table}{Best-known density bounds for $\pi, \gamma$, and $\gamma^+$.}\label{previous bounds}
\end{center}
\end{table}

Kamčev, Letzter, and Pokrovskiy~\cite{Kamcev2024} proved that the Turán density of longer tight cycles $C_\ell$ is $2\sqrt{3} - 3$, when $\ell$ is not multiple of three and sufficiently large (when $\ell $ is divisible by three, then $C_\ell$ is $3$-partite, hence its Tur\'an density is $0$).
Similarly, Balogh and Luo~\cite{balogh2024turandensitylongtight}
proved for $\ell$ sufficiently large and not divisible by $3$ that the Tur\'an density of $C_\ell^-$ is $1/4$.
Recently, this was proved for every $\ell\ge 5$ by Lidick\'y, Mattes, and Pfender~\cite{lidicky2024c5-}.
That the co-degree density is $1/3$ for tight cycles of length at least $10$ and not divisible by $3$ was proved by Piga, Sanhueza-Matamala, and Schacht~\cite{piga2024codegreeturandensity3uniform} and Ma~\cite{ma2024codegreeturandensity3uniform}.
On the other hand, the fact that the co-degree density of $C_\ell^-$ is $0$ is due to Piga, Sales, and Sch\"ulke~\cite{piga2023smallcodegree}.

Given $r \geq 2$, we call $\alpha \in [0,1]$ an \textit{achievable value} (for $\gamma^+$) if there exists a family $\mathcal{F}$ of $r$-graphs such that $\gamma^+(\mathcal{F}) = \alpha$. 
Before our article, the only known achievable values of $\gamma^+$ for $r=3$ were $0,1/3, 1/2,$ and $2/3$. 

Our goal is to understand the positive co-degree density by demonstrating additional $\gamma^+$-jumps for every $r$, as well as by expanding the known list of achievable values of $\gamma^+$ for $r =3$.
Our paper is organized as follows. In Section~\ref{new results}, we summarize our main results. In Section~\ref{prelims}, we state some additional definitions and lemmas which will be used in our proofs. In Section~\ref{proofs}, we demonstrate further $\gamma^+$-jumps for every $r$, and for $r =3$ establish an infinite set of achievable values for $\gamma^+$. 
In Section~\ref{J4}, we use flag algebras to exactly determine $\gamma^+$ for two additional $3$-graphs, hence adding two values to the list of achievable values of $\gamma^+$. Finally, in Section~\ref{conclusion}, we have some concluding remarks and list a wide variety of open problems in the area.

\section{New results}\label{new results}

Our first main theorem extends the range of jumps described in Proposition~\ref{first jump}. We will need the following definition.
An \emph{$r$-triangle}, denoted by $T^r$, is an $r$-graph with $r+1$ vertices and three $r$-edges. Notice that $T^r$ can be obtained from the $2$-graph triangle $T^2$ by adding $r-2$ vertices and including them to each of the three edges.
Such hypergraphs are sometimes called daises in the literature.
The $T^2$ part of a $T^r$ is called the \emph{base} of the $T^r$.
Notice that $T^3$ is the same as $K_4^{3-}$.

\begin{restatable}{theorem}{secondjump}\label{general jumps}
Let $\mathcal{F}$ be a family of $r$-graphs for $r \geq 2$. Then 

\[ \gamma^+(\mathcal{F}) \in \left\{0, \frac{1}{r}\right\} \cup \left[\frac{2}{2r-1}, 1\right]. \]
Thus, every $\alpha \in[0, \frac{2}{2r-1})$ is a $\gamma^+$-jump.

Moreover, $\gamma^+(\mathcal{F}) = 0$ if and only if some member of $\mathcal{F}$ is $r$-partite, and $\gamma^+(\mathcal{F}) = \frac{1}{r}$ if and only if no member of $\mathcal{F}$ is $r$-partite, but some member of $\mathcal{F}$ is contained in a blow-up of some $T^r$. 
\end{restatable}

For $r=3$, we also provide an infinite set of achievable values for $\gamma^+$, based on forbidden families involving the following $3$-graphs. For $k \geq 3$, let $J_k$ be the $(k+1)$-vertex $3$-graph, on vertex set $[k+1]$, with $3$-edges of the form $1ij$ for every $i , j \in \{2, \dots, k+1\}$. Let $K_4^3$ denote the complete $4$-vertex $3$-graph, and let $F_{3,2}$ denote the $5$-vertex $3$-graph, on vertex set $\{1,2,3,4,5\}$, with edge set $\{123, 124, 125, 345\}$. See Table~\ref{fig-3graphs} for an illustration of these and other relevant $3$-graphs.

\begin{figure}
\begin{center}
  \begin{tabular}{ | l|| c |  Sc | Sc  |   }
    \hline
    $K_4^3$ 
    & 
\vc{\begin{tikzpicture}
\outercycle{5}{4}
\drawhyperedge{0}{4}
\drawhypervertex{0}{0}
\drawhypervertex{1}{0}
\drawhypervertex{2}{0}
\drawhyperedge{1}{4}
\drawhypervertex{0}{1}
\drawhypervertex{1}{1}
\drawhypervertex{3}{1}
\drawhyperedge{2}{4}
\drawhypervertex{0}{2}
\drawhypervertex{2}{2}
\drawhypervertex{3}{2}
\drawhyperedge{3}{4}
\drawhypervertex{1}{3}
\drawhypervertex{2}{3}
\drawhypervertex{3}{3}
\end{tikzpicture} 
}    
    &
 \vc{
\begin{tikzpicture}[scale=0.55]
\foreach \i in {1,2,3,4}{
\draw (90*\i-45:1.5) coordinate(\i);
}

\foreach \r in {0-45,90-45,180-45,270-45}{
\begin{scope}[rotate=\r]
\draw[hyperedge] 
(0:1.5) to[out=140,in=275,looseness=1.2] (90:1.5)  to[out=265,in=40,looseness=1.2] (180:1.5)  to[out=35,in=145,looseness=1.2] (0:1.5) 
;
\end{scope}
}

\draw 
(1) node[vtx,label=right:{\tiny 1}]{}
(2) node[vtx,label=left:{\tiny 2}]{}
(3) node[vtx,label=left:{\tiny 3}]{}
(4) node[vtx,label=right:{\tiny 4}]{}
;
\end{tikzpicture}
}   
   \\
   \hline
$F_{3,2}$ &
\vc{\begin{tikzpicture}\outercycle{6}{5}
\drawhyperedge{0}{5}
\drawhypervertex{0}{0}
\drawhypervertex{1}{0}
\drawhypervertex{2}{0}
\drawhyperedge{1}{5}
\drawhypervertex{0}{1}
\drawhypervertex{3}{1}
\drawhypervertex{4}{1}
\drawhyperedge{2}{5}
\drawhypervertex{1}{2}
\drawhypervertex{3}{2}
\drawhypervertex{4}{2}
\drawhyperedge{3}{5}
\drawhypervertex{2}{3}
\drawhypervertex{3}{3}
\drawhypervertex{4}{3}
\end{tikzpicture} 
}
&
\vc{
\begin{tikzpicture}
\draw
(0,-0.7) coordinate(1) node[vtx,label=left:{\tiny $1$}](a){}
(0,0) coordinate(2) node[vtx,label=below:{\tiny $2$}](b){}
(0,0.7) coordinate(3) node[vtx,label=left:{\tiny $3$}](c){}
(1,-0.5) coordinate(4) node[vtx,label=right:{\tiny $5$}](d){}
(1,0.5) coordinate(5) node[vtx,label=right:{\tiny $4$}](d){}
;
\draw[hyperedge] (1) to[bend left] (2) to[bend left] (3) to[bend right] (1);
\draw[hyperedge] (1) to[out=0,in=190] (4) to[out=190,in=250] (5) to[out=250,in=0] (1);
\draw[hyperedge] (2) to[out=0,in=140] (4) to[out=140,in=210,looseness=1.5] (5) to[out=210,in=0] (2);
\draw[hyperedge] (3) to[out=0,in=120,looseness=1.2] (4) to[out=120,in=170,looseness=1.2] (5) to[out=170,in=0,looseness=1.2] (3);
\end{tikzpicture}
}
\\
\hline
$F_{4,2}$ &
\vc{ 
  \begin{tikzpicture}[flag_pic]
\drawhyperedge{0}{6}
\drawhypervertex{0}{0}
\drawhypervertex{1}{0}
\drawhypervertex{2}{0}
\drawhyperedge{1}{6}
\drawhypervertex{0}{1}
\drawhypervertex{1}{1}
\drawhypervertex{3}{1}
\drawhyperedge{2}{6}
\drawhypervertex{0}{2}
\drawhypervertex{2}{2}
\drawhypervertex{3}{2}
\drawhyperedge{3}{6}
\drawhypervertex{0}{3}
\drawhypervertex{4}{3}
\drawhypervertex{5}{3}
\drawhyperedge{4}{6}
\drawhypervertex{1}{4}
\drawhypervertex{4}{4}
\drawhypervertex{5}{4}
\drawhyperedge{5}{6}
\drawhypervertex{2}{5}
\drawhypervertex{4}{5}
\drawhypervertex{5}{5}
\drawhyperedge{6}{6}
\drawhypervertex{3}{6}
\drawhypervertex{4}{6}
\drawhypervertex{5}{6}
\end{tikzpicture} } 
&
\vc{
\begin{tikzpicture}
\draw
(0,-0.7) coordinate(1) node[vtx,label=below:{\tiny $1$}](a){}
(0,-0.2) coordinate(2) node[vtx,label=below:{\tiny $2$}](b){}
(0, 0.2) coordinate(6) node[vtx,label=above:{\tiny $3$}](b){}
(0,0.7) coordinate(3) node[vtx,label=left:{\tiny $4$}](c){}
(1,-0.5) coordinate(4) node[vtx,label=right:{\tiny $6$}](d){}
(1,0.5) coordinate(5) node[vtx,label=right:{\tiny $5$}](d){}
;
\draw[hyperedge] (1) to[bend left] (2) to[bend left] (3) to[bend right] (1);
\draw[hyperedge] (1) to[bend left=40] (6) to[bend left=40] (3) to[bend right=40] (1);
\draw[hyperedge] (1) to[bend left=70,looseness=2] (2) to[bend left=70,looseness=2] (6) to[bend right=70,looseness=2] (1);
\draw[hyperedge] (1) to[out=0,in=190] (4) to[out=190,in=250] (5) to[out=250,in=0] (1);
\draw[hyperedge] (2) to[out=0,in=140] (4) to[out=140,in=210,looseness=1.5] (5) to[out=210,in=0] (2);
\draw[hyperedge] (6) to[out=0,in=140] (4) to[out=140,in=210,looseness=1.5] (5) to[out=210,in=0] (6);
\draw[hyperedge] (3) to[out=0,in=120,looseness=1.2] (4) to[out=120,in=170,looseness=1.2] (5) to[out=170,in=0,looseness=1.2] (3);
\end{tikzpicture}
} 
\\\hline
$F_{3,3}$
 &
\vc{\begin{tikzpicture}\outercycle{7}{6}
\drawhyperedge{0}{6}
\drawhypervertex{0}{0}
\drawhypervertex{1}{0}
\drawhypervertex{2}{0}
\drawhyperedge{1}{6}
\drawhypervertex{0}{1}
\drawhypervertex{3}{1}
\drawhypervertex{4}{1}
\drawhyperedge{2}{6}
\drawhypervertex{0}{2}
\drawhypervertex{3}{2}
\drawhypervertex{5}{2}
\drawhyperedge{3}{6}
\drawhypervertex{0}{3}
\drawhypervertex{4}{3}
\drawhypervertex{5}{3}
\drawhyperedge{4}{6}
\drawhypervertex{1}{4}
\drawhypervertex{3}{4}
\drawhypervertex{4}{4}
\drawhyperedge{5}{6}
\drawhypervertex{1}{5}
\drawhypervertex{3}{5}
\drawhypervertex{5}{5}
\drawhyperedge{6}{6}
\drawhypervertex{1}{6}
\drawhypervertex{4}{6}
\drawhypervertex{5}{6}
\drawhyperedge{7}{6}
\drawhypervertex{2}{7}
\drawhypervertex{3}{7}
\drawhypervertex{4}{7}
\drawhyperedge{8}{6}
\drawhypervertex{2}{8}
\drawhypervertex{3}{8}
\drawhypervertex{5}{8}
\drawhyperedge{9}{6}
\drawhypervertex{2}{9}
\drawhypervertex{4}{9}
\drawhypervertex{5}{9}
\end{tikzpicture} 
}
 &
\vc{
\begin{tikzpicture}
\draw
(0,-0.7) coordinate(1) node[vtx,label=left:{\tiny $3$}](a){}
(0,0) coordinate(2) node[vtx,label=left:{\tiny $2$}](b){}
(0,0.7) coordinate(3) node[vtx,label=left:{\tiny $1$}](c){}
(1,-0.7) coordinate(4) node[vtx,label=right:{\tiny $6$}](d){}
(1,0.0) coordinate(5) node[vtx,label=right:{\tiny $5$}](d){}
(1,0.7) coordinate(6) node[vtx,label=right:{\tiny $4$}](d){}
;
\draw[hyperedge] (1) to[bend left] (2) to[bend left] (3) to[bend right] (1);
\draw[hyperedge] (1) to[out=0,in=190] (4) to[out=190,in=250] (5) to[out=250,in=0] (1);
\draw[hyperedge] (2) to[out=0,in=140] (4) to[out=140,in=210,looseness=1.5] (5) to[out=210,in=0] (2);
\draw[hyperedge] (3) to[out=-60,in=120,looseness=1.2] (4) to[out=120,in=170,looseness=1.2] (5) to[out=170,in=-60,looseness=1.2] (3);
\draw[hyperedge] (1) to[out=60,in=190] (5) to[out=190,in=250] (6) to[out=250,in=60] (1);
\draw[hyperedge] (2) to[out=0,in=140] (5) to[out=140,in=210,looseness=1.5] (6) to[out=210,in=0] (2);
\draw[hyperedge] (3) to[out=0,in=120,looseness=1.2] (5) to[out=120,in=170,looseness=1.2] (6) to[out=170,in=0,looseness=1.2] (3);
\draw[hyperedge] (1) to[out=10,in=170,looseness=1.5] (4) to[out=170,in=190] (6) to[out=190,in=10] (1);
\draw[hyperedge] (2) to[out=0,in=140] (4) to[out=140,in=210,looseness=1.5] (6) to[out=210,in=0] (2);
\draw[hyperedge] (3) to[out=-10,in=120,looseness=1.2] (4) to[out=120,in=190,looseness=1.2] (6) to[out=190,in=-10,looseness=1.2] (3);
\end{tikzpicture}
}
   \\
   \hline
   $J_4$
   &
\vc{\begin{tikzpicture}\outercycle{6}{5}
\drawhyperedge{0}{5}
\drawhypervertex{0}{0}
\drawhypervertex{1}{0}
\drawhypervertex{2}{0}
\drawhyperedge{1}{5}
\drawhypervertex{0}{1}
\drawhypervertex{1}{1}
\drawhypervertex{3}{1}
\drawhyperedge{2}{5}
\drawhypervertex{0}{2}
\drawhypervertex{1}{2}
\drawhypervertex{4}{2}
\drawhyperedge{3}{5}
\drawhypervertex{0}{3}
\drawhypervertex{2}{3}
\drawhypervertex{3}{3}
\drawhyperedge{4}{5}
\drawhypervertex{0}{4}
\drawhypervertex{2}{4}
\drawhypervertex{4}{4}
\drawhyperedge{5}{5}
\drawhypervertex{0}{5}
\drawhypervertex{3}{5}
\drawhypervertex{4}{5}
\end{tikzpicture} 
}   
   &
\vc{
\begin{tikzpicture}[scale=0.75]
\clip (-2.5,-1.1) rectangle (1.5,1.1);
\draw
(-2,0) coordinate(1) node[vtx,label=left:{\tiny $1$}](b){}
(45:1) coordinate(2) node[vtx,label=right:{\tiny $2$}](c){}
(135:1) coordinate(3) node[vtx,label=left:{\tiny $3$}](d){}
(225:1) coordinate(4) node[vtx,label=left:{\tiny $4$}](d){}
(315:1) coordinate(5) node[vtx,label=right:{\tiny $5$}](d){}
;
\draw[hyperedge] (1) to[out=70,in=140,looseness=0.9] (2) to[out=140,in=60,looseness=0.8]
(3) to[out=120,in=70,looseness=0.8] (1);
\draw[hyperedge] (1) to[out=0,in=250,looseness=1.3] (3) to[out=250,in=110] 
(4) to[out=110,in=0,looseness=1.3] (1);
\draw[hyperedge] (1) to[out=-70,in=240,looseness=0.8] (4) to[out=300,in=220,looseness=0.8] 
(5) to[out=220,in=-70,looseness=0.9] (1);
\draw[hyperedge] (1) to[out=0,in=110,looseness=0.7] (5) to[out=110,in=250,looseness=1.1] 
(2) to[out=250,in=0,looseness=0.7] (1);
\draw[hyperedge] (1) to[out=-25,in=205,looseness=0.8] (2) to[out=205,in=110,looseness=0.8] 
(4) to[out=110,in=-25,looseness=0.9] (1);
\draw[hyperedge] (1) to[out=25,in=250,looseness=0.9] (3) to[out=250,in=155,looseness=0.8] 
(5) to[out=155,in=25,looseness=0.8] (1);
\draw[line width = 0.5pt]
(3) to[out=250,in=155,looseness=0.8] (5)
(2) to[out=205,in=110,looseness=0.8] (4)
(5) to[out=110,in=250,looseness=1.1] (2)
(4) to[out=300,in=220,looseness=0.8] (5)
(3) to[out=250,in=110] (4)
(2) to[out=140,in=60,looseness=0.8] (3)
;
\end{tikzpicture}
}   
   \\
   \hline
   $J_5$
   &
\vc{\begin{tikzpicture}\outercycle{7}{6}
\drawhyperedge{0}{6}
\drawhypervertex{0}{0}
\drawhypervertex{1}{0}
\drawhypervertex{2}{0}
\drawhyperedge{1}{6}
\drawhypervertex{0}{1}
\drawhypervertex{1}{1}
\drawhypervertex{3}{1}
\drawhyperedge{2}{6}
\drawhypervertex{0}{2}
\drawhypervertex{1}{2}
\drawhypervertex{4}{2}
\drawhyperedge{3}{6}
\drawhypervertex{0}{3}
\drawhypervertex{1}{3}
\drawhypervertex{5}{3}
\drawhyperedge{4}{6}
\drawhypervertex{0}{4}
\drawhypervertex{2}{4}
\drawhypervertex{3}{4}
\drawhyperedge{5}{6}
\drawhypervertex{0}{5}
\drawhypervertex{2}{5}
\drawhypervertex{4}{5}
\drawhyperedge{6}{6}
\drawhypervertex{0}{6}
\drawhypervertex{2}{6}
\drawhypervertex{5}{6}
\drawhyperedge{7}{6}
\drawhypervertex{0}{7}
\drawhypervertex{3}{7}
\drawhypervertex{4}{7}
\drawhyperedge{8}{6}
\drawhypervertex{0}{8}
\drawhypervertex{3}{8}
\drawhypervertex{5}{8}
\drawhyperedge{9}{6}
\drawhypervertex{0}{9}
\drawhypervertex{4}{9}
\drawhypervertex{5}{9}
\end{tikzpicture} 
}
&
\vc{
\begin{tikzpicture}
\draw
(-2,0) coordinate(1) node[vtx,label=left:{\tiny $1$}](1){}
\foreach \x in {2,...,6}
{
(-144+72*\x:1) coordinate(\x) node[vtx](v\x){}
}
;
\draw[hyperedge,opacity=0.2] (1)--(2)--(3)--(1);
\draw[hyperedge,opacity=0.2] (1)--(2)--(4)--(1);
\draw[hyperedge,opacity=0.2] (1)--(2)--(5)--(1);
\draw[hyperedge,opacity=0.2] (1)--(2)--(6)--(1);
\draw[hyperedge,opacity=0.2] (1)--(3)--(4)--(1);
\draw[hyperedge,opacity=0.2] (1)--(3)--(5)--(1);
\draw[hyperedge,opacity=0.2] (1)--(3)--(6)--(1);
\draw[hyperedge,opacity=0.2] (1)--(4)--(5)--(1);
\draw[hyperedge,opacity=0.2] (1)--(4)--(6)--(1);
\draw[hyperedge,opacity=0.2] (1)--(5)--(6)--(1);
\draw [line width = 0.5pt] (2)--(3)--(4)--(5)--(6)--(2)--(4)--(6)--(3)--(5)--(2)  ;
\draw (-2,0) coordinate(1) node[vtx](1){}
\foreach \x in {2,...,6}
{
(-144+72*\x:1) coordinate(\x) node[vtx](v\x){}
};
\end{tikzpicture}
}\\
\hline
\end{tabular}
\end{center}
\captionof{table}{Small 3-uniform hypergraphs.}\label{fig-3graphs}
\end{figure}

\begin{restatable}{theorem}{densities}\label{densities}

For every $k \geq 4$,
\[ \gamma^+(\{K_4^3, F_{3,2}, J_k\}) = \frac{k-2}{2k-3}.\]

\end{restatable}

We also investigate whether the densities exhibited by Theorem~\ref{densities} can be achieved by forbidding a single $3$-graph. For $k = 4$, we have $\gamma^+(\{K_4^3, F_{3,2}, J_4\}) = \frac{2}{5}$.
Let $F_1$ be a 7-vertex 3-graph 
with edges $\{125, 135, 235, 126, 146, 246, 347 \}$,
see Figure~\ref{F1F2}.
We show that $F_1$ has positive co-degree density $\frac{2}{5}$.

\begin{restatable}{theorem}{singlegraph}\label{single graph}

$\gamma^+(F_1) = \frac{2}{5}$.

\end{restatable}

Finally, we utilize flag algebras to determine the positive co-degree densities of two graphs, which exhibit new achievable values of $\gamma^+$ larger than  $\frac{1}{2}$. We determine $\gamma^+(J_4)$, and show that the (asymptotic) extremal construction is the blow-up of the complement of the Fano plane. 
Notice that $J_4$ is a 3-daisy; see~\cite{ellis2024daisies} for recent progress on Tur\'an densities of $r$-daisies.
\begin{restatable}{theorem}{singlegraphJfour}\label{J4}
$\gamma^+(J_4) = \frac{4}{7}$.
\end{restatable}

We introduce another $3$-graph, which will have a different new density. 
Let $F_{4,2}$ be the $6$-vertex $3$-graph with edges $\{123, 124, 134, 156, 256, 356, 456 \}$, depicted in Table~\ref{fig-3graphs}.
Note that in $F_{4,2}$ the common neighborhood of 5 and 6 is $\{1,2,3,4\}$, and $\{1,2,3,4\}$ spans a $K_4^{3-}$. We determine $\gamma^+(F_{4,2})$, and show that its (asymptotic) extremal construction is the balanced blow-up of $K_5^3$.

\begin{restatable}{theorem}{Fdensity}\label{F42density}
$\gamma^+(F_{4,2}) = \frac{3}{5}$.
\end{restatable}

We also include some non-tight results obtained using flag algebras. 
\begin{restatable}{theorem}{FlagDensity}\label{thm:FlagDensity}
The following bounds hold.
\[
\pi(F_{4,2}) \leq 0.4933327,
\quad
\gamma(F_{4,2}) \leq 0.4185,
\quad
\gamma^+(F_{3,3}) \leq 0.616.
\]
\end{restatable}

The lower bounds in Table~\ref{previous bounds} for $\pi(F_{4,2})$ and $\gamma(F_{4,2})$ come from the lower bound constructions for $F_{3,2}$.

\section{Preliminaries}\label{prelims}

We begin by stating some results related to supersaturation and a hypergraph removal lemma. The hypergraph removal lemma states that an $r$-graph containing only few copies of some subhypergraph $F$ can be made $F$-free by the deletion of only few $r$-edges. For a discussion of removal lemmas, including the below formulation, see~\cite{ConlonFox}.

\begin{lemma}\label{hypergraph removal} 
 Fix $\alpha > 0$ and let $F$ be an $r$-graph. There exists $\delta > 0$ such that if $H$ is an $n$-vertex $r$-graph containing at most $\delta n^{|V(F)|}$ copies of $F$, then there exists $E' \subset E(H)$ such that $|E'| \leq \alpha n^r$ and $H - E'$ is $F$-free.

\end{lemma}

Although an $r$-graph with $o(n^{|V(F)|})$ copies of $F$ can be made $F$-free by deleting $o(n^r)$ $r$-edges, it is not obvious that the deletion would change the minimum positive co-degree by only $o(n)$. The following ``clean-up'' lemma due to Halfpap, Lemons, and Palmer~\cite{halfpap2024positive} allows us to apply the hypergraph removal lemma to minimum positive co-degree problems. Roughly, this lemma guarantees that any positive co-degree drop arising from the deletion of a small set of $r$-edges can be mitigated by the deletion of another set of $r$-edges.

\begin{lemma}[Halfpap-Lemons-Palmer~\cite{halfpap2024positive}]\label{positive codegree cleanup} 

Let $H$ be an $n$-vertex $r$-graph and fix $0< \varepsilon <1$ small enough that $(r+1)! \varepsilon^{1/{2^{r-1}}}n^r < |E(H)|$. Let $H_1$ be a subhypergraph of $H$ obtained by the deletion of at most $\varepsilon n^r$ $r$-edges. Then $H_1$ has a subhypergraph $H_2$ 
with $\delta_{r-1}^+(H_2) \geq \delta_{r-1}^+(H) - 2^r r! \varepsilon^{1/{2^{r-1}}}n$.

\end{lemma}

In practice it is not difficult to fulfill
the condition in Lemma~\ref{positive codegree cleanup} that $(r+1)! \varepsilon^{1/{2^{r-1}}}n^r < |E(H)|$,  since $\delta_{r-1}^+(H)$ can be used to give a lower bound on $|E(H)|$.

\begin{lemma}[Halfpap-Lemons-Palmer~\cite{halfpap2024positive}]\label{edge approx}
Fix $c>0$ and
suppose ${H}$ is an $r$-graph with $\delta_{r-1}^+(H) \geq cn$. Then, for $n$ large enough, $|E( {H})| \geq \frac{1}{2}\frac{c^r}{r!} n^r$.
\end{lemma}

Thus, for an $n$-vertex $r$-graph $H$ with $n$ sufficiently large and $\delta_{r-1}^+(H) \geq cn$, we can choose $\varepsilon$ in Lemma~\ref{positive codegree cleanup} as a function of $c$ and $r$  alone.

Lemmas~\ref{hypergraph removal} and \ref{positive codegree cleanup} can be used to prove supersaturation and related properties for minimum positive co-degree problems. In particular, we have the following basic formulation.

\begin{theorem}[Halfpap-Lemons-Palmer~\cite{halfpap2024positive}]\label{supersaturation}

Fix $\varepsilon > 0$ and let $F$ be an $r$-graph. Then there exists $\delta > 0$ such that, if $H$ is an $n$-vertex $r$-graph with 
\[\delta_{r-1}^+(H) >\mathrm{co^+ex}(F) + \varepsilon n,\]
then $H$ contains at least $\delta n^{|V(F)|}$ copies of $F$. 
\end{theorem}

By a standard argument (see, e.g., \cite{Keevashsurvey}), if $\delta > 0$ and $t \in \mathbb{N}$ are fixed and $n$ is sufficiently large, then every $n$-vertex $r$-graph $H$ containing $\delta n^{|V(F)|}$ copies of $F$ must contain \acopyof $F[t]$. 
Thus, as an immediate consequence of Theorem~\ref{supersaturation}, we have blow-up invariance for $\mathrm{co^+ex}(n,F)$.

\begin{cor}[Halfpap-Lemons-Palmer~\cite{halfpap2024positive}]\label{blow up invariance}

Let $F$ be an $r$-graph and $t$ a positive integer. Then
$$\mathrm{co^+ex}(n,F) \leq \mathrm{co^+ex}(n,F[t]) \leq \mathrm{co^+ex}(n,F) + o(n).$$
\end{cor}

We remark that Corollary~\ref{blow up invariance} is useful because it implies that if $H$ and $F$ are $r$-graphs such that $F$ is contained in $H[t]$ for some $t$, then $\gamma^+(F) \leq \gamma^+(H)$. For example, when paired with the facts that a single $3$-edge $e$ has $\gamma^+(e) = 0$ and  $\gamma^+(C_5) = \frac{1}{2}$, and with an appropriate lower bound construction, Corollary~\ref{blow up invariance} implies that $\gamma^+(C_{\ell}) = 0$ if $\ell \equiv 0 \,\,(\textrm{mod } 3)$ and $\gamma^+(C_{\ell}) = \frac{1}{2}$ if $\ell \not\equiv 0 \,\,(\textrm{mod } 3)$ and $C_{\ell}$ is contained in a blow-up of $C_5$. In fact, this resolves the positive co-degree density of every tight cycle except for $C_4 = K_4^3$ and $C_7$.
Halfpap~\cite{halfpap2025} proved $\gamma^+(C_7)=1/2$. 

In Section~\ref{proofs}, we will apply essentially the same idea, finding one construction whose blow-up contains another and then relating their positive co-degree densities. However, for our purposes, a somewhat different formulation from the above statements will be desirable. The proof ideas for Theorem~\ref{supersaturation} and Corollary~\ref{blow up invariance} are used to derive the following lemma.

\begin{lemma}\label{family removal}

Let $F$ be a fixed $r$-graph on $f$ vertices and $\mathcal{F} = \{F_1, F_2, \dots ,F_k\}$ a finite family of $r$-graphs such that $F_i[f]$ contains $F$ for every $i \in [k]$. For every $d \in [0,1)$, if $\gamma^+(F) > d$, then for some $\beta > 0$ and $n$ sufficiently large there exists an  $n$-vertex, $\{F\} \cup \mathcal{F}$-free $r$-graph $H$ with $\delta_{r-1}^+(H) > \left( d + \beta \right)n$.
\end{lemma}

\begin{proof}

Given $d$ as stated, choose $\beta > 0$ so that $\gamma^+(F) \geq d + 3\beta$, and $\alpha > 0$ such that $2^rr!\alpha^{1/{2^{r-1}}} < \beta$ and $ (r+1)!\alpha^{1/{2^{r-1}}} n^r < \frac{1}{2}\frac{(2\beta)^r}{r!} n^r$.
For each $i \in [k]$, take $\delta_i > 0$ as guaranteed by Lemma~\ref{hypergraph removal} such that any hypergraph $H$ containing fewer than $\delta_i n^{|V(F_i)|}$ copies of $F_i$ can be made $F_i$-free by the deletion of at most $\frac{\alpha}{k} |E(H)|$ $r$-edges. Choose $N \in \mathbb{N}$ sufficiently large such that for all $n \geq N$, we have:  
\begin{itemize}[itemsep=1pt, parsep=0pt]
\item $\mathrm{co^+ex}(n,F) \geq (d + 2\beta)n$;
\item for each $i \in [k]$, if $H$ is an $n$-vertex graph containing $\delta_i n^{|V(F_i)|}$ copies of $F_i$, then $H$ contains \acopyof $F_i[f]$.
\end{itemize}

Fix $n \geq N$ and let $H$ be an $n$-vertex $r$-graph with $\delta_{r-1}^+(H) = \mathrm{co^+ex}(n,F)$. Then $H$ contains fewer than $\delta_i n^{|V(F_i)|}$ copies of $F_i$ for every $i \in [k]$, and thus can be made $\mathcal{F}$-free by deletion of at most $\alpha |E(H)|$ edges by repeated application of Lemma~\ref{hypergraph removal}. Since we have 
\[(r+1)!\alpha^{1/{2^{r-1}}}  n^r < \frac{1}{2}\frac{(2\beta)^r}{r!} n^r \leq |E(H)|\] 
by Lemma~\ref{edge approx} and the definition of $\alpha$, we can now apply Lemma~\ref{positive codegree cleanup} to delete an additional set of $r$-edges, resulting in an $\{F\} \cup \mathcal{F}$-free, $n$-vertex $r$-graph $H'$ with 
\[\delta_{r-1}^+(H') \geq \delta_{r-1}^+(H) - 2^rr!\alpha^{1/{2^{r-1}}}n >\delta_{r-1}^+(H) - \beta n \geq (d + \beta)n. \qedhere
\] 
\end{proof}

Lemma~\ref{family removal} is an important tool in finding positive co-degree densities as it essentially allows us to expand our list of forbidden configurations.

We conclude this section with some relevant definitions and notation. 
We often consider $4$-vertex cliques with one edge removed; these may be $2$-uniform or $3$-uniform. To avoid ambiguity, we denote by $K_4^-$ the $2$-graph on $4$ vertices and $5$ edges, and by $K_4^{3-}$ the $3$-graph on $4$ vertices and three $3$-edges.

Some of the hypergraphs we consider can be naturally described as arising from lower-uniformity hypergraphs. We define the following operation, which increases the uniformity of a hypergraph by one.
Given an $r$-graph $H$, the \textit{suspension} $\widehat{H}$ is the $(r+1)$-graph with vertex set consisting of $V(H)$ and one new vertex $v$, and $(r+1)$-edges
\[ 
E(\widehat{H}) = \{ e \cup \{v\} : e \in E(H) \}.
\]

We call $V(H)$ the \textit{$r$-graph vertices} and $v$ the \textit{spike vertex}.
Notice that the $(r+1)$-triangle $T^{r+1}$ is a suspension of the $r$-triangle $T^{r}$.

Let $H$ be an $r$-graph. For a subset $X=\{x_1,\ldots,x_{r-1}\}$ of size $r-1$ of vertices of $H$, denote by $N(X)$ the set of all vertices $v$ such that $X\cup\{v\} \in E(H)$. We use $d(X) := |N(X)|$, which is the co-degree of $X$. 
To simplify notation we use $N(x_1,\ldots,x_{r-1}) := N(X)$ and $d(x_1,\ldots,x_{r-1}) := d(X)$. We call $N(X)$ the \textit{neighborhood} of $X$.

In a $3$-graph $H$, the \textit{link graph} $L(v)$ of a vertex $v \in V(H)$ is the auxiliary $2$-graph on $V(H) - \{v\}$ where $xy$ is an edge if and only if $vxy$ is a $3$-edge of $H$.

For 3-graphs $G$ and $H$ the \emph{density} of $G$ in $H$, denoted by $d(G,H)$, is the number of subgraphs of $H$ isomorphic to $G$ divided by $\binom{|V(H)|}{|V(G)|}$. Notice that the density is always in $[0,1]$.

\section{Jumps and positive co-degree densities below $\frac{1}{2}$}\label{proofs}

\begin{proof}[Proof of Theorem~\ref{general jumps}]
Let $\mathcal{F}$ be a family of $r$-graphs. 
By Proposition~\ref{first jump}, $\gamma^+(\mathcal{F}) \in \{0\} \cup [\frac{1}{r}, 1]$, so it is sufficient to show that there is no family $\mathcal{F}$ with $\gamma^+(\mathcal{F}) \in (\frac{1}{r}, \frac{2}{2r-1})$.

First assume that forbidding $\mathcal{F}$ implies that some blow-up $T^r[t]$ of $T^r$ is also forbidden (recall, $T^r$ denotes the triangle with three edges on $r+1$ vertices).
Corollary~\ref{blow up invariance} gives $\gamma^+(T^r[t]) = \gamma^+(T^r)$.
As 
\[
\gamma^+(\mathcal{F}) \leq \gamma^+(T^r[t]) = \gamma^+(T^r),
\]
it is sufficient to show that $\gamma^+(T^r) \leq \frac1r$.

Suppose that $H$ is an $n$-vertex $r$-graph with $\delta_{r-1}^+(H) > \frac{n}{r}$, and let $v_1v_2\dots v_r$ be an $r$-edge of $H$. Consider the $r$ vertex sets, each  of size $r-1$, contained in the $r$-edge $v_1v_2\dots v_r$. Each of them is a set with positive co-degree, hence each has neighborhood of size greater than $\frac{n}{r}$. Since there are $r$ such sets, there must be a vertex, say $v_{r+1}$, which is contained in at least two such neighborhoods. Relabeling if needed, we may assume that 
\[
v_{r+1} \in N(v_1, \dots, v_{r-2}, v_{r-1}) \cap N(v_1, \dots, v_{r-2}, v_{r}).
\]

The three $r$-edges $v_1v_2\dots v_{r-2}v_{r-1}v_r$, $v_1v_2\dots v_{r-2}v_{r-1}v_{r+1}$, and $v_1v_2\dots v_{r-2}v_{r}v_{r+1}$ form \acopyof $T^r$. Thus, $\frac{1}{r} \geq \gamma^+\left(T^r\right) \geq  \gamma^+\left(\mathcal{F}\right)$. This implies that blow-ups of $T^r$ are $\mathcal{F}$-free.

Now assume that forbidding $\mathcal{F}$ does not exclude any blow-up of $T^r$.
Consider the following $n$-vertex blow-up of $T^r$. The three vertices corresponding to the base triangle $T^2$ are blown up to classes of size $\frac{n}{2r-1}$. All other vertices are blown up to classes of size $\frac{2n}{2r-1}$. In total, we have $3$ classes of size $\frac{n}{2r-1}$ and $r-2$ classes of size $\frac{2n}{2r -1}$, for a total of $n$ vertices, as desired. A set of $r-1$ vertices whose intersection with any class has at least two vertices will have co-degree $0$. A set of $r-1$ vertices which intersects all three classes of size $\frac{n}{2r-1}$ will also have co-degree $0$. All other sets of size $r-1$ have neighborhood of size exactly $\frac{2n}{2r-1}$, corresponding either to two classes of size $\frac{n}{2r-1}$ or one class of size $\frac{2n}{2r-1}$.
This construction implies $\gamma^+(\mathcal{F}) \geq \frac{2}{2r-1}$.

Now, we can characterize families of $r$-graphs $\mathcal{F}$ with $\gamma^+(\mathcal{F}) \in \{0, \frac{1}{r}\}$. Corollary~\ref{blow up invariance} establishes that if $F 
\in \mathcal{F}$ is $r$-partite, then $\gamma^+(F) = 0$, so $\gamma^+(\mathcal{F}) = 0$ as well. If no $F \in \mathcal{F}$ is  $r$-partite, then any blow-up of an $r$-edge is $\mathcal{F}$-free, so $\gamma^+(\mathcal{F}) \geq \frac{1}{r}$. Thus, $\gamma^+(\mathcal{F}) = 0$ if and only if some $F \in \mathcal{F}$ is $r$-partite. Similarly, if an $F \in \mathcal{F}$ is contained in a blow-up of $T^r$, then 
\[
 \gamma^+(\mathcal{F}) 
\leq \gamma^+(F) \leq \gamma^+(T^r) \leq \frac{1}{r}.\]
Thus, if some $F \in \mathcal{F}$ is contained in some $T^r$ blow-up but no member of $\mathcal{F}$ is $r$-partite, then we have $\gamma^+(F) = \frac{1}{r}$. On the other hand, if no member of $\mathcal{F}$ is contained in any  blow-up of $T^r$, then the above-described blow-up of $T^r$ establishes that $\gamma^+(\mathcal{F}) \geq \frac{2}{2r - 1}$.
\end{proof}

\begin{remark}
The characterization in Theorem~\ref{general jumps} implies the positive co-degree densities for a variety of natural $3$-graphs.  It is straightforward to verify that $C_{\ell}^-$ (the $\ell$-vertex (tight) cycle with one edge deleted) is contained in a sufficiently large blow-up of $K_4^{3-}$. Moreover, $C_{\ell}^-$ is $3$-partite if and only if $\ell \equiv 0 \pmod 3$. Thus, $\gamma^+(C_{\ell}^-) = 0$ when $\ell \equiv 0 \pmod 3$, and $\gamma^+(C_{\ell}^-) = \frac{1}{3}$ otherwise. This generalizes the result of Wu~\cite{Wu} that $\gamma^+(C_5^-) = \frac{1}{3}$.
\end{remark}

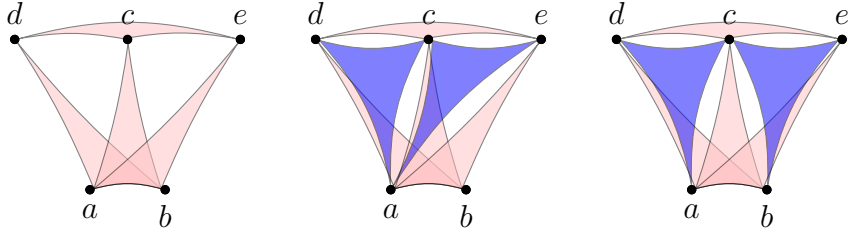
\begin{figure}[h]
    \centering

\begin{tikzpicture}


\draw (0,0) node[vtx](a1){};
\draw (1,0) node[vtx](b1){};
\draw (-1,2) node[vtx](c1){};
\draw (0.5, 2) node[vtx](d1){};    
\draw (2, 2) node[vtx](e1){}; 
\draw[fill=pink,opacity=0.5]
(a1)to[bend left = 15] (b1)to[bend left=5] (c1)to[bend left=5] cycle
;
\draw[fill=pink,opacity=0.5]
(a1)to[bend left = 15] (b1)to[bend left=5] (d1)to[bend left=5] cycle
;
\draw[fill=pink,opacity=0.5]
(a1)to[bend left = 15] (b1)to[bend left=5] (e1)to[bend left=5] cycle
;
\draw[fill=pink,opacity=0.5]
(c1)to[bend left = 10] (d1)to[bend left=10] (e1)to[bend right=15] cycle
;
\draw (0,0) node[vtx,label=below:$a$](a1){};
\draw (1,0) node[vtx,label=below:$b$](b1){};
\draw (-1,2) node[vtx,label=above:$d$](c1){};
\draw (0.5, 2) node[vtx,label=above:$c$](d1){};    
\draw (2, 2) node[vtx,label=above:$e$](e1){}; 


\draw (4,0) node[vtx](a2){};
\draw (5,0) node[vtx](b2){};
\draw (3,2) node[vtx](c2){};
\draw (4.5, 2) node[vtx](d2){};    
\draw (6, 2) node[vtx](e2){}; 
\draw[fill=pink,opacity=0.5]
(a2)to[bend left = 15] (b2)to[bend left=5] (c2)to[bend left=5] cycle
;

\draw[fill=pink,opacity=0.5]
(a2)to[bend left = 15] (b2)to[bend left=5] (d2)to[bend left=5] cycle
;

\draw[fill=pink,opacity=0.5]
(a2)to[bend left = 15] (b2)to[bend left=5] (e2)to[bend left=5] cycle
;

\draw[fill=pink,opacity=0.5]
(c2)to[bend left = 10] (d2)to[bend left=10] (e2)to[bend right=15] cycle
;

\draw[fill=blue,opacity=0.5]
(c2)to[bend right = 15] (d2)to[bend right=15] (a2)to[bend right=15] cycle
;

\draw[fill=blue,opacity=0.5]
(d2)to[bend right = 15] (e2)to[bend right=12] (a2)to[bend right=12] cycle
;

\draw (4,0) node[vtx, label=below:$a$](a2){};
\draw (5,0) node[vtx, label=below:$b$](b2){};
\draw (3,2) node[vtx, label=above:$d$](c2){};
\draw (4.5, 2) node[vtx, label=above:$c$](d2){};    
\draw (6, 2) node[vtx, label=above:$e$](e2){};


\draw (8,0) node[vtx](a3){};
\draw (9,0) node[vtx](b3){};
\draw (7,2) node[vtx](c3){};
\draw (8.5, 2) node[vtx](d3){};    
\draw (10, 2) node[vtx](e3){}; 
\draw[fill=pink,opacity=0.5]
(a3)to[bend left = 15] (b3)to[bend left=5] (c3)to[bend left=5] cycle
;

\draw[fill=pink,opacity=0.5]
(a3)to[bend left = 15] (b3)to[bend left=5] (d3)to[bend left=5] cycle
;

\draw[fill=pink,opacity=0.5]
(a3)to[bend left = 15] (b3)to[bend left=5] (e3)to[bend left=5] cycle
;

\draw[fill=pink,opacity=0.5]
(c3)to[bend left = 10] (d3)to[bend left=10] (e3)to[bend right=15] cycle
;

\draw[fill=blue,opacity=0.5]
(c3)to[bend right = 15] (d3)to[bend right=15] (a3)to[bend right=15] cycle
;

\draw[fill=blue,opacity=0.5]
(d3)to[bend right = 15] (e3)to[bend right=15] (b3)to[bend right=15] cycle
;

\draw (8,0) node[vtx, label=below:$a$](a3){};
\draw (9,0) node[vtx, label=below:$b$](b3){};
\draw (7,2) node[vtx, label=above:$d$](c3){};
\draw (8.5, 2) node[vtx, label=above:$c$](d3){};    
\draw (10, 2) node[vtx, label=above:$e$](e3){};

\end{tikzpicture}

    \caption{$F_{3,2}$ and two superhypergraphs, $F_{3,2}^{++1}$ and $F_{3,2}^{++2}$.}
    \label{F32}
\end{figure}

The next natural question is whether Theorem \ref{general jumps} is best possible. That is, does there exist some family $\mathcal{F}$ of $r$-graphs with $\gamma^+(\mathcal{F}) = \frac{2}{2r - 1}$? We answer this question in the affirmative when $r = 3$. Furthermore, we show that an infinite number of densities in the interval $[\frac{2}{5}, \frac{1}{2}]$ are achievable when $r = 3$. We begin by exhibiting a family $\mathcal{F}$ with $\gamma^+(\mathcal{F}) = \frac{2}{5}$.  We first define $F_{3,2}^{++1}$ and $F_{3,2}^{++2}$, two superhypergraphs of $F_{3,2}$. Each is created by adding two edges to $F_{3,2}$. Let $F_{3,2}$ have vertex set $\{a,b,c,d,e\}$ and edge set $\{abc, abd, abe, cde\}$. Then we create $F_{3,2}^{++1}$ by adding $acd$ and $ace$, and we create $F_{3,2}^{++2}$ by adding $acd$ and $bce$, see Figure~\ref{F32}.

\begin{prop}\label{2/5 family}
For $\mathcal{F} = \{K_4^3, F_{3,2}^{++1},  F_{3,2}^{++2}, J_4\}$ and $\mathcal{F'} = \{K_4^3, F_{3,2}, J_4\}$ we have
\[ 
\gamma^+(\mathcal{F}) = \gamma^+(\mathcal{F'}) = \frac{2}{5}.
\]

\end{prop}

\begin{proof}
Since $F_{3,2}$ is a subgraph of both $F_{3,2}^{++1},  F_{3,2}^{++2}$, hence
$\gamma^+(\mathcal{F}) \geq \gamma^+(\mathcal{F'})$.
As noted in the proof of Theorem~\ref{general jumps}, an appropriately balanced $n$-vertex blow-up of $T^3=K_4^{3-}$ has minimum positive co-degree $\frac{2n}{5} + O(1)$. 
Any blow-up of $K_4^{3-}$ is $J_4$-free and $F_{3,2}$-free since $K_4^{3-}$ is $4$-partite and neither of  $F_{3,2}$ and $J_4$ are. 
Since  blow-ups  $K_4^{3-}$ are also $K_4^{3}$-free, we have
$\frac{2}{5} \leq \gamma^+(\mathcal{F'}) \leq \gamma^+(\mathcal{F})$.

Next we show $\gamma^+(\mathcal{F}) \leq \frac{2}{5}$, which completes the proof.
Fix an $\varepsilon > 0$ and suppose that $H$ is an $n$-vertex $3$-graph with $\delta_2^+
(H) \geq \left( \frac{2}{5} + \varepsilon \right)n$ for sufficiently large $n$. 
Since $\gamma^+(K_4^{3-})=\frac13$, $H$ contains \acopyof $K_4^{3-}$, say with vertex set $\{a,b,c,d\}$ and edge set $\{abc, abd, acd\}$. 

Consider the following five positive co-degree pairs: $ab, ac, ad, bc,$ and $bd$. Since $\delta_2^+(H) \geq \left(\frac{2}{5} + \varepsilon \right)n,$ there exists some vertex $e$ that is in the neighborhood of at least three of these pairs. Note that $c$ and $d$ are symmetric; up to this symmetry, we have six cases based on which three of these five pairs form an edge with $e$.

\begin{enumerate}[label=(\arabic*),itemsep=0.9pt, parsep=0pt]

\item If $abe, ace, ade \in E(H)$ then $abc, abd, acd, abe, ace, ade$ form a $J_4$.

\item If $abe, ace, bce \in E(H)$ then $abc, abe, ace, bce$ form a $K_4^3$.

\item If $abe, ace, bde \in E(H)$ 
then $acb$, $acd$, $ace$, $bde$ form an $F_{3,2}$. Moreover, $abd$ and $abe$ are $3$-edges, so $\{a,b,c,d,e\}$  spans an \acopyof $F_{3,2}^{++1}$.

\item If $ace, ade, bce \in E(H)$ 
then $adb, adc, ade, bce$ form an $F_{3,2}$. Moreover, $abc$ and $ace$ are $3$-edges, so $\{a,b,c,d,e\}$ spans an \acopyof $F_{3,2}^{++1}$. 

\item If $abe, bce, bde \in E(H)$ 
then $bea, bec, bed, acd$ form an $F_{3,2}$. Moreover, $bac$ and $bad$ are $3$-edges, so $\{a,b,c,d,e\}$  spans an \acopyof $F_{3,2}^{++1}$.

\item If $ace, bce, bde \in E(H)$ 
then $acb, acd, ace, bde$ form an $F_{3,2}$. Moreover, $abd$ and $cbe$ are $3$-edges, so $\{a,b,c,d,e\}$ spans an  \acopyof $F_{3,2}^{++2}$.
\end{enumerate}

This implies $\gamma^+(\mathcal{F}) \leq \frac{2}{5}$.
\end{proof}

We use Proposition~\ref{2/5 family} as the base of an inductive argument, showing that an infinite number of positive co-degree densities in $[\frac{2}{5}, \frac{1}{2}]$ are achievable by families of $3$-graphs.

\begin{proof}[Proof of Theorem~\ref{densities}]
For a lower bound, observe that any blow-up of $J_{k-1}$ is $\{K_4^3, F_{3,2}, J_k\}$-free, and that an appropriately balanced $n$-vertex blow-up of $J_{k-1}$ (with class sizes $\frac{k-2}{2k-3},$ $ \frac{1}{2k-3},\ldots, \frac{1}{2k-3}$) has minimum positive co-degree $\left(\frac{k-2}{2k-3}\right)n + O(1)$. Thus, $\gamma^+(\{K_4^3, F_{3,2}, J_k\}) \geq \frac{k-2}{2k-3}$ for all $k \geq 4$.

To show that this lower bound is best-possible, we use induction on $k$. 
Proposition~\ref{2/5 family} yields the statement for $k=4$.
We assume the statement holds for $k \geq 4$, and prove  it  for $k + 1$. Suppose that $H$ is an $n$-vertex $3$-graph with minimum positive co-degree $\delta_2^+(H) >  \frac{k-1}{2k-1} n$ where is $n$ sufficiently large. If $H$ contains one of $K_4^3$ or $F_{3,2}$, then we are done, so suppose not. Then by the inductive hypothesis, $\gamma^+(\{K_4^3, F_{3,2}, J_k\}) = \frac{k-2}{2k-3} < \frac{k-1}{2k -1}$, so $H$ contains a $J_k$, say $J$ (here we use that $n$ is large enough). Let $V(J) = \{v_1, \dots, v_{k+1}\}$, where $v_1$ is the universal vertex of $J$ (i.e., $E(J) = \{v_1v_iv_j : 2 \leq i < j \leq k+1\}$).  Define 
\[ S = \{(v_1, v_i) : 2 \leq i \leq k+1 \} \cup \{(v_i,v_{i+1}): 2 \leq i \leq k\}.\]

Observe that $|S| = 2k - 1$ and every vertex pair in $S$ has positive co-degree in $H$. Since $\delta_2^+(H) >  \frac{k-1}{2k-1} n$, there exists a vertex $u \in V(H)$ which is in the neighborhood of at least $k$ pairs  $(v_i,v_j) \in S$. Note that $u$ may be an element of $V(J)$. However, if $(v_i,v_j) \in S$ is a pair such that $u \in N(v_i,v_j)$, then $u \not\in \{v_i,v_j\}$.

If $u \in N(v_1, v_i)$ for every $i \in \{2, \dots, k+1\}$ then $u \not\in V(J)$, and $V(J) \cup \{u\}$ spans \acopyof $J_{k+1}$, with universal vertex $v_1$.

Hence, we may  assume that there is some $i \geq 2$ for which $u \in N(v_i, v_{i+1})$.
Our goal in this case is to find $F_{3,2}$. Observe that there must be a $j$ such that $u \in N(v_1, v_j)$. If there is such a $j \not \in \{i, i+1\}$, then $u \not\in \{v_1,v_i, v_{i+1}, v_j\}$, and $\{v_1, v_i, v_{i+1}, v_j, u\}$ will span an $F_{3,2}$, using $3$-edges $v_1v_jv_i, v_1v_jv_{i+1}, v_1v_ju, v_iv_{i+1}u$. Hence, we may assume that every such $j$ is in $\{i, i+1\}$.  In particular, $u$ is in at most two neighborhoods of the form $N(v_1,v_j)$. Since $k \geq 4$ and $u$ is in the 
neighborhood of at least $k$ pairs from $S$, we  have $u \in N(v_{\ell},v_{\ell + 1})$ for some $\ell \neq i$. Now, if $u \in N(v_1, v_i) $ and $u \in N(v_1, v_{i+1})$, then $\{v_1,v_i,v_{i+1},u\}$ spans a $K_4^3$. If not, then without loss of generality $u \in N(v_1, v_i) $ and $u$ is in the neighborhood of the three pairs  $(v_i, v_{i+1}), (v_{\ell}, v_{\ell + 1}), (v_m, v_{m+1})$ for some $\ell,m$. One of these pairs, say $(v_{\ell}, v_{\ell + 1})$, must be disjoint from $(v_1, v_i)$, in which case we observe that $\{v_1,v_i, v_\ell, v_{\ell+1},u\}$ will span an $F_{3,2}$. 
\end{proof}

Although infinitely many positive co-degree densities in $[\frac{2}{5}, \frac{1}{2}]$ can be achieved by forbidden families of $3$-graphs, it remains unclear whether all of these densities can be achieved by forbidding a single $3$-graph. Using Proposition~\ref{2/5 family}, we now show that $\frac{2}{5}$ is indeed achievable by a single $3$-graph. We begin by defining two new $3$-graphs $F_1$ and $F_2$ as follows.

\begin{align*}V(F_1) &  = \{a,b,c,d,e,f,g\}, \quad  & E(F_1) & = \{abe, ace, bce, abf, adf, bdf, cdg\},\\
V(F_2) &  = \{a,b,c,d,e,g\}, \quad  &E(F_2) & = \{abe, ace, bce, ade, bde, cdg\}.
\end{align*}

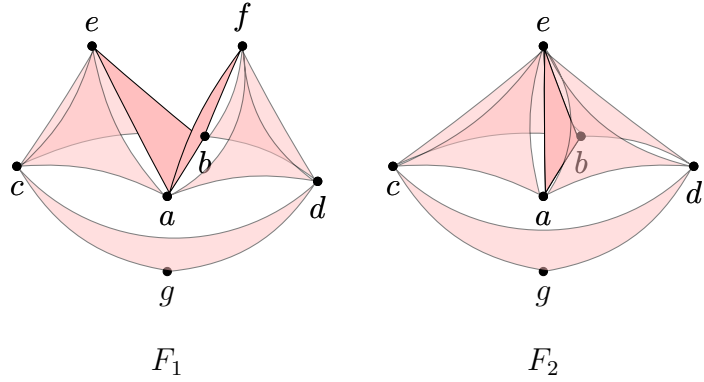
\begin{figure}[h]
    \centering

\begin{tikzpicture}
    
\draw (5,0.4) node[vtx,label=below:$c$](x1){};
\draw (7,0) node[vtx,label=below:$a$](z1){};
\draw (7.5,0.8) node[vtx,label=below:$b$](w1){};
\draw (9, 0.4) node[vtx,label=below:$d$](y1){};    
\draw (7, 2) node[vtx,label=above:$e$](f){}; 

\draw (7, -1) node[vtx,label=below:$g$](g1){};

\draw[fill=pink,opacity=0.5]
(x1) to[bend left=15] (w1)to[bend left=15] (f) to[bend left=2] cycle
;
\draw[fill=pink,opacity=0.5]
(x1) to[bend left=15] (z1)to[bend left=15] (f) to[bend left=15] cycle
;

\draw[fill=pink,opacity=1]
(z1) to[bend left=1] (w1)to[bend left=1] (f) to[bend right=1] cycle
;

\draw[fill=pink,opacity=0.5]
(y1) to[bend right=15] (w1)to[bend left = 25] (f) to[bend right=2] cycle
;

\draw[fill=pink,opacity=1]
(z1) to[bend left=1] (w1)to[bend left=1] (f) to[bend right=1] cycle
;

\draw (7.5,0.8) node[vtx,label = below:$b$](w1){};

\draw[fill=pink,opacity=0.5]
(y1) to[bend right=15] (z1)to[bend right=35] (f) to[bend right=25] cycle
;

\draw[fill=pink,opacity=0.5]
(x1) to[bend right=45] (y1) to[bend left=25] (g1) to[bend left=25] cycle
;

\draw (5,0.4) node[vtx,label=below:$c$](x1){};
\draw (7,0) node[vtx,label=below:$a$](z1){};

\draw (9, 0.4) node[vtx,label=below:$d$](y1){};    
\draw (7, 2) node[vtx,label=above:$e$](f){}; 

\draw (0,0.4) node[vtx,label=below:$c$](x2){};
\draw (2,0) node[vtx,label=below:$a$](z2){};
\draw (2.5,0.8) node[vtx,label=below:$b$](w2){};
\draw (4, 0.2) node[vtx,label=below:$d$](y2){};    
\draw (1, 2) node[vtx,label=above:$e$](f){}; 
\draw (3, 2) node[vtx,label=above:$f$](c){}; 

\draw (2, -1) node[vtx,label=below:$g$](g2){}; 

\draw[fill=pink,opacity=0.5]
(x2) to[bend left=15] (w2)to[bend left=15] (f) to[bend left=2] cycle
;
\draw[fill=pink,opacity=0.5]
(x2) to[bend left=15] (z2)to[bend left=15] (f) to[bend left=15] cycle
;

\draw[fill=pink,opacity=1]
(z2) to[bend left=1] (w2)to[bend left=1] (f) to[bend right=1] cycle
;

\draw[fill=pink,opacity=0.5]
(y2) to[bend right=15] (w2)to[bend right = 15] (c) to[bend right=2] cycle
;

\draw[fill=pink,opacity=1]
(z2) to[bend right=1] (w2)to[bend left=1] (c) to[bend right=10] cycle
;

\draw[fill=pink,opacity=0.5]
(x2) to[bend right=45] (y2)to[bend left=25] (g2) to[bend left=25] cycle
;

\draw (2.5,0.8) node[vtx,label=below:$b$](w2){};

\draw[fill=pink,opacity=0.5]
(y2) to[bend right=15] (z2)to[bend right=35] (c) to[bend right=25] cycle
;

\draw (0,0.4) node[vtx,label=below:$c$](x2){};
\draw (2,0) node[vtx,label=below:$a$](z2){};

\draw (4, 0.2) node[vtx,label=below:$d$](y2){};    
\draw (1, 2) node[vtx,label=above:$e$](f){}; 
\draw (3, 2) node[vtx,label=above:$f$](c){};


\draw (7, -2.2) node{$F_2$}; 
\draw (2, -2.2) node{$F_1$};

\end{tikzpicture}
    \caption{$F_1$ and $F_2$.}
    \label{F1F2}
\end{figure}

We depict $F_1$ and $F_2$ in Figure~\ref{F1F2}.
We remark that each of $F_1,F_2$ can be viewed as a partial identification of two $K_4^{3-}$ copies, along with an extra $3$-edge (using $g$) which ensures that $c$ and $d$ have positive co-degree.

Note that in both $F_1$ and $F_2$, the vertex $g$ has very little structural interaction with the other vertices.
Suppose $H$ is isomorphic to the subhypergraph of $F_1$ (resp.\ $F_2$) induced on $V(F_1) \setminus \{g\}$ (resp.\ $V(F_2) \setminus \{g\}$). Then $H$ is guaranteed to extend to \acopyof $F_1$ (resp.\ $F_2$) if $d(c,d) > 4$.
We will be working in $3$-graphs with minimum positive co-degree much larger than $4$, so to find copies of $F_1$ or $F_2$, it will suffice to find $F_1 \setminus \{g\}$ or $F_2 \setminus \{g\}$ and to demonstrate that the pair $\{c,d\}$ has positive co-degree.

Note that $F_2$ is $5$-partite, so $\gamma^+(F_2) \geq \frac{1}{2}$. However, $F_1$ is $4$-partite, and is contained in a blow-up of $J_4$, but it is not contained in a blow-up of $J_3=K_4^{3-}$. Thus, $\gamma^+(F_1) \geq \frac{2}{5}$. 
We now prove that $\gamma^+(F_1) = \frac{2}{5}$.

\begin{proof}[Proof of Theorem~\ref{single graph}]
For a contradiction, assume that there exists $\varepsilon > 0$ such that $\gamma^+(F_1) \geq \frac{2}{5} + 3\varepsilon$. We begin with a claim that will allow us to expand our forbidden family.

\begin{claim}
The $3$-blow-up of each member of the following family
\[\mathcal{F} = \{K_4^3, J_4, F_{3,2}^{++1}, F_2 \}\]
 contains $F_1$.
 \end{claim}

 \begin{proof}
Observe that $F_2$ can be obtained from $F_1$ by identifying $e$ and $f$.
Using the notation
\begin{align*}
V(K_4^3) &=\{1,2,3,4\}, &E&(K_4^3) \quad =\{123,124,134,234\},\\
V(J_4) &=\{1,2,3,4,5\}, &E&(J_4)\quad \ =\{123, 124, 125, 134, 135, 145\},\\
V(F_{3,2}^{++1})&=\{1,2,3,4,5\},  &E&(F_{3,2}^{++1})=\{123,  124,  125,  345,  134,  135\}
\end{align*}
the following maps prove the claim for $K_4^3$, $J_4$, and $F_{3,2}^{++1}$.


$f_1: F_1\to K_4^3:\  a\to 1,\  b\to 2,\ c \to 3, \ d \to 4, \ e \to 4,\ f \to 3, \ g \to 1$.

$f_2: F_1\to J_4: \ \ a\to 2, \ b\to 3,\ c \to 4,\ d \to 5, \ e \to 1, \ f \to 1, \ g \to 1$.

$f_3: F_1\to F_{3,2}^{++1}: a\to 2,\ b\to 3,\ c \to 4, \ d \to 5,\  e \to 1, \ f \to 1, \ g \to 3$.
%

\end{proof}

 Thus, by Lemma~\ref{family removal}, for $n$ large enough there exists a $3$-graph $H$ with $\delta_2^+(H) \geq \left(\frac{2}{5} + \varepsilon\right)n$ which is $\{F_1\} \cup \mathcal{F}$-free. By Proposition~\ref{2/5 family}, 
\[
\gamma^+(\{K_4^3, J_4, F_{3,2}^{++1}, F_{3,2}^{++2}\}) = \frac{2}{5},
\]
so  $H$ contains \acopyof $F_{3,2}^{++2}$, which we shall call $F$. Put 
\[V(F) = \{a,b,c,d,e\},\quad \text{ with } \quad E(F) = \{abc, abd, abe, cde, acd, bce\}.\] 

Observe that every pair of vertices in $V(F)$ has positive co-degree. Since $\delta_2^+(H) \geq \left( \frac{2}{5} + \varepsilon \right)n$, there exists some vertex $f$ that is in the neighborhood of at least five pairs from $V(F)$.
Let $G_f$ be the link graph of vertex $f$ induced on vertex set $V(F)$ (i.e.\ $G_f = L(f)[V(F)]$).

\begin{claim}\label{VFclaims}
The following statements hold:
\begin{enumerate}[label=(\roman*),itemsep=1pt, parsep=0pt]
    \item  For every $xyz \in E(H)$,  at most two of $xy, xz, yz$ are in $E(G_f)$. 
    
\item There is no isolated vertex in $G_f$.

 \item At most one of $ad, cd$ is in $E(G_f)$ and
 at most one of $be, ce$ is in $E(G_f)$.

\item At most one of $cd, de$ is in $E(G_f)$ and  at most one of $ce, de$ is in $E(G_f)$. 

\item At most one of $ab, ac$ is in $E(G_f)$ and at most one of $ab, bc$ is in $E(G_f)$.

\item At most one of $ab, ad$ is in $E(G_f)$ and  at most one of $ab, be$ is in $E(G_f)$.

\item At most one of $ac, ad$ is in $E(G_f)$ and  at most one of $bc, be$ is in $E(G_f)$.

\item $G_f$ is triangle-free.

\end{enumerate}


\end{claim}
\begin{proof}
To prove each statement, we shall assume that the statement does not hold and use this assumption to find some forbidden structure in $H$. 
Notice that the two parts in each of the statements (iii)--(vii) are symmetric so it is sufficient to prove only the first part of each of those claims. Here the symmetry is coming from the fact that $a\to b,\ b\to a,\ c\to c,\ d\to e, \ e\to d$ is an automorphism of $F$.

For (i), observe that if $xyz \in E(H)$ and $xy, xz, yz \in E(G_f)$, then $\{x,y,z,f\}$ spans a  $K_4^3$ in $H$.

For (ii), assume that $G_f$ contains a vertex of degree 0. Then the other four vertices in $V(F)$ span (at least) five edges of $G_f$, so $G_f$ contains a \acopyof $K_4^-$. As illustrated in Figure~\ref{k4-link}, a \acopyof $K_4^-$ 

in $G_f$ implies that $H$ contains $F_2$ if the appropriate pair of vertices ($x$ and $y$, in the figure) has positive co-degree.
 Since all pairs of vertices in $F$ have positive co-degree, \acopyof a $K_4^-$ in $G_f$ indeed implies that $H$ contains \acopyof $F_2$.

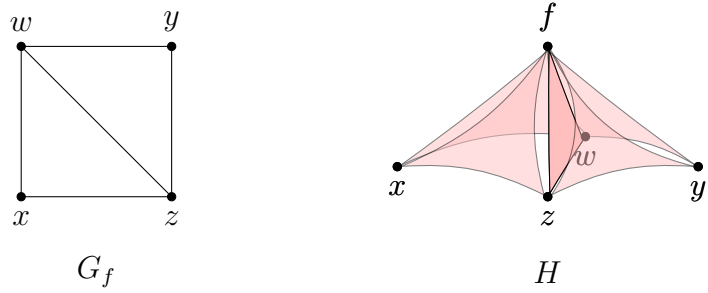
\begin{figure}[h]
\centering 

\begin{tikzpicture}

\draw (0,0) node[vtx,label=below:$x$](x){};
\draw (2,0) node[vtx,label=below:$z$](z){};
\draw (0,2) node[vtx,label=above:$w$](w){};
\draw (2,2) node[vtx,label=above:$y$](y){};

\draw 
    (x) to (z)
    (z) to (y)
    (z) to (w)
    (w) to (x)
    (w) to (y);

\draw (1, -1) node{$G_f$};

\draw (5,0.4) node[vtx,label=below:$x$](x1){};
\draw (7,0) node[vtx,label=below:$z$](z1){};
\draw (7.5,0.8) node[vtx,label=below:$w$](w1){};
\draw (9, 0.4) node[vtx,label=below:$y$](y1){};    
\draw (7, 2) node[vtx,label=above:$f$](f){}; 

\draw[fill=pink,opacity=0.5]
(x1) to[bend left=15] (w1)to[bend left=15] (f) to[bend left=2] cycle
;
\draw[fill=pink,opacity=0.5]
(x1) to[bend left=15] (z1)to[bend left=15] (f) to[bend left=15] cycle
;

\draw[fill=pink,opacity=1]
(z1) to[bend left=1] (w1)to[bend left=1] (f) to[bend right=1] cycle
;

\draw[fill=pink,opacity=0.5]
(y1) to[bend right=15] (w1)to[bend left = 25] (f) to[bend right=2] cycle
;

\draw[fill=pink,opacity=1]
(z1) to[bend left=1] (w1)to[bend left=1] (f) to[bend right=1] cycle
;

\draw (7.5,0.8) node[vtx,label = below:$w$](w1){};

\draw[fill=pink,opacity=0.5]
(y1) to[bend right=15] (z1)to[bend right=35] (f) to[bend right=25] cycle
;

\draw (5,0.4) node[vtx,label=below:$x$](x1){};
\draw (7,0) node[vtx,label=below:$z$](z1){};

\draw (9, 0.4) node[vtx,label=below:$y$](y1){};    
\draw (7, 2) node[vtx,label=above:$f$](f){}; 

\draw (7,-1) node{$H$};

\end{tikzpicture}
\caption{\acopyof $K_4^-$ 
in $G_f$ and the resulting structure in $H$.} \label{k4-link}

\end{figure}

For (iii), observe that if $ad, cd \in E(G_f)$, then $\{a,c,d,f\}$ spans \acopyof $K_4^{3-}$ 
with spike vertex $d$, and $\{a, b, c, e\}$  spans \acopyof $K_4^{3-}$  with spike vertex $b$, so $\{a,b,c,d,e, f\}$ spans \acopyof $F_1$ (since  $d(e,f) > 0$ by (ii)).


For (iv), if $cd, de \in E(G_f)$, then $\{c,d,e,f\}$ spans \acopyof $K_4^{3-}$ with spike vertex $d$, and $\{a,b,c,e\}$ spans \acopyof $K_4^{3-}$ with spike vertex $b$, so $\{a,b,c,d,e,f\}$ spans \acopyof $F_1$ (since $d(a,f) > 0$ by (ii)).


For (v), if $ab, ac \in E(G_f)$, then  $\{a,b,c,d\}$ and $\{a, b, c, f\}$ span copies of $K_4^{3-}$  with spike vertex $a$, so $\{a,b,c,d,f\}$ spans \acopyof $F_2$ (since $d(d,f) > 0$ by (ii)).


For (vi), if $ab, ad \in E(G_f)$, then  $\{a,b,c,d\}$ and $\{a, b, d, f\}$ span copies of $K_4^{3-}$ with spike vertex $a$, so $\{a,b,c,d,f\}$ spans \acopyof $F_2$ (since $d(c,f) > 0$ by (ii)).


For (vii), if $ac, ad \in E(G_f)$, then $\{a,b,c,d\}$ and $\{a, c, d, f\}$ span copies of $K_4^{3-}$ with spike vertex $a$, so $\{a,b,c,d,f\}$ spans \acopyof $F_2$ (since $d(b,f) > 0$ by (ii)).


Finally, we show (viii). From (i), we know that $\{c,d,e\}$ does not form a triangle in $G_f$, and from (v) and (vi) it follows that $ab$ is not contained in a triangle in $G_f$. Thus, any triangle in $G_f$ must include one edge spanned by $\{c,d,e\}$, and two edges with one vertex in $\{a,b\}$ and the other in $\{c,d,e\}$. By (i), $\{a,c,d\}$ does not span a triangle (nor does $\{b, c, e\}$). Up to symmetry, there are two other potential triangles: $\{a, d, e\}$, and $\{a,c,e\}$. Observe that if $\{a,d,e\}$ spans a triangle in $G_f$, then $\{a,d,e, f\}$ spans \acopyof $K_4^{3-}$ with spike vertex $f$, and $\{a, b, c, e\}$ spans \acopyof $K_4^{3-}$ with spike vertex $b$, so $\{a,b,c,d,e,f\}$ spans \acopyof $F_1$ (since $d(c,d) > 0$). 

Next, suppose $\{a,c,e\}$ spans a triangle in $G_f$. By (iii), (iv), (v), and (vii), none of $be, de, ab, ad$ are in $E(G_f)$. Thus, two of $bc, bd, cd$ are in $E(G_f)$. We shall show that this is impossible. Observe that if $bc \in E(G_f)$, then $\{a,b,c,f\}$ and $\{b,c,e,f\}$ each spans a $K_4^{3-}$  with spike vertex $c$, so $\{a,b,c,e,f\}$ spans \acopyof $F_2$ (since $d(a,e) > 0$). Moreover, if $cd \in E(G_f)$, then $\{a,c,d,f\}$ and $\{c,d,e,f\}$ span $K_4^{3-}$ with spike vertex $c$, so $\{a,c,d,e,f\}$ spans \acopyof $F_2$ (since $d(a,e) > 0$). Thus, $G_f$ is triangle-free.
\end{proof}

With Claim~\ref{VFclaims} established, we are ready to determine the structure of $G_f$. By (i), at most two out of five edges of $G_f$ are spanned by $\{c,d,e\}$, so one of $a,b$ has degree at least $2$. Without loss of generality, $d(a) \geq d(b)$ and $d(a) \geq 2$. 

\begin{observation}
    $ab \not\in E(G_f)$. 
\end{observation}

\begin{proof}
Suppose to the contrary that $ab \in E(G_f)$. By (v) and (vi), $ac$ and $ad$ are not in $E(G_f)$, so we must have $ae \in E(G_f)$. By (v) and (vi), we also have $bc, be \not\in E(G_f)$. Since $G_f$ has five edges, at least three of $bd, cd, ce, de$ are in $E(G_f)$. By (iv), at most one of $ce, de$ is in $E(G_f)$, so $bd, cd \in E(G_f)$. Also by (iv), if $cd \in E(G_f)$, then $de \not \in E(G_f)$, so $ce \in E(G_f)$. 
Since $ae \in E(G_f)$, there exists a vertex $x$ not in $V(F)$ that is a common neighbor of $a$ and $e$.
The set of edges $\{adb,dce, fab, fbd, fdc, fce, aex\}$ forms $F_1$. See Figure~\ref{noab} for an illustration of $G_f$, $F$, and $F_1$.
\end{proof}

\begin{figure}[h]
    \centering

\begin{tikzpicture}

\draw (0,0) node[vtx,label=above:$b$](b){};
\draw (-2,0) node[vtx,label=above:$a$](a){};
\draw (-1,-3) node[vtx,label=below:$c$](c){};
\draw (1,-3) node[vtx,label=below:$e$](e){};
\draw (-3,-3) node[vtx,label=below:$d$](d){};

\draw[fill=pink,opacity=0.5]
(a) to[bend right=15] (b)to[bend right=15] (c) to[bend right=15] cycle
;

\draw[fill=pink,opacity=0.5]
(a) to[bend right=15] (b)to[bend right=15] (d) to[bend right=15] cycle
;

\draw[fill=pink,opacity=0.5]
(a) to[bend right=15] (b)to[bend right=15] (e) to[bend right=15] cycle
;

\draw[fill=pink,opacity=0.5]
(d) to[bend left=15] (c)to[bend left=50] (a) to[bend right=15] cycle
;

\draw[fill=pink,opacity=0.5]
(e) to[bend right=15] (c)to[bend right=50] (b) to[bend left=15] cycle
;

\draw[fill=pink,opacity=0.5]
(c) to[bend left=15] (d)to[bend right=20] (e) to[bend left=15] cycle
;

\draw (0,0) node[vtx,label=above:$b$](b){};
\draw (-2,0) node[vtx,label=above:$a$](a){};
\draw (-1,-3) node[vtx,label=below:$c$](c){};
\draw (1,-3) node[vtx,label=below:$e$](e){};
\draw (-3,-3) node[vtx,label=below:$d$](d){};

\draw[thick, blue] (a) to[bend right = 25] (b);

\draw[thick, blue] (b) to[bend right = 20] (d);

\draw[thick, blue] (a) to[bend left = 20] (e);

\draw[thick, blue] (c) to[bend right = 10] (e);

\draw[thick, blue] (d) to[bend right = 10] (c);

\draw (-1, -4) node{$G_f$ and $F$};

\begin{scope}[yshift = 1cm]
\draw (5,-2.6) node[vtx,label=below:$a$](x1){};
\draw (7,-3) node[vtx,label=below:$d$](z1){};
\draw (7.5,-2.2) node[vtx,label=below:$f$](w1){};
\draw (9, -2.8) node[vtx,label=below:$e$](y1){};    
\draw (6, -1) node[vtx,label=above:$b$](f){}; 
\draw (8, -1) node[vtx,label=above:$c$](c){}; 

\draw (7, -4) node[vtx,label=below:$x$](g1){}; 
\draw[fill=pink,opacity=0.5]
(x1) to[bend right=45] (y1) to[bend left=25] (g1) to[bend left=25] cycle
;

\draw[fill=pink,opacity=0.5]
(x1) to[bend left=15] (w1)to[bend left=15] (f) to[bend left=2] cycle
;
\draw[fill=pink,opacity=0.5]
(x1) to[bend left=15] (z1)to[bend left=15] (f) to[bend left=15] cycle
;

\draw[fill=pink,opacity=1]
(z1) to[bend left=1] (w1)to[bend left=1] (f) to[bend right=1] cycle
;

\draw[fill=pink,opacity=0.5]
(y1) to[bend right=15] (w1)to[bend right = 15] (c) to[bend right=2] cycle
;

\draw[fill=pink,opacity=1]
(z1) to[bend right=1] (w1)to[bend left=1] (c) to[bend right=10] cycle
;

\draw (7.5,-2.2) node[vtx,label=below:$f$](w1){};

\draw[fill=pink,opacity=0.5]
(y1) to[bend right=15] (z1)to[bend right=35] (c) to[bend right=25] cycle
;

\draw (5,-2.6) node[vtx,label=below:$a$](x1){};
\draw (7,-3) node[vtx,label=below:$d$](z1){};
\draw (9, -2.8) node[vtx,label=below:$e$](y1){};    
\draw (6, -1) node[vtx,label=above:$b$](f){}; 
\draw (8, -1) node[vtx,label=above:$c$](c){}; 
\end{scope}

\draw (7,-4) node{$F_1$};

\end{tikzpicture}

    \caption{The configuration when $ab \in E(G_f)$, and the resulting copy of $F_1$.}
    \label{noab}
\end{figure}
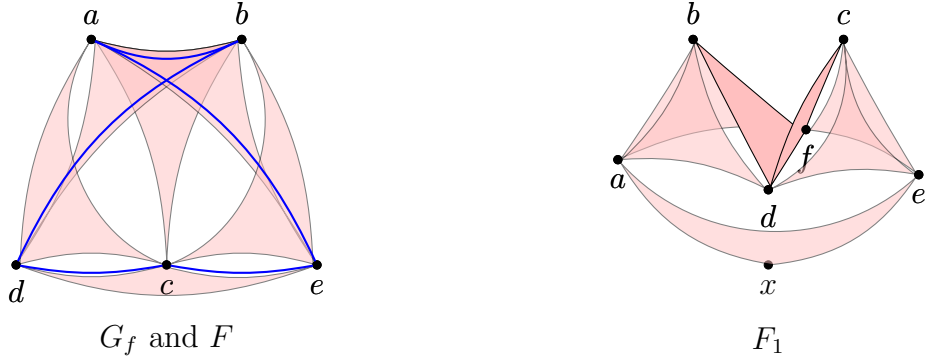

 Since $d(a) \geq 2$, two of $ac, ad, ae$ are in $E(G_f)$. By (vii), at most one of $ac, ad$ is in $E(G_f)$, so $E(G_f)$ contains $ae$ and exactly one of $ac, ad$. 
Suppose first that $ad \in E(G_f)$. By (iii), $cd\not\in E(G_f)$, and by (viii), $de \not\in E(G_f)$. Since $d(b) \leq d(a) = 2$, we must have $ce \in E(G_f)$. Now, by (iii), $be \not\in E(G_f)$, so we must have $bc, bd \in E(G_f)$. 
Since $ae \in E(G_f)$, there exists a vertex $x$ not in $V(F)$ that is a common neighbor of $a$ and $e$.
The set of edges $\{adb,bce, fad, fbd, fbc, fce, aex\}$ forms $F_1$. See Figure~\ref{noad} for an illustration of $G_f$, $F$, and $F_1$.

\begin{figure}[h]
    \centering

\begin{tikzpicture}

\draw (0,0) node[vtx,label=above:$b$](b){};
\draw (-2,0) node[vtx,label=above:$a$](a){};
\draw (-1,-3) node[vtx,label=below:$c$](c){};
\draw (1,-3) node[vtx,label=below:$e$](e){};
\draw (-3,-3) node[vtx,label=below:$d$](d){};

\draw[fill=pink,opacity=0.5]
(a) to[bend right=15] (b)to[bend right=15] (c) to[bend right=15] cycle
;

\draw[fill=pink,opacity=0.5]
(a) to[bend right=15] (b)to[bend right=15] (d) to[bend right=15] cycle
;

\draw[fill=pink,opacity=0.5]
(a) to[bend right=15] (b)to[bend right=15] (e) to[bend right=15] cycle
;

\draw[fill=pink,opacity=0.5]
(d) to[bend left=15] (c)to[bend left=50] (a) to[bend right=15] cycle
;

\draw[fill=pink,opacity=0.5]
(e) to[bend right=15] (c)to[bend right=50] (b) to[bend left=15] cycle
;

\draw[fill=pink,opacity=0.5]
(c) to[bend left=15] (d)to[bend right=20] (e) to[bend left=15] cycle
;

\draw (0,0) node[vtx,label=above:$b$](b){};
\draw (-2,0) node[vtx,label=above:$a$](a){};
\draw (-1,-3) node[vtx,label=below:$c$](c){};
\draw (1,-3) node[vtx,label=below:$e$](e){};
\draw (-3,-3) node[vtx,label=below:$d$](d){};

\draw[thick, blue] (a) to[bend left = 20] (e);

\draw[thick, blue] (b) to[bend right = 20] (d);

\draw[thick, blue] (a) to[bend left = 15] (d);

\draw[thick, blue] (c) to[bend right = 10] (e);

\draw[thick, blue] (c) to[bend right = 55] (b);

\draw (-1, -4) node{$G_f$ and $F$};

\begin{scope}[yshift = 1cm]
\draw (5,-2.6) node[vtx,label=below:$a$](x1){};
\draw (7,-3) node[vtx,label=below:$b$](z1){};
\draw (7.5,-2.2) node[vtx,label=below:$f$](w1){};
\draw (9, -2.8) node[vtx,label=below:$e$](y1){};    
\draw (6, -1) node[vtx,label=above:$d$](f){}; 
\draw (8, -1) node[vtx,label=above:$c$](c){}; 

\draw (7, -4) node[vtx,label=below:$x$](g1){}; 
\draw[fill=pink,opacity=0.5]
(x1) to[bend right=45] (y1) to[bend left=25] (g1) to[bend left=25] cycle
;

\draw[fill=pink,opacity=0.5]
(x1) to[bend left=15] (w1)to[bend left=15] (f) to[bend left=2] cycle
;
\draw[fill=pink,opacity=0.5]
(x1) to[bend left=15] (z1)to[bend left=15] (f) to[bend left=15] cycle
;

\draw[fill=pink,opacity=1]
(z1) to[bend left=1] (w1)to[bend left=1] (f) to[bend right=1] cycle
;

\draw[fill=pink,opacity=0.5]
(y1) to[bend right=15] (w1)to[bend right = 15] (c) to[bend right=2] cycle
;

\draw[fill=pink,opacity=1]
(z1) to[bend right=1] (w1)to[bend left=1] (c) to[bend right=10] cycle
;

\draw (7.5,-2.2) node[vtx,label=below:$f$](w1){};

\draw[fill=pink,opacity=0.5]
(y1) to[bend right=15] (z1)to[bend right=35] (c) to[bend right=25] cycle
;

\draw (5,-2.6) node[vtx,label=below:$a$](x1){};
\draw (7,-3) node[vtx,label=below:$b$](z1){};
\draw (9, -2.8) node[vtx,label=below:$e$](y1){};    
\draw (6, -1) node[vtx,label=above:$d$](f){}; 
\draw (8, -1) node[vtx,label=above:$c$](c){}; 
\end{scope}

\draw (7,-4) node{$F_1$};

\end{tikzpicture}

    \caption{The configuration when $ad, ae \in E(G_f)$, and the resulting copy of $F_1$.}
    \label{noad}
\end{figure}
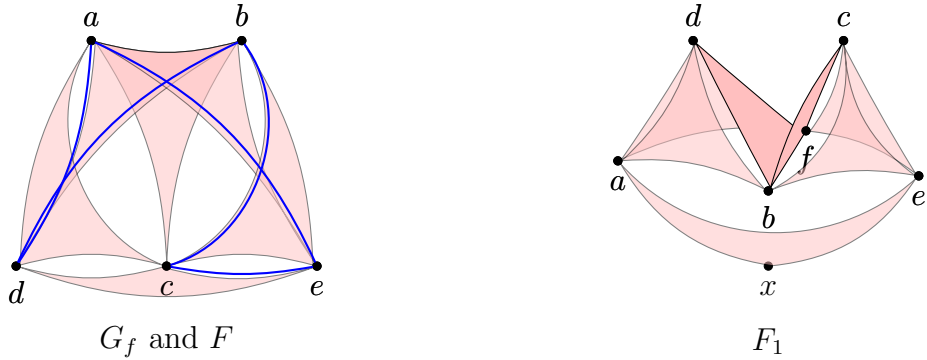

Thus, we have $ac, ae \in E(G_f)$ and $ad \not\in E(G_f)$. By (viii), $ce \not\in E(G_f)$, and by (iv), at most one of $cd, de$ is in $E(G_f)$. Thus, $d(b) = 2$. 
By (vii), at most one of $bc$, $be$ is in $E(G_f)$,
hence $bd \in E(G_f)$. Suppose for a contradiction that $be \in E(G_f)$. Then by (viii), $de\not\in E(G_f)$. Since $G_f$ has at least 5 edges, $dc \in E(G_f)$.
Since $bd \in E(G_f)$, there exists a vertex $x$ not in $V(F)$ that is a common neighbor of $b$ and $d$.
The set of edges $\{abe,adc, fbe, fae, fac, fcd, bdx\}$ forms $F_1$. See Figure~\ref{fig:nobd} for an illustration of $G_f$, $F$, and $F_1$.

\begin{figure}[h]
    \centering

\begin{tikzpicture}

\draw (0,0) node[vtx,label=above:$a$](b){};
\draw (-2,0) node[vtx,label=above:$b$](a){};
\draw (-1,-3) node[vtx,label=below:$c$](c){};
\draw (1,-3) node[vtx,label=below:$d$](e){};
\draw (-3,-3) node[vtx,label=below:$e$](d){};

\draw[fill=pink,opacity=0.5]
(a) to[bend right=15] (b)to[bend right=15] (c) to[bend right=15] cycle
;

\draw[fill=pink,opacity=0.5]
(a) to[bend right=15] (b)to[bend right=15] (d) to[bend right=15] cycle
;

\draw[fill=pink,opacity=0.5]
(a) to[bend right=15] (b)to[bend right=15] (e) to[bend right=15] cycle
;

\draw[fill=pink,opacity=0.5]
(d) to[bend left=15] (c)to[bend left=50] (a) to[bend right=15] cycle
;

\draw[fill=pink,opacity=0.5]
(e) to[bend right=15] (c)to[bend right=50] (b) to[bend left=15] cycle
;

\draw[fill=pink,opacity=0.5]
(c) to[bend left=15] (d)to[bend right=20] (e) to[bend left=15] cycle
;

\draw[thick, blue] (a) to[bend left = 20] (e);

\draw[thick, blue] (b) to[bend right = 20] (d);

\draw[thick, blue] (a) to[bend left = 15] (d);

\draw[thick, blue] (c) to[bend right = 10] (e);

\draw[thick, blue] (c) to[bend right = 55] (b);

\draw (c) node[vtx,label=below:$c$](c){};

\draw (-1, -4) node{$G_f$ and $F$};

\begin{scope}[yshift = 1cm]
\draw (5,-2.6) node[vtx,label=below:$b$](x1){};
\draw (7,-3) node[vtx,label=below:$a$](z1){};
\draw (7.5,-2.2) node[vtx,label=below:$f$](w1){};
\draw (9, -2.8) node[vtx,label=below:$d$](y1){};    
\draw (6, -1) node[vtx,label=above:$e$](f){}; 
\draw (8, -1) node[vtx,label=above:$c$](c){}; 

\draw (7, -4) node[vtx,label=below:$x$](g1){}; 
\draw[fill=pink,opacity=0.5]
(x1) to[bend right=45] (y1) to[bend left=25] (g1) to[bend left=25] cycle
;

\draw[fill=pink,opacity=0.5]
(x1) to[bend left=15] (w1)to[bend left=15] (f) to[bend left=2] cycle
;
\draw[fill=pink,opacity=0.5]
(x1) to[bend left=15] (z1)to[bend left=15] (f) to[bend left=15] cycle
;

\draw[fill=pink,opacity=1]
(z1) to[bend left=1] (w1)to[bend left=1] (f) to[bend right=1] cycle
;

\draw[fill=pink,opacity=0.5]
(y1) to[bend right=15] (w1)to[bend right = 15] (c) to[bend right=2] cycle
;

\draw[fill=pink,opacity=1]
(z1) to[bend right=1] (w1)to[bend left=1] (c) to[bend right=10] cycle
;

\draw (7.5,-2.2) node[vtx,label=below:$f$](w1){};

\draw[fill=pink,opacity=0.5]
(y1) to[bend right=15] (z1)to[bend right=35] (c) to[bend right=25] cycle
;

\end{scope}

\draw (7,-4) node{$F_1$};

\end{tikzpicture}

    \caption{The configuration when $be\in E(G_f)$, and the resulting copy of $F_1$.}
    \label{fig:nobd}
\end{figure}
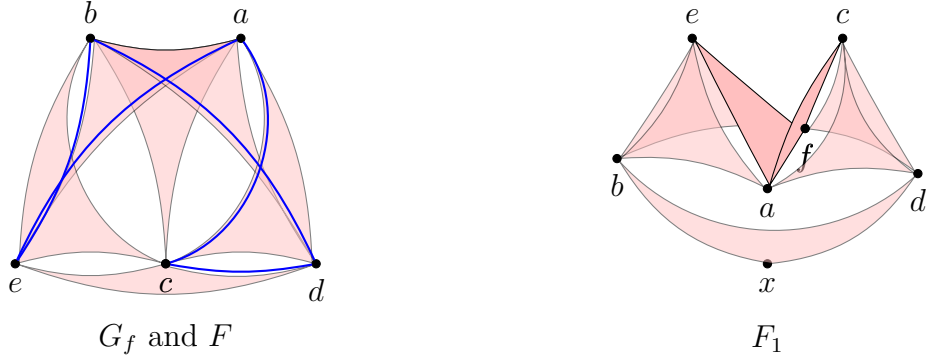


Thus, $bc \in E(G_f)$. By (viii), $cd \not\in E(G_f)$, so we must have $de \in E(G_f)$.
Hence $E(G_f) = \{ ae,ac,bc,bd,de \}$; see Figure~\ref{goodC5} for an illustration.

\begin{figure}[h]

\centering
\begin{tikzpicture}

\draw (0,0) node[vtx,label=above:$b$](b){};
\draw (-2,0) node[vtx,label=above:$a$](a){};
\draw (-1,-3) node[vtx,label=below:$c$](c){};
\draw (1,-3) node[vtx,label=below:$e$](e){};
\draw (-3,-3) node[vtx,label=below:$d$](d){};

\draw[fill=pink,opacity=0.5]
(a) to[bend right=15] (b)to[bend right=15] (c) to[bend right=15] cycle
;

\draw[fill=pink,opacity=0.5]
(a) to[bend right=15] (b)to[bend right=15] (d) to[bend right=15] cycle
;

\draw[fill=pink,opacity=0.5]
(a) to[bend right=15] (b)to[bend right=15] (e) to[bend right=15] cycle
;

\draw[fill=pink,opacity=0.5]
(d) to[bend left=15] (c)to[bend left=50] (a) to[bend right=15] cycle
;

\draw[fill=pink,opacity=0.5]
(e) to[bend right=15] (c)to[bend right=50] (b) to[bend left=15] cycle
;

\draw[fill=pink,opacity=0.5]
(c) to[bend left=15] (d)to[bend right=20] (e) to[bend left=15] cycle
;

\draw[thick, blue] (a) to[bend left = 20] (c);

\draw[thick, blue] (b) to[bend right = 20] (d);

\draw[thick, blue] (a) to[bend left = 20] (e);

\draw[thick, blue] (d) to[bend right = 20] (e);

\draw[thick, blue] (b) to[bend right = 20] (c);

\draw (0,0) node[vtx,label=above:$b$](b){};
\draw (-2,0) node[vtx,label=above:$a$](a){};
\draw (-1,-3) node[vtx,label=below:$c$](c){};
\draw (1,-3) node[vtx,label=below:$e$](e){};
\draw (-3,-3) node[vtx,label=below:$d$](d){};

\end{tikzpicture}
\caption{$F$ and $G_f$ when $ac, ae \in E(G_f)$.}\label{goodC5}

\end{figure}
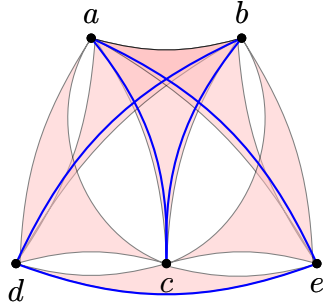

Unlike in the previous cases, we cannot immediately find \acopyof $F_1$ (or any other forbidden hypergraph) in Figure~\ref{goodC5}. 
However, we now have that the subhypergraph of $H$ induced on $\{a,b,c,d,e,f\}$ has edge set $\{abc, abd, abe, adc, bce, cde, fac, fcb, fbd, fde, fea  \}$. 
We call this subhypergraph $F'$, and we shall use $F'$ to find \acopyof some forbidden hypergraph. 

Observe first that each of the $15$ pairs of vertices in $F'$ have positive co-degree. Thus, there exists some $g \in V(H)$ that is in the neighborhood of at least $15\left( \frac{2}{5} + \varepsilon \right) > 6$ pairs. Let $G_g$ be the link graph of vertex $g$ induced on $V(F')$ (i.e.\ $G_g = L(g)[V(F')]$). 
Again, we begin with some observations on $G_g$. 

\begin{claim}\label{Gg}
$\delta(G_g) \geq 1$. Moreover, in $G_g$,
$N(f) = \{c,d,e\}$.
\end{claim}

\begin{proof}
As in the proof of Claim~\ref{VFclaims}, since all pairs from $\{a,b,c,d,e,f\}$ have positive co-degree, $G_g$ cannot contain \acopyof $K_4^-$.

 If $G_g$ has an isolated vertex $v$, then $G_g - v$ is a component on $5$ vertices and $7$ edges that necessarily contains a \acopyof $K_4^-$. Thus, every vertex of $G_g$ has positive degree.

Next, observe that the previous analysis of $G_f$ in fact shows that if some vertex $v \in V(H) \setminus V(F)$ has positive co-degree with at least $5$ pairs from $V(F)$, then the link graph of $v$ induced on $V(F)$ must be equal to $G_f$. In particular, either $G_g$ contains $G_f$ or $g$ has positive co-degree with at most $4$ pairs from $V(F)$. Observe that if $G_g$ contains $G_f$, then $\{a,b,c,f\}$ and $\{a,b,c,g\}$ span copies of $K_4^{3-}$ with spike vertex $c$. This implies that $\{a,b,c,f,g\}$ spans \acopyof $F_2$, since we have argued that $d(f,g) > 0$. Thus, $g$ has positive co-degree with at most $4$ pairs from $V(F)$, which implies that $f$ has degree at least $3$ in $G_g$, since $|E(G_g)| \geq 7$.

Finally, to determine $N(f)$ in $G_g$, we consider the interaction between $G_g$ and $G_f$. Suppose that $x,y$ are neighbors of $f$ in $G_g$ such that $xy \in E(G_f)$. Then $\{x,y,f,g\}$ forms \acopyof $K_4^{3-}$ with spike vertex $f$. We consider the possible values of $x,y$; we know $xy \in \{ac,ae, bc, bd, de\}$. Since $\{a,b,c,d\}$ spans \acopyof $K_4^{3-}$ with spike vertex $a$, and we have established that $g$ has positive co-degree with every vertex in $V(F)$, we can find \acopyof $F_1$ in $H$ on vertex set $\{a,b,c,d,f,g\}$ if $xy \in \{bc, bd\}$. Similarly, if $xy \in \{ac, ae\}$, then we can find \acopyof $F_1$ on $\{a,b,c,e,f,g\}$. Thus, in $G_g$, either $N(f)$ is an independent set or $N(f)$ contains precisely the edge $de$. Recall that $|N(f)| \geq 3$, so the first outcome is impossible, as the independence number of $G_f$ is $2$. The second outcome occurs only if 
$N(f) = \{c,d,e\}$.
\end{proof}

With Claim~\ref{Gg} established, we conclude that $|E(G_g)| = 7$, with exactly 4 edges of $G_g$ spanned by $\{a,b,c,d,e\}$. Given that $cf, df, ef  \in E(G_g)$, we shall prove that this is not possible.

First, observe that none of $cd,ce,de$ are in $E(G_g)$. Indeed, if $de \in E(G_g)$, then $\{d,e,f,g\}$ spans \acopyof $K_4^3$ (recall that $de \in E(G_f)$). If either $cd$ or $ce$ is in $E(G_g)$, then either $\{c,d,f,g\}$ or $\{c,e,f,g\}$ spans \acopyof $K_4^{3-}$ with spike vertex $g$; if $\{c,d,f,g\}$ spans $K_4^{3-}$, then $\{a,b,c,d,f,g\}$ spans an \acopyof $F_1$, while if $\{c,e,f,g\}$ spans $K_4^{3-}$, then $\{a,b,c,e,f,g\}$ spans an \acopyof $F_1$.

Notice that Claim~\ref{VFclaims} applies also to $G_g$ since Claim~\ref{Gg} proves (ii) and  it is the only statement that used the number of edges of $G_f$.
Next, we prove that $ab \not\in E(G_g)$. Indeed, by (v) and (vi), 
if $ab \in E(G_g)$, then $ac,bc,ad,be$ are not in $E(G_g)$. 
Recall $cd,ce,de$ are not in $E(G_g)$.
So if $ab \in E(G_g)$, then at most three edges of $G_g$ are spanned by $\{a,b,c,d,e\}$. However, by Claim~\ref{Gg}, $d(f) = 3$ in $G_g$, so we have $|E(G_g)| \leq 6$, a contradiction. Thus, $ab \not\in E(G_g)$.

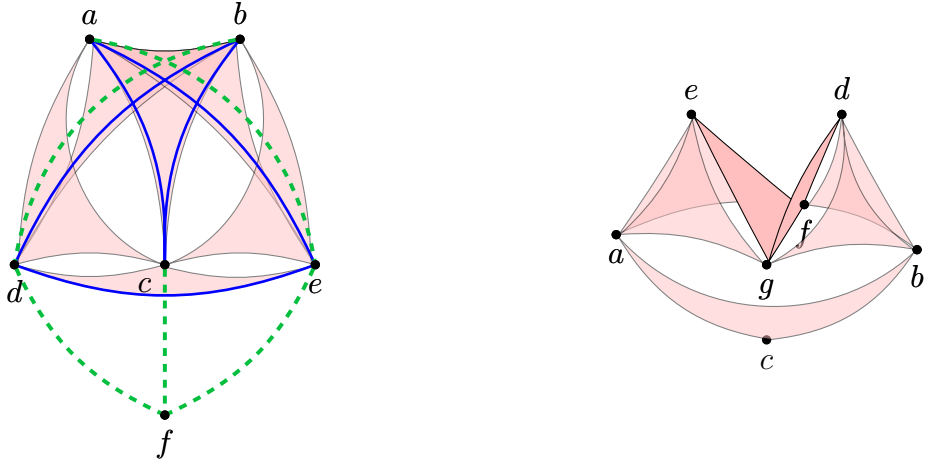
\begin{figure}[h!]

\centering
\begin{tikzpicture}

\draw (0,0) node[vtx,label=above:$b$](b){};
\draw (-2,0) node[vtx,label=above:$a$](a){};
\draw (-1,-3) node[vtx,label=below left:$c$](c){};
\draw (1,-3) node[vtx,label=below:$e$](e){};
\draw (-3,-3) node[vtx,label=below:$d$](d){};
\draw (-1,-5) node[vtx,label=below:$f$](f){};

\draw[fill=pink,opacity=0.5]
(a) to[bend right=15] (b)to[bend right=15] (c) to[bend right=15] cycle
;

\draw[fill=pink,opacity=0.5]
(a) to[bend right=15] (b)to[bend right=15] (d) to[bend right=15] cycle
;

\draw[fill=pink,opacity=0.5]
(a) to[bend right=15] (b)to[bend right=15] (e) to[bend right=15] cycle
;

\draw[fill=pink,opacity=0.5]
(d) to[bend left=15] (c)to[bend left=50] (a) to[bend right=15] cycle
;

\draw[fill=pink,opacity=0.5]
(e) to[bend right=15] (c)to[bend right=50] (b) to[bend left=15] cycle
;

\draw[fill=pink,opacity=0.5]
(c) to[bend left=15] (d)to[bend right=20] (e) to[bend left=15] cycle
;

\draw[line width=1pt, blue] 
 (a) to[bend left = 20] (c)
 (b) to[bend right = 20] (d)
 (a) to[bend left = 20] (e)
 (d) to[bend right = 20] (e)
 (b) to[bend right = 20] (c);

\draw[line width=1.5pt, darkpastelgreen, dashed] 
(d) to[bend right = 20] (f)
(e) to[bend left = 20] (f)
(c) to (f)
(b) to[bend right = 35] (d)
(a) to[bend left = 35] (e);

\draw (0,0) node[vtx,label=above:$b$](b){};
\draw (-2,0) node[vtx,label=above:$a$](a){};
\draw (-1,-3) node[vtx,label=below left:$c$](c){};
\draw (1,-3) node[vtx,label=below:$e$](e){};
\draw (-3,-3) node[vtx,label=below:$d$](d){};

\draw (-1,-5) node[vtx,label=below:$f$](f){};

\draw (5,-2.6) node[vtx,label=below:$a$](x1){};
\draw (7,-3) node[vtx,label=below:$g$](z1){};
\draw (7.5,-2.2) node[vtx,label=below:$f$](w1){};
\draw (9, -2.8) node[vtx,label=below:$b$](y1){};    
\draw (6, -1) node[vtx,label=above:$e$](s1){}; 
\draw (8, -1) node[vtx,label=above:$d$](t1){}; 

\draw (7, -4) node[vtx,label=below:$c$](g1){}; 
\draw[fill=pink,opacity=0.5]
(x1) to[bend right=45] (y1) to[bend left=25] (g1) to[bend left=25] cycle
;

\draw[fill=pink,opacity=0.5]
(x1) to[bend left=15] (w1)to[bend left=15] (s1) to[bend left=2] cycle
;
\draw[fill=pink,opacity=0.5]
(x1) to[bend left=15] (z1)to[bend left=15] (s1) to[bend left=15] cycle
;

\draw[fill=pink,opacity=1]
(z1) to[bend left=1] (w1)to[bend left=1] (s1) to[bend right=1] cycle
;

\draw[fill=pink,opacity=0.5]
(y1) to[bend right=15] (w1)to[bend right = 15] (t1) to[bend right=2] cycle
;

\draw[fill=pink,opacity=1]
(z1) to[bend right=1] (w1)to[bend left=1] (t1) to[bend right=10] cycle
;

\draw (7.5,-2.2) node[vtx,label=below:$f$](w1){};

\draw[fill=pink,opacity=0.5]
(y1) to[bend right=15] (z1)to[bend right=35] (t1) to[bend right=25] cycle
;

\draw (5,-2.6) node[vtx,label=below:$a$](x1){};
\draw (7,-3) node[vtx,label=below:$g$](z1){};
\draw (9, -2.8) node[vtx,label=below:$b$](y1){};    
\draw (6, -1) node[vtx,label=above:$e$](s1){}; 
\draw (8, -1) node[vtx,label=above:$d$](t1){};


\end{tikzpicture}
\caption{$F'$ (edges involving $f$ indicated by $G_f$ in thick blue) and a configuration in $G_g$ (in dashed green) yielding $F_1$.}\label{finalconfig}
\end{figure}

Finally, by (vii), at most two of $ac, ad, bc, be$ are in $E(G_g)$. So in order to have $|E(G_g)| = 7$, both $ae$ and $bd$ must be in $E(G_g)$. However, this results in \acopyof $F_1$ on $\{a,b,d,e,g,f\}$ with edges $\{age, afe, abc, gfe, gfd, gbd, bdf  \}$ as depicted in Figure~\ref{finalconfig}. 
\end{proof}

\section{Positive co-degree densities from flag algebras}\label{sec:J4}

This section contains calculations using flag algebras, as introduced by Razborov~\cite{Raz07}. 
For an introduction to flag algebras and formal definitions, see~\cite{FalgasK4-,Raz07,FALGAS-RAVRY_VAUGHAN_2013}. 
Computer code is available at \oururl.

\begin{figure}[h!]
\begin{center}
\begin{tikzpicture}
\draw[scale=5,line width=6pt]
(0,0)  coordinate (x7) -- ++(1,0)  coordinate[pos=0.5](x1) coordinate(x2)
(x7) -- ++(60:1)  coordinate[pos=0.5](x3) coordinate(x4)
(x2) -- ++(120:0.5) coordinate(x5)
(x7) -- ++(30:0.57735) coordinate(x6) 
(x4)--(x5) 
(x2)--(x3)
(x6)--(x5)
(x4)--(x1)
(x6) circle (0.2886751345 cm)  
;
\draw
\foreach \x/\y in {1/1,2/2,3/6,4/7,5/3,6/4,7/5}{
(x\x) node[circle,draw,fill=white]{$X_\y$}
}
;
\end{tikzpicture}
\end{center}
    \caption{A blow-up of the Fano plane. In  $\overline{\mathbb{F}}_n$, $3$-edges span triples of classes which do not span $3$-edges in the blow-up of $\mathbb{F}$.} 
    \label{fig:Fano}
\end{figure}
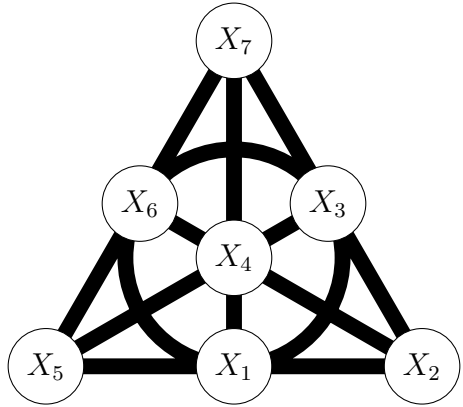

We use flag algebras to prove Theorem~\ref{J4} that $\gamma^+(J_4) = 4/7$, with the asymptotically unique construction being the balanced blow-up of the complement of the Fano plane. In Figure~\ref{fig:Fano}, we illustrate (a blow-up of) the Fano plane, with a particular labeling of classes to which we will later refer. We denote the complement of the Fano plane by $\overline{\mathbb{F}}$ and the $n$-vertex, balanced blow-up of the complement of the Fano plane by $\overline{\mathbb{F}}_n$.

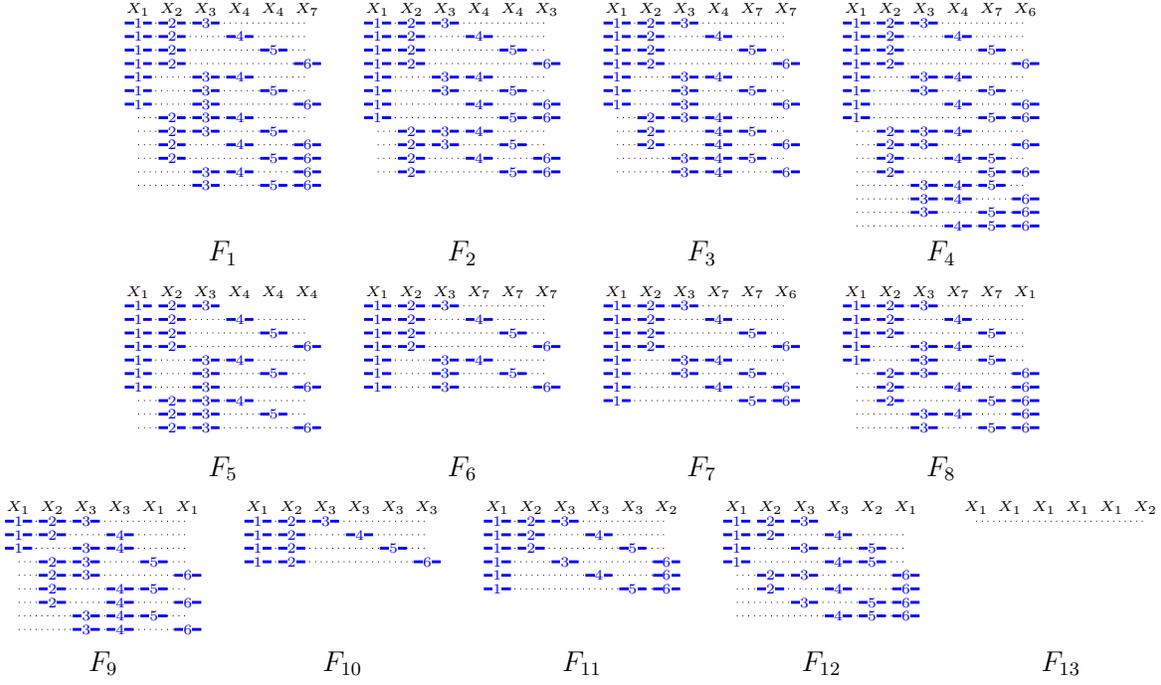
\begin{figure}
\begin{center}

\def\xshift{1.25}
\def\yshift{-3.6}
\def\yshiftB{-2.6}
\def\yshiftC{-2.3}
\def\yshiftD{-2.3}

\def\ourscale{0.9}

\scalebox{\ourscale}{
 \vc{ 
  \begin{tikzpicture}[flag_pic]\outercycle{6}{0}
      \draw(\xshift,\yshift) node{$F_1$};
\draw
(x0) node{\tiny$X_1$}
(x1) node{\tiny$X_2$}
(x2) node{\tiny$X_3$}
(x3) node{\tiny$X_4$}
(x4) node{\tiny$X_4$}
(x5) node{\tiny$X_7$}
;
\drawhyperedge{0}{6}
\drawhypervertex{0}{0}
\drawhypervertex{1}{0}
\drawhypervertex{2}{0}
\drawhyperedge{1}{6}
\drawhypervertex{0}{1}
\drawhypervertex{1}{1}
\drawhypervertex{3}{1}
\drawhyperedge{2}{6}
\drawhypervertex{0}{2}
\drawhypervertex{1}{2}
\drawhypervertex{4}{2}
\drawhyperedge{3}{6}
\drawhypervertex{0}{3}
\drawhypervertex{1}{3}
\drawhypervertex{5}{3}
\drawhyperedge{4}{6}
\drawhypervertex{0}{4}
\drawhypervertex{2}{4}
\drawhypervertex{3}{4}
\drawhyperedge{5}{6}
\drawhypervertex{0}{5}
\drawhypervertex{2}{5}
\drawhypervertex{4}{5}
\drawhyperedge{6}{6}
\drawhypervertex{0}{6}
\drawhypervertex{2}{6}
\drawhypervertex{5}{6}
\drawhyperedge{7}{6}
\drawhypervertex{1}{7}
\drawhypervertex{2}{7}
\drawhypervertex{3}{7}
\drawhyperedge{8}{6}
\drawhypervertex{1}{8}
\drawhypervertex{2}{8}
\drawhypervertex{4}{8}
\drawhyperedge{9}{6}
\drawhypervertex{1}{9}
\drawhypervertex{3}{9}
\drawhypervertex{5}{9}
\drawhyperedge{10}{6}
\drawhypervertex{1}{10}
\drawhypervertex{4}{10}
\drawhypervertex{5}{10}
\drawhyperedge{11}{6}
\drawhypervertex{2}{11}
\drawhypervertex{3}{11}
\drawhypervertex{5}{11}
\drawhyperedge{12}{6}
\drawhypervertex{2}{12}
\drawhypervertex{4}{12}
\drawhypervertex{5}{12}
\end{tikzpicture} } 
 \vc{ 
  \begin{tikzpicture}[flag_pic]\outercycle{6}{0}
      \draw(\xshift,\yshift) node{$F_2$};
\draw
(x0) node{\tiny$X_1$}
(x1) node{\tiny$X_2$}
(x2) node{\tiny$X_3$}
(x3) node{\tiny$X_4$}
(x4) node{\tiny$X_4$}
(x5) node{\tiny$X_3$}
;
\drawhyperedge{0}{6}
\drawhypervertex{0}{0}
\drawhypervertex{1}{0}
\drawhypervertex{2}{0}
\drawhyperedge{1}{6}
\drawhypervertex{0}{1}
\drawhypervertex{1}{1}
\drawhypervertex{3}{1}
\drawhyperedge{2}{6}
\drawhypervertex{0}{2}
\drawhypervertex{1}{2}
\drawhypervertex{4}{2}
\drawhyperedge{3}{6}
\drawhypervertex{0}{3}
\drawhypervertex{1}{3}
\drawhypervertex{5}{3}
\drawhyperedge{4}{6}
\drawhypervertex{0}{4}
\drawhypervertex{2}{4}
\drawhypervertex{3}{4}
\drawhyperedge{5}{6}
\drawhypervertex{0}{5}
\drawhypervertex{2}{5}
\drawhypervertex{4}{5}
\drawhyperedge{6}{6}
\drawhypervertex{0}{6}
\drawhypervertex{3}{6}
\drawhypervertex{5}{6}
\drawhyperedge{7}{6}
\drawhypervertex{0}{7}
\drawhypervertex{4}{7}
\drawhypervertex{5}{7}
\drawhyperedge{8}{6}
\drawhypervertex{1}{8}
\drawhypervertex{2}{8}
\drawhypervertex{3}{8}
\drawhyperedge{9}{6}
\drawhypervertex{1}{9}
\drawhypervertex{2}{9}
\drawhypervertex{4}{9}
\drawhyperedge{10}{6}
\drawhypervertex{1}{10}
\drawhypervertex{3}{10}
\drawhypervertex{5}{10}
\drawhyperedge{11}{6}
\drawhypervertex{1}{11}
\drawhypervertex{4}{11}
\drawhypervertex{5}{11}
\end{tikzpicture} } 
 \vc{ 
  \begin{tikzpicture}[flag_pic]\outercycle{6}{0}
      \draw(\xshift,\yshift) node{$F_3$};
\draw
(x0) node{\tiny$X_1$}
(x1) node{\tiny$X_2$}
(x2) node{\tiny$X_3$}
(x3) node{\tiny$X_4$}
(x4) node{\tiny$X_7$}
(x5) node{\tiny$X_7$}
;
\drawhyperedge{0}{6}
\drawhypervertex{0}{0}
\drawhypervertex{1}{0}
\drawhypervertex{2}{0}
\drawhyperedge{1}{6}
\drawhypervertex{0}{1}
\drawhypervertex{1}{1}
\drawhypervertex{3}{1}
\drawhyperedge{2}{6}
\drawhypervertex{0}{2}
\drawhypervertex{1}{2}
\drawhypervertex{4}{2}
\drawhyperedge{3}{6}
\drawhypervertex{0}{3}
\drawhypervertex{1}{3}
\drawhypervertex{5}{3}
\drawhyperedge{4}{6}
\drawhypervertex{0}{4}
\drawhypervertex{2}{4}
\drawhypervertex{3}{4}
\drawhyperedge{5}{6}
\drawhypervertex{0}{5}
\drawhypervertex{2}{5}
\drawhypervertex{4}{5}
\drawhyperedge{6}{6}
\drawhypervertex{0}{6}
\drawhypervertex{2}{6}
\drawhypervertex{5}{6}
\drawhyperedge{7}{6}
\drawhypervertex{1}{7}
\drawhypervertex{2}{7}
\drawhypervertex{3}{7}
\drawhyperedge{8}{6}
\drawhypervertex{1}{8}
\drawhypervertex{3}{8}
\drawhypervertex{4}{8}
\drawhyperedge{9}{6}
\drawhypervertex{1}{9}
\drawhypervertex{3}{9}
\drawhypervertex{5}{9}
\drawhyperedge{10}{6}
\drawhypervertex{2}{10}
\drawhypervertex{3}{10}
\drawhypervertex{4}{10}
\drawhyperedge{11}{6}
\drawhypervertex{2}{11}
\drawhypervertex{3}{11}
\drawhypervertex{5}{11}
\end{tikzpicture} } 
 \vc{ 
  \begin{tikzpicture}[flag_pic]\outercycle{6}{0}
    \draw(\xshift,\yshift) node{$F_4$};
\draw
(x0) node{\tiny$X_1$}
(x1) node{\tiny$X_2$}
(x2) node{\tiny$X_3$}
(x3) node{\tiny$X_4$}
(x4) node{\tiny$X_7$}
(x5) node{\tiny$X_6$}
;
\drawhyperedge{0}{6}
\drawhypervertex{0}{0}
\drawhypervertex{1}{0}
\drawhypervertex{2}{0}
\drawhyperedge{1}{6}
\drawhypervertex{0}{1}
\drawhypervertex{1}{1}
\drawhypervertex{3}{1}
\drawhyperedge{2}{6}
\drawhypervertex{0}{2}
\drawhypervertex{1}{2}
\drawhypervertex{4}{2}
\drawhyperedge{3}{6}
\drawhypervertex{0}{3}
\drawhypervertex{1}{3}
\drawhypervertex{5}{3}
\drawhyperedge{4}{6}
\drawhypervertex{0}{4}
\drawhypervertex{2}{4}
\drawhypervertex{3}{4}
\drawhyperedge{5}{6}
\drawhypervertex{0}{5}
\drawhypervertex{2}{5}
\drawhypervertex{4}{5}
\drawhyperedge{6}{6}
\drawhypervertex{0}{6}
\drawhypervertex{3}{6}
\drawhypervertex{5}{6}
\drawhyperedge{7}{6}
\drawhypervertex{0}{7}
\drawhypervertex{4}{7}
\drawhypervertex{5}{7}
\drawhyperedge{8}{6}
\drawhypervertex{1}{8}
\drawhypervertex{2}{8}
\drawhypervertex{3}{8}
\drawhyperedge{9}{6}
\drawhypervertex{1}{9}
\drawhypervertex{2}{9}
\drawhypervertex{5}{9}
\drawhyperedge{10}{6}
\drawhypervertex{1}{10}
\drawhypervertex{3}{10}
\drawhypervertex{4}{10}
\drawhyperedge{11}{6}
\drawhypervertex{1}{11}
\drawhypervertex{4}{11}
\drawhypervertex{5}{11}
\drawhyperedge{12}{6}
\drawhypervertex{2}{12}
\drawhypervertex{3}{12}
\drawhypervertex{4}{12}
\drawhyperedge{13}{6}
\drawhypervertex{2}{13}
\drawhypervertex{3}{13}
\drawhypervertex{5}{13}
\drawhyperedge{14}{6}
\drawhypervertex{2}{14}
\drawhypervertex{4}{14}
\drawhypervertex{5}{14}
\drawhyperedge{15}{6}
\drawhypervertex{3}{15}
\drawhypervertex{4}{15}
\drawhypervertex{5}{15}
\end{tikzpicture} } 
}

\scalebox{\ourscale}{
 \vc{ 
  \begin{tikzpicture}[flag_pic]
    \draw(\xshift,\yshiftB) node{$F_5$};
    \outercycle{6}{0}
\draw
(x0) node{\tiny$X_1$}
(x1) node{\tiny$X_2$}
(x2) node{\tiny$X_3$}
(x3) node{\tiny$X_4$}
(x4) node{\tiny$X_4$}
(x5) node{\tiny$X_4$}
;
\drawhypervertex{0}{0}
\drawhypervertex{1}{0}
\drawhypervertex{2}{0}
\drawhyperedge{1}{6}
\drawhypervertex{0}{1}
\drawhypervertex{1}{1}
\drawhypervertex{3}{1}
\drawhyperedge{2}{6}
\drawhypervertex{0}{2}
\drawhypervertex{1}{2}
\drawhypervertex{4}{2}
\drawhyperedge{3}{6}
\drawhypervertex{0}{3}
\drawhypervertex{1}{3}
\drawhypervertex{5}{3}
\drawhyperedge{4}{6}
\drawhypervertex{0}{4}
\drawhypervertex{2}{4}
\drawhypervertex{3}{4}
\drawhyperedge{5}{6}
\drawhypervertex{0}{5}
\drawhypervertex{2}{5}
\drawhypervertex{4}{5}
\drawhyperedge{6}{6}
\drawhypervertex{0}{6}
\drawhypervertex{2}{6}
\drawhypervertex{5}{6}
\drawhyperedge{7}{6}
\drawhypervertex{1}{7}
\drawhypervertex{2}{7}
\drawhypervertex{3}{7}
\drawhyperedge{8}{6}
\drawhypervertex{1}{8}
\drawhypervertex{2}{8}
\drawhypervertex{4}{8}
\drawhyperedge{9}{6}
\drawhypervertex{1}{9}
\drawhypervertex{2}{9}
\drawhypervertex{5}{9}
\end{tikzpicture} } 
 \vc{ 
  \begin{tikzpicture}[flag_pic]\outercycle{6}{0}
      \draw(\xshift,\yshiftB) node{$F_6$};
\draw
(x0) node{\tiny$X_1$}
(x1) node{\tiny$X_2$}
(x2) node{\tiny$X_3$}
(x3) node{\tiny$X_7$}
(x4) node{\tiny$X_7$}
(x5) node{\tiny$X_7$}
;
\drawhyperedge{0}{6}
\drawhypervertex{0}{0}
\drawhypervertex{1}{0}
\drawhypervertex{2}{0}
\drawhyperedge{1}{6}
\drawhypervertex{0}{1}
\drawhypervertex{1}{1}
\drawhypervertex{3}{1}
\drawhyperedge{2}{6}
\drawhypervertex{0}{2}
\drawhypervertex{1}{2}
\drawhypervertex{4}{2}
\drawhyperedge{3}{6}
\drawhypervertex{0}{3}
\drawhypervertex{1}{3}
\drawhypervertex{5}{3}
\drawhyperedge{4}{6}
\drawhypervertex{0}{4}
\drawhypervertex{2}{4}
\drawhypervertex{3}{4}
\drawhyperedge{5}{6}
\drawhypervertex{0}{5}
\drawhypervertex{2}{5}
\drawhypervertex{4}{5}
\drawhyperedge{6}{6}
\drawhypervertex{0}{6}
\drawhypervertex{2}{6}
\drawhypervertex{5}{6}
\end{tikzpicture} } 
 \vc{ 
  \begin{tikzpicture}[flag_pic]\outercycle{6}{0}
      \draw(\xshift,\yshiftB) node{$F_7$};
\draw
(x0) node{\tiny$X_1$}
(x1) node{\tiny$X_2$}
(x2) node{\tiny$X_3$}
(x3) node{\tiny$X_7$}
(x4) node{\tiny$X_7$}
(x5) node{\tiny$X_6$}
;
\drawhyperedge{0}{6}
\drawhypervertex{0}{0}
\drawhypervertex{1}{0}
\drawhypervertex{2}{0}
\drawhyperedge{1}{6}
\drawhypervertex{0}{1}
\drawhypervertex{1}{1}
\drawhypervertex{3}{1}
\drawhyperedge{2}{6}
\drawhypervertex{0}{2}
\drawhypervertex{1}{2}
\drawhypervertex{4}{2}
\drawhyperedge{3}{6}
\drawhypervertex{0}{3}
\drawhypervertex{1}{3}
\drawhypervertex{5}{3}
\drawhyperedge{4}{6}
\drawhypervertex{0}{4}
\drawhypervertex{2}{4}
\drawhypervertex{3}{4}
\drawhyperedge{5}{6}
\drawhypervertex{0}{5}
\drawhypervertex{2}{5}
\drawhypervertex{4}{5}
\drawhyperedge{6}{6}
\drawhypervertex{0}{6}
\drawhypervertex{3}{6}
\drawhypervertex{5}{6}
\drawhyperedge{7}{6}
\drawhypervertex{0}{7}
\drawhypervertex{4}{7}
\drawhypervertex{5}{7}
\end{tikzpicture} } 
 \vc{ 
  \begin{tikzpicture}[flag_pic]\outercycle{6}{0}
      \draw(\xshift,\yshiftB) node{$F_8$};
\draw
(x0) node{\tiny$X_1$}
(x1) node{\tiny$X_2$}
(x2) node{\tiny$X_3$}
(x3) node{\tiny$X_7$}
(x4) node{\tiny$X_7$}
(x5) node{\tiny$X_1$}
;
\drawhyperedge{0}{6}
\drawhypervertex{0}{0}
\drawhypervertex{1}{0}
\drawhypervertex{2}{0}
\drawhyperedge{1}{6}
\drawhypervertex{0}{1}
\drawhypervertex{1}{1}
\drawhypervertex{3}{1}
\drawhyperedge{2}{6}
\drawhypervertex{0}{2}
\drawhypervertex{1}{2}
\drawhypervertex{4}{2}
\drawhyperedge{3}{6}
\drawhypervertex{0}{3}
\drawhypervertex{2}{3}
\drawhypervertex{3}{3}
\drawhyperedge{4}{6}
\drawhypervertex{0}{4}
\drawhypervertex{2}{4}
\drawhypervertex{4}{4}
\drawhyperedge{5}{6}
\drawhypervertex{1}{5}
\drawhypervertex{2}{5}
\drawhypervertex{5}{5}
\drawhyperedge{6}{6}
\drawhypervertex{1}{6}
\drawhypervertex{3}{6}
\drawhypervertex{5}{6}
\drawhyperedge{7}{6}
\drawhypervertex{1}{7}
\drawhypervertex{4}{7}
\drawhypervertex{5}{7}
\drawhyperedge{8}{6}
\drawhypervertex{2}{8}
\drawhypervertex{3}{8}
\drawhypervertex{5}{8}
\drawhyperedge{9}{6}
\drawhypervertex{2}{9}
\drawhypervertex{4}{9}
\drawhypervertex{5}{9}
\end{tikzpicture} } 
}

\scalebox{\ourscale}{
 \vc{ 
  \begin{tikzpicture}[flag_pic]\outercycle{6}{0}
      \draw(\xshift,\yshiftC) node{$F_9$};
\draw
(x0) node{\tiny$X_1$}
(x1) node{\tiny$X_2$}
(x2) node{\tiny$X_3$}
(x3) node{\tiny$X_3$}
(x4) node{\tiny$X_1$}
(x5) node{\tiny$X_1$}
;
\drawhyperedge{0}{6}
\drawhypervertex{0}{0}
\drawhypervertex{1}{0}
\drawhypervertex{2}{0}
\drawhyperedge{1}{6}
\drawhypervertex{0}{1}
\drawhypervertex{1}{1}
\drawhypervertex{3}{1}
\drawhyperedge{2}{6}
\drawhypervertex{0}{2}
\drawhypervertex{2}{2}
\drawhypervertex{3}{2}
\drawhyperedge{3}{6}
\drawhypervertex{1}{3}
\drawhypervertex{2}{3}
\drawhypervertex{4}{3}
\drawhyperedge{4}{6}
\drawhypervertex{1}{4}
\drawhypervertex{2}{4}
\drawhypervertex{5}{4}
\drawhyperedge{5}{6}
\drawhypervertex{1}{5}
\drawhypervertex{3}{5}
\drawhypervertex{4}{5}
\drawhyperedge{6}{6}
\drawhypervertex{1}{6}
\drawhypervertex{3}{6}
\drawhypervertex{5}{6}
\drawhyperedge{7}{6}
\drawhypervertex{2}{7}
\drawhypervertex{3}{7}
\drawhypervertex{4}{7}
\drawhyperedge{8}{6}
\drawhypervertex{2}{8}
\drawhypervertex{3}{8}
\drawhypervertex{5}{8}
\end{tikzpicture} } 
 \vc{ 
  \begin{tikzpicture}[flag_pic]\outercycle{6}{0}
        \draw(\xshift,\yshiftC) node{$F_{10}$};
\draw
(x0) node{\tiny$X_1$}
(x1) node{\tiny$X_2$}
(x2) node{\tiny$X_3$}
(x3) node{\tiny$X_3$}
(x4) node{\tiny$X_3$}
(x5) node{\tiny$X_3$}
;
\drawhyperedge{0}{6}
\drawhypervertex{0}{0}
\drawhypervertex{1}{0}
\drawhypervertex{2}{0}
\drawhyperedge{1}{6}
\drawhypervertex{0}{1}
\drawhypervertex{1}{1}
\drawhypervertex{3}{1}
\drawhyperedge{2}{6}
\drawhypervertex{0}{2}
\drawhypervertex{1}{2}
\drawhypervertex{4}{2}
\drawhyperedge{3}{6}
\drawhypervertex{0}{3}
\drawhypervertex{1}{3}
\drawhypervertex{5}{3}
\end{tikzpicture} } 
 \vc{ 
  \begin{tikzpicture}[flag_pic]\outercycle{6}{0}
      \draw(\xshift,\yshiftC) node{$F_{11}$};
\draw
(x0) node{\tiny$X_1$}
(x1) node{\tiny$X_2$}
(x2) node{\tiny$X_3$}
(x3) node{\tiny$X_3$}
(x4) node{\tiny$X_3$}
(x5) node{\tiny$X_2$}
;
\drawhyperedge{0}{6}
\drawhypervertex{0}{0}
\drawhypervertex{1}{0}
\drawhypervertex{2}{0}
\drawhyperedge{1}{6}
\drawhypervertex{0}{1}
\drawhypervertex{1}{1}
\drawhypervertex{3}{1}
\drawhyperedge{2}{6}
\drawhypervertex{0}{2}
\drawhypervertex{1}{2}
\drawhypervertex{4}{2}
\drawhyperedge{3}{6}
\drawhypervertex{0}{3}
\drawhypervertex{2}{3}
\drawhypervertex{5}{3}
\drawhyperedge{4}{6}
\drawhypervertex{0}{4}
\drawhypervertex{3}{4}
\drawhypervertex{5}{4}
\drawhyperedge{5}{6}
\drawhypervertex{0}{5}
\drawhypervertex{4}{5}
\drawhypervertex{5}{5}
\end{tikzpicture} } 
 \vc{ 
  \begin{tikzpicture}[flag_pic]\outercycle{6}{0}
      \draw(\xshift,\yshiftC) node{$F_{12}$};
\draw
(x0) node{\tiny$X_1$}
(x1) node{\tiny$X_2$}
(x2) node{\tiny$X_3$}
(x3) node{\tiny$X_3$}
(x4) node{\tiny$X_2$}
(x5) node{\tiny$X_1$}
;
\drawhypervertex{0}{0}
\drawhypervertex{1}{0}
\drawhypervertex{2}{0}
\drawhyperedge{1}{6}
\drawhypervertex{0}{1}
\drawhypervertex{1}{1}
\drawhypervertex{3}{1}
\drawhyperedge{2}{6}
\drawhypervertex{0}{2}
\drawhypervertex{2}{2}
\drawhypervertex{4}{2}
\drawhyperedge{3}{6}
\drawhypervertex{0}{3}
\drawhypervertex{3}{3}
\drawhypervertex{4}{3}
\drawhyperedge{4}{6}
\drawhypervertex{1}{4}
\drawhypervertex{2}{4}
\drawhypervertex{5}{4}
\drawhyperedge{5}{6}
\drawhypervertex{1}{5}
\drawhypervertex{3}{5}
\drawhypervertex{5}{5}
\drawhyperedge{6}{6}
\drawhypervertex{2}{6}
\drawhypervertex{4}{6}
\drawhypervertex{5}{6}
\drawhyperedge{7}{6}
\drawhypervertex{3}{7}
\drawhypervertex{4}{7}
\drawhypervertex{5}{7}
\end{tikzpicture} } 
 \vc{ 
  \begin{tikzpicture}[flag_pic]\outercycle{6}{0}
      \draw(\xshift,\yshiftD) node{$F_{13}$};
  \drawhyperedge{0}{6}
\draw
(x0) node{\tiny$X_1$}
(x1) node{\tiny$X_1$}
(x2) node{\tiny$X_1$}
(x3) node{\tiny$X_1$}
(x4) node{\tiny$X_1$}
(x5) node{\tiny$X_2$}
;
\end{tikzpicture} }
}
\end{center}
\caption{The family $\mathcal{F}$ of thirteen $6$-vertex $3$-graphs.}\label{J4 positive density set}
\end{figure}

Let $\mathcal{F}$ be the family of thirteen 3-graphs 6-vertex induced subgraphs of $\overline{\mathbb{F}}_n$, depicted in Figure~\ref{J4 positive density set}. We include labels to indicate which classes of $\overline{\mathbb{F}}_n$ each vertex belongs to, corresponding to the labeling in Figure~\ref{fig:Fano}. Note that due to the symmetries of $\overline{\mathbb{F}}_n$, the indicated labelings are not the unique ways a subgraph can be obtained.

Let $\mathcal{F}_6$ be the family of $J_4$-free $6$-vertex $3$-graphs. 

The following claim shows that if a $J_4$-free $n$-vertex $3$-graph $G$ has minimum positive co-degree at least $\frac{4}{7}n$ then every $3$-graph on six vertices with positive density in $G$ is in $\mathcal{F}$.

\begin{claim}\label{claim:fa47}
    For every fixed $\delta > 0$, there exists $n_0$ such that for every $n \geq n_0$, if $G_n$ is a $J_4$-free $n$-vertex 3-graph with $\delta_2^+(G_n) \geq \frac{4n}{7}$, then
\[
\sum_{H \in \mathcal{F}_6 \setminus \mathcal{F}} d(H,G_n) \leq \delta.
\]    
\end{claim}
\begin{proof}
The claim is proved by a standard application of flag algebras. 
The notable part is encoding the condition $\delta_2^+(G_n) \geq \frac{4n}{7}$ into flag algebras by \eqref{eq:f47}.
 \begin{align}
 0 \leq
    \vc{ 
  \begin{tikzpicture}[flag_pic]\outercycle{3}{2}
\draw (x0) node[labeled_vertex]{};\draw (x1) node[labeled_vertex]{};\draw (x2) node[unlabeled_vertex]{};
\labelvertex{0}\labelvertex{1}
\drawhyperedge{0}{3}
\drawhypervertex{0}{0}
\drawhypervertex{1}{0}
\drawhypervertex{2}{0}
\end{tikzpicture} }   \times  \Bigg(  7  \vc{ 
  \begin{tikzpicture}[flag_pic]\outercycle{3}{2}
\draw (x0) node[labeled_vertex]{};\draw (x1) node[labeled_vertex]{};\draw (x2) node[unlabeled_vertex]{};
\labelvertex{0}\labelvertex{1}\drawhyperedge{0}{3}
\drawhypervertex{0}{0}
\drawhypervertex{1}{0}
\drawhypervertex{2}{0}
\end{tikzpicture} }  - 4  \vc{ 
  \begin{tikzpicture}[flag_pic]\outercycle{2}{2}
\draw (x0) node[labeled_vertex]{};\draw (x1) node[labeled_vertex]{};
\labelvertex{0}\labelvertex{1}\end{tikzpicture} }  \Bigg). 
\label{eq:f47}
\end{align}
If two labeled vertices 1 and 2 have zero co-degree, the right-hand side of \eqref{eq:f47} is zero because of the first term. If they have positive co-degree, both terms in the product are non-negative. 

The calculation is computer-assisted and too large to fit in this paper. The details are available at \oururl.
\end{proof}

In the proof of Theorem~\ref{J4} we will use an induced removal lemma.

\begin{theorem}[Induced Removal Lemma, \cite{MR2299696, MR2604289}]\label{thm:removal}
Let $r, C \in \mathbb{Z}$ and $\varepsilon >0$ be fixed.
For every family of $r$-graphs $\mathcal{F}$ on at most $C$ vertices, there exists $\delta > 0$ such that every sufficiently large $n$, every $r$-graph on $n$ vertices,
which contains at most $\delta n^{|V(F)|}$ induced copies of $F$ for every $F \in \mathcal{F}$,  can be made induced $\mathcal{F}$-free by adding and/or deleting at most $\varepsilon n^r$ hyperedges.
\end{theorem}


Claim~\ref{claim:fa47} together with Theorem~\ref{thm:removal} forces 3-graphs with positive co-degree at least $\frac{4}{7}n$ to be highly structured, which is the core of the proof of Theorem~\ref{J4}.


\begin{proof}[Proof of Theorem~\ref{J4}]

Recall that $\overline{\mathbb{F}}_n$ is the blow-up of the complement of the Fano plane (see Figure~\ref{fig:Fano}) on vertices $x_1,\ldots,x_7$, where each vertex $x_i$ is blown-up to
$\frac{n}{7}$ vertices into the set $X_i$.
Let $x \in X_i$ and $y \in X_j$ for some $i,j \in [7]$. If $i=j$ then $d(x,y) = 0$. If $i \neq j$ then 
$d(x,y) = \frac{4}{7}n$. 
As $J_4 \not\subseteq \overline{\mathbb{F}}_n$, $\mathrm{co^+ex}(n,J_4) \geq  \frac{4}{7}n + o(n)$, which implies $\gamma^+(J_4) \geq 4/7$.

Next, we show that $\gamma^+(J_4) \leq 4/7$. Fix $\beta > 0$ and $\varepsilon > 0$ small enough such  that $24 \varepsilon^{1/4} < \frac{1}{12}  \left(\frac{4}{7}\right)^3 \approx 0.015$ and $48\varepsilon^{1/4} < \beta$. 
We shall fix $n$ sufficiently large such that Claim~\ref{claim:fa47} and Theorem~\ref{thm:removal} imply that every $J_4$-free, $n$-vertex $3$-graph $G$ has a subgraph  $G'$ at edit distance at most $\varepsilon n^3$ from $G$ such that 
every 6-vertex induced subgraph of $G'$ is in $\mathcal{F}$.
We will also take $n$ sufficiently large that Lemma~\ref{edge approx} may be applied with $c = \frac{4}{7}$ and any $n$-vertex $3$-graph $H$ with $\delta_2^+(H) \geq \left(\frac{4}{7} - 48 \varepsilon^{1/4}\right)n$ contains a $K_4^3$. This last condition is possible because $\gamma^+(K_4^3) \leq 0.543$ (see Table~\ref{previous bounds}), and $0.543 <\frac{4}{7} - 48 \varepsilon^{1/4}$ by the choice of $\varepsilon$.

Fix an $n$-vertex $J_4$-free $3$-graph $G_n$ with $\delta_2^+(G_n) \geq \frac{4}{7}n$. Our goal is to show that $\delta_2^+(G_n) < \left( \frac{4}{7} + \beta\right)n$, which will establish  $\gamma^+(J_4) \leq \frac{4}{7}$. The proof will also imply that in fact, the balanced blow-up of the complement of the Fano plane is the asymptotically unique extremal construction.

By Lemma \ref{edge approx} and the choice of $\varepsilon$, $G_n$ contains more than $24 \varepsilon^{1/4} n^3$ hyperedges. 
We can thus apply Lemma~\ref{positive codegree cleanup} to conclude that $G'_n$ contains an $n$-vertex subgraph $G''_n$ with $\delta_2^+(G''_n) \geq \left(\frac{4}{7} - 48 \varepsilon^{1/4}\right)n$. Since $G'_n$ may not be a subgraph of $G_n$, we apply Lemma~\ref{positive codegree cleanup} to the subgraph of $G'_n$ obtained by deleting any $3$-edges that were added to $G_n$ by the application of Theorem~\ref{thm:removal}.
Thus, $G''_n$ contains $K_4^3$, say on vertices $v_1,v_2,v_3,v_4$.
Hence $\{v_1,v_2,v_3,v_4\}$ also spans $K_4^3$ in $G'_n$.
A search through $\mathcal{F}$ shows that there are only two possible subgraphs on $5$ vertices containing $K_4^3$, which we label $A$ and $B$ below.
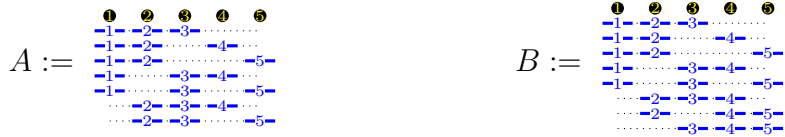
\begin{align*}
A&:=
 \vc{ 
  \begin{tikzpicture}[flag_pic]\outercycle{5}{0}
\draw (x0) node[unlabeled_vertex]{};\draw (x1) node[unlabeled_vertex]{};\draw (x2) node[unlabeled_vertex]{};\draw (x3) node[unlabeled_vertex]{};\draw (x4) node[unlabeled_vertex]{};
\labelvertex{0}\labelvertex{1}\labelvertex{2}\labelvertex{3}\labelvertex{4}\drawhyperedge{0}{5}
\drawhypervertex{0}{0}
\drawhypervertex{1}{0}
\drawhypervertex{2}{0}
\drawhyperedge{1}{5}
\drawhypervertex{0}{1}
\drawhypervertex{1}{1}
\drawhypervertex{3}{1}
\drawhyperedge{2}{5}
\drawhypervertex{0}{2}
\drawhypervertex{1}{2}
\drawhypervertex{4}{2}
\drawhyperedge{3}{5}
\drawhypervertex{0}{3}
\drawhypervertex{2}{3}
\drawhypervertex{3}{3}
\drawhyperedge{4}{5}
\drawhypervertex{0}{4}
\drawhypervertex{2}{4}
\drawhypervertex{4}{4}
\drawhyperedge{5}{5}
\drawhypervertex{1}{5}
\drawhypervertex{2}{5}
\drawhypervertex{3}{5}
\drawhyperedge{6}{5}
\drawhypervertex{1}{6}
\drawhypervertex{2}{6}
\drawhypervertex{4}{6}
\end{tikzpicture} } 
&
B&:=
 \vc{ 
  \begin{tikzpicture}[flag_pic]\outercycle{5}{0}
\draw (x0) node[unlabeled_vertex]{};\draw (x1) node[unlabeled_vertex]{};\draw (x2) node[unlabeled_vertex]{};\draw (x3) node[unlabeled_vertex]{};\draw (x4) node[unlabeled_vertex]{};
\labelvertex{0}\labelvertex{1}\labelvertex{2}\labelvertex{3}\labelvertex{4}\drawhyperedge{0}{5}
\drawhypervertex{0}{0}
\drawhypervertex{1}{0}
\drawhypervertex{2}{0}
\drawhyperedge{1}{5}
\drawhypervertex{0}{1}
\drawhypervertex{1}{1}
\drawhypervertex{3}{1}
\drawhyperedge{2}{5}
\drawhypervertex{0}{2}
\drawhypervertex{1}{2}
\drawhypervertex{4}{2}
\drawhyperedge{3}{5}
\drawhypervertex{0}{3}
\drawhypervertex{2}{3}
\drawhypervertex{3}{3}
\drawhyperedge{4}{5}
\drawhypervertex{0}{4}
\drawhypervertex{2}{4}
\drawhypervertex{4}{4}
\drawhyperedge{5}{5}
\drawhypervertex{1}{5}
\drawhypervertex{2}{5}
\drawhypervertex{3}{5}
\drawhyperedge{6}{5}
\drawhypervertex{1}{6}
\drawhypervertex{3}{6}
\drawhypervertex{4}{6}
\drawhyperedge{7}{5}
\drawhypervertex{2}{7}
\drawhypervertex{3}{7}
\drawhypervertex{4}{7}
\end{tikzpicture} } 
\end{align*}

We now partition $V(G'_n)$ into 7 sets $X_1,\ldots,X_7$ as follows. 
Put $v_i \in X_i$ for $i \in [4]$.
For $v \in V(G'_n) \setminus \{v_1,v_2,v_3,v_4 \}$, define $G_v :=  G_n[v_1,v_2,v_3,v_4,v]$.
If $G_v$ is isomorphic to $A$ then put $v \in X_i$ where $d_{G_v}(v,v_i) = 0$ for some $v_i$.
If $G_v$ is isomorphic to $B$ then put 
$v \in X_5$ if  $v_1v_2v, v_3v_4v \not\in E(G_v)$ or
$v \in X_6$ if  $v_1v_3v, v_2v_4v \not\in E(G_v)$ or 
$v \in X_7$ if  $v_1v_4v, v_2v_3v \not\in E(G_v)$.
This addresses every vertex $v$, so we have a partition of $V(G'_n)$ into $X_1,\ldots,X_7$.
This labeling matches that in Figure~\ref{fig:Fano}.

Now, we show that $G'_n$ is in fact a blow-up of the complement of the Fano plane with classes $X_1, \dots, X_7$.

\begin{claim}\label{no bad edges}
$G'_n$ satisfies the following conditions.

(i) For every $i \in [7]$, no edge of $G'_n$ intersects $X_i$ in more than one vertex.

(ii) If $a \in X_i, b \in X_j, c \in X_k$ for classes $X_i,X_j,X_k$ which correspond to an edge of the Fano plane as labeled in Figure~\ref{fig:Fano}, then $abc \not\in E(G'_n)$.

(iii) If $a \in X_i, b \in X_j, c \in X_k$ for distinct classes $X_i,X_j,X_k$ which do not correspond to an edge of the Fano plane as labeled in Figure~\ref{fig:Fano}, then $abc \in E(G'_n)$.

In particular, $G'_n$ is a blow-up of the complement of the Fano plane.

\end{claim}

\begin{proof}
We prove conditions (i), (ii), and (iii) one by one; in each case, we will argue that if the condition is not satisfied, then $G'_n$ must contain some subgraph that is not in $\mathcal{F}$. Throughout, refer to Figure~\ref{J4 positive density set} for the labeled members of $\mathcal{F}$.
We implemented this check by a computer to reduce the number of cases needed to be done by hand. 

For (i), suppose for a contradiction that there exists an edge $abc$ such that $a,b$ are in the same class of $G'_n$. First note that by the definition of $X_1, \dots, X_7$, if $abc$ either intersects $\bigcup_{i=1}^4 X_i$ in at most two vertices or $abc$ is contained in $X_i$ for some $i \in [4]$, then there exists $\{i,j,k\} \subset [4]$ such that two of $N(x_ix_j)$, $N(x_ix_k)$, and $N(x_jx_k)$ contain $a,b,c$. Up to symmetry, we may assume $a,b,c \in N(x_1x_2)$ and $N(x_1x_3)$. 
The the set of edges of an induced subgraph $H$ of $G'_n$ on $\{x_1,x_2,x_3,a,b,c\}$
includes 
\[E = \{x_1x_2x_3,\,  x_1x_2a,\, x_1x_2b,\, x_1x_2c,\, x_1x_3a,\, x_1x_3b,\, x_1x_3c,\,  abc \}.\]
Since no graph of $\mathcal{F}$ contains a subset of edges isomorphic to $E$, we have $H \not\in \mathcal{F}$, a contradiction.  We use an analogous claim repeatedly. We use computer for verification as well as arguments by hand.

Thus, if edge $abc$ exists, we must have $a,b \in X_i$ and $c \in X_j$ for some $i \neq j$ with $i,j \in [4]$. 
Note that at most one of $a,b,c$ is in $\{x_1,x_2,x_3,x_4\}$, since by the definition of the classes, no edge containing $x_i$ and $x_j$ for $\{i,j\} \in [4]$ intersects $X_i$ or $X_j$.
We have (up to symmetry of the classes) two cases.

 \textbf{Case 1:} $a,b \in X_3$, $c \in X_4$, and $x_3 \not\in \{a,b\}$.
    Define  $H=G'_n[x_1,x_2,x_3, a, b,c]$. 
    The set of edges of $H$ 
includes 
\[E = \{x_1x_2x_3,\,  x_1x_2a,\, x_1x_2b,\, x_1x_2c,\, x_1x_3c,\,  abc \}\] 
and avoids edges $N=\{x_1x_3a,\,  x_1x_3b,\,  x_2x_3a,\,  x_2x_3b\}$.
    Since no graph in $\mathcal{F}$ contains a subset of edges isomorphic to $E$ and avoids $N$, we have $H \not\in \mathcal{F}$, a contradiction.

   \textbf{Case 2:} $a = x_3$, $b \in X_3$, and $c \in X_4$.
    Consider the subgraph $H$ of $G'_n$ induced on $\{x_1,x_2,x_3,x_4, b,c\}$. 
    The set of edges of $H$ 
includes 
\[E = \{x_1x_2x_3,\, x_1x_2x_4,\, x_1x_3x_4,\, x_2x_3x_4,\, x_1x_2b,\, x_1x_4b,\, x_2x_4b,\, x_1x_2c,\, x_1x_3c,\, x_2x_3c,\, abc \}.\] 
    Since no graph in $\mathcal{F}$ contains a subset of edges isomorphic to $E$, we have $H \not\in \mathcal{F}$, a contradiction.


For (ii), suppose for a contradiction that there exists an edge $abc$ spanning three classes that correspond to an edge in Figure~\ref{fig:Fano}. Up to symmetry, there are two cases.


\textbf{Case 1:} The edge $abc$ intersects $\bigcup_{i = 1}^4 X_i$ in two vertices. Without loss of generality, $a \in X_1, b\in X_4$, and $c \in X_7$.
Observe that we cannot have both $a = x_1$ and $b = x_4$ as otherwise $c \not\in X_7$ by definition of the classes, $N(x_1,x_4)$ is disjoint from $X_7$. 
Without loss of generality, $b \neq x_4$. We consider the subgraph $H$ of $G'_n$ induced on $\{x_2,x_3,x_4, a, b, c\}$. 
 The set of edges of $H$ 
includes 
\[E = \{x_2x_3x_4,\, x_2x_3a,\, x_2x_4a,\, x_3x_4a,\, x_2x_3b,\, x_2x_4c,\, x_3x_4c,\, abc \}\] 
and avoids edges containing both $x_4$ and $b$.
    Since no graph in $\mathcal{F}$ contains a subset of edges isomorphic to $E$ while avoiding edges containing both $x_4$ and $b$,  we have $H \not\in \mathcal{F}$, a contradiction.

\textbf{Case 2:} The edge $abc$ does not intersect $\bigcup_{i = 1}^4 X_i$. Without loss of generality, $a \in X_5, b \in X_6,$ and $c \in X_7$.
We define the subgraph $H$ of $G'_n$ induced on $\{x_1,x_2,x_3, a,b,c\}$. 
 The  set of edges of $H$ 
includes 
\[E = \{x_1x_2x_3,\, x_1x_2b,\, x_1x_2c,\, x_1x_3a,\, x_1x_3c,\, x_2x_3a,\, x_2x_3b,\, abc \}\] 
and avoids edges $N=\{x_1x_2a$, $x_2x_3c$, $x_1x_3b\}$.
    Since no graph in  $\mathcal{F}$ contains a subset of edges isomorphic to $E$ while avoiding $N$ we have $H \not\in \mathcal{F}$, a contradiction.

Finally, for (iii), suppose for a contradiction that there exist vertices $a \in X_i, b\in X_j, c \in X_k$ such that $X_i,X_j,X_k$ do not correspond to an edge in Figure~\ref{fig:Fano} and $abc \not\in E(G'_n)$. Up to symmetry of classes, there are three cases.


\textbf{Case 1:}
We have $i,j,k \in [4]$. Without loss of generality, $a \in X_1, b \in X_2,$ and $ c \in X_3$.
By the definition of $X_1,X_2,X_3$, note that at most one of $a,b,c$ is in $\{x_1,x_2,x_3\}$. Without loss of generality, $b \neq x_2$ and $c \neq x_3$. We define the subgraph $H$ of $G'_n$ induced on $\{x_2,x_3,x_4, a, b, c\}$. 
 The set of edges of $H$ 
includes 
\[E = \{x_2x_3x_4,\, x_2x_3a,\, x_2x_4a,\, x_3x_4a,\, x_3x_4b,\, x_2x_4c \}\]
and avoids edges $N = \{abc,\, x_2x_3b,\, x_2x_4b,\,  x_2x_3c,\, x_3x_4c  \}$.
    Since, no graph of $\mathcal{F}$ contains a subset of edges isomorphic to $E$ while avoiding edges in $N$, $H \not\in \mathcal{F}$, a contradiction.

\textbf{Case 2:}
Precisely two of $i,j,k$ are in $[4]$. Without loss of generality, $a \in X_1, b\in X_2$, and  $c \in X_6$.
Observe that by the definition of $X_6$, at most one of $a,b$ is in $\{x_1,x_2\}$; without loss of generality, $b \neq x_2$. We define the subgraph $H$ of $G'_n$ induced on $\{x_2,x_3,x_4,a,b,c\}$. 
 The set of edges of $H$ 
includes 
\[E = \{x_2x_3x_4,\, x_2x_3a,\, x_2x_4a,\, x_3x_4a,\, x_3x_4b,\, x_2x_3c,\, x_3x_4c \}\] 
and avoids edges $N = \{abc,\, x_2x_3b,\, x_2x_4b,\, x_2x_4c\}$.
    Since no graph in $\mathcal{F}$ contains a subset of edges isomorphic to $E$ while avoiding edges in $N$, we have $H \not\in \mathcal{F}$, a contradiction.

\textbf{Case 3:}
Exactly one of $i,j,k$ is in $[4]$. Without loss of generality, $a \in X_1$, $b \in X_5$, and $c \in X_6$.

We define the subgraph $H$ of $G'_n$ induced on $\{x_2,x_3,x_4, a, b, c\}$. 
 The set of edges of $H$ 
includes 
\[E = \{x_2x_3x_4,\, x_2x_3a,\, x_2x_4a,\, x_3x_4a,\, x_2x_3b,\, x_2x_4b,\, x_3x_4b,\, x_2x_3c,\, x_3x_4c \}\] 
and avoids edges $N = \{abc,\, x_2x_4c,\, x_3x_4b \}$.
    Since, no graph of $\mathcal{F}$ contains a subset of edges isomorphic to $E$ while avoiding edges in $N$, $H \not\in \mathcal{F}$, a contradiction.


We conclude that all conditions (i), (ii), and (iii) hold, i.e., $G'_n$ is a blow-up of $\overline{\mathbb{F}}$.
\end{proof}

Finally, we show that $G'_n$ is almost balanced. Recall that $G'_n$ contains a subgraph $G''_n$
with $\delta_2^{+}(G''_n) \geq (\frac{4}{7} - 48 \varepsilon^{1/4})n$. Fix $\alpha \geq 0$ so that a largest class in $G'_n$ has size at least $\left(\frac{1}{7} + \alpha\right) n$. Without loss of generality, $X_1$ is a largest class. We bound the co-degree in $G''_n$ of pairs containing vertices in $X_1$. Observe that there are three sets of classes disjoint from $X_1$ that appear as neighborhoods of vertex pairs in $G''_n$. Namely,
\[N(x_1, x_2) \subseteq X_3 \cup X_4 \cup X_6 \cup X_7;\]
\[N(x_1, x_3) \subseteq X_2 \cup X_4 \cup X_5 \cup X_7;\]
\[N(x_1, x_4) \subseteq X_2 \cup X_3 \cup X_5 \cup X_6.\]
Thus, in $G''_n$, we have
\[d(x_1,x_2) + d(x_1,x_3) + d(x_1,x_4) \leq \sum_{i = 2}^7 2|X_i| \leq \left(\frac{12}{7} - 2\alpha\right)n.\]
By averaging, one of $d(x_1,x_2), d(x_1,x_3), d(x_1,x_4)$ is at most $\left( \frac{4}{7} - \frac{2\alpha}{3} \right)n$. Thus, $\delta_2^+(G''_n) \leq \frac{4n}{7}$, a contradiction if $\delta_2^{+}(G_n) \geq \left( \frac{4}{7} + \beta\right)n$. We conclude that $\gamma^+(J_4) \leq \frac{4}{7}$.
To see that $G''_n$ should be approximately balanced, note that by the minimum positive co-degree condition on $G''_n$, we thus must have $\frac{2\alpha}{3} \leq 48\varepsilon^{1/4}$, i.e., $G''_n$ (and $G'_n$) contains no class of size larger than $\left(\frac{1}{7} + 72 \varepsilon^{1/4}\right)n$. This upper bound implies that $G'_n$ contains no class of size smaller than $\left(\frac{1}{7} - 432 \varepsilon^{1/4}\right)n$.
\end{proof}

We determine the positive co-degree density and the asymptotically unique extremal construction for $F_{4,2}$. 
Since the proof is analogous to (but simpler than) the proof of Theorem~\ref{J4}, we only include a sketch.


\begin{proof}[Sketch of the proof of Theorem~\ref{J4}]
Observe first that $F_{4,2}$ is $6$-partite, so it is not contained in the balanced blow-up of $K_5^3$, which implies $\gamma^+(F_{4,2}) \geq \frac{3}{5}$. We now  show that $\gamma^+(F_{4,2}) \leq \frac{3}{5}$. Let $\mathcal{F}$ be the family of seven $6$-vertex $3$-graphs depicted in Figure~\ref{seven-six-graphs}.

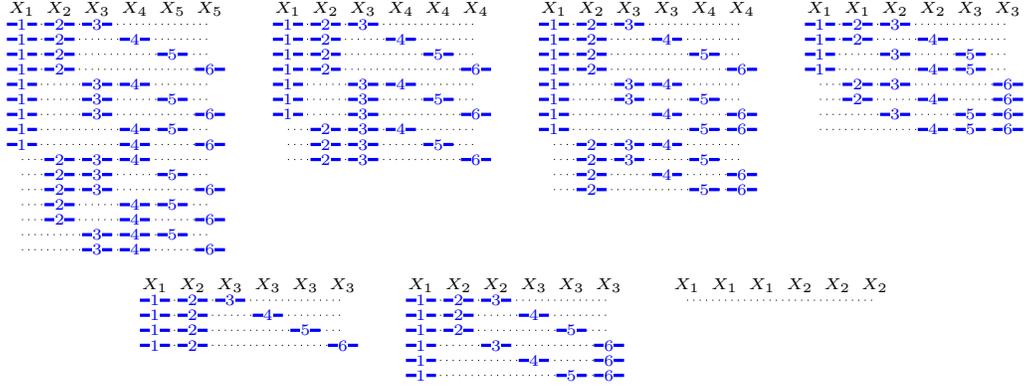
\begin{figure}

\begin{center}
 \vc{ 
  \begin{tikzpicture}[flag_pic]\outercycle{6}{0}
\draw[opacity=0](0,0)--(0,-3.4);
\draw
(x0) node{\tiny$X_1$}
(x1) node{\tiny$X_2$}
(x2) node{\tiny$X_3$}
(x3) node{\tiny$X_4$}
(x4) node{\tiny$X_5$}
(x5) node{\tiny$X_5$}
;
\drawhyperedge{0}{6}
\drawhypervertex{0}{0}
\drawhypervertex{1}{0}
\drawhypervertex{2}{0}
\drawhyperedge{1}{6}
\drawhypervertex{0}{1}
\drawhypervertex{1}{1}
\drawhypervertex{3}{1}
\drawhyperedge{2}{6}
\drawhypervertex{0}{2}
\drawhypervertex{1}{2}
\drawhypervertex{4}{2}
\drawhyperedge{3}{6}
\drawhypervertex{0}{3}
\drawhypervertex{1}{3}
\drawhypervertex{5}{3}
\drawhyperedge{4}{6}
\drawhypervertex{0}{4}
\drawhypervertex{2}{4}
\drawhypervertex{3}{4}
\drawhyperedge{5}{6}
\drawhypervertex{0}{5}
\drawhypervertex{2}{5}
\drawhypervertex{4}{5}
\drawhyperedge{6}{6}
\drawhypervertex{0}{6}
\drawhypervertex{2}{6}
\drawhypervertex{5}{6}
\drawhyperedge{7}{6}
\drawhypervertex{0}{7}
\drawhypervertex{3}{7}
\drawhypervertex{4}{7}
\drawhyperedge{8}{6}
\drawhypervertex{0}{8}
\drawhypervertex{3}{8}
\drawhypervertex{5}{8}
\drawhyperedge{9}{6}
\drawhypervertex{1}{9}
\drawhypervertex{2}{9}
\drawhypervertex{3}{9}
\drawhyperedge{10}{6}
\drawhypervertex{1}{10}
\drawhypervertex{2}{10}
\drawhypervertex{4}{10}
\drawhyperedge{11}{6}
\drawhypervertex{1}{11}
\drawhypervertex{2}{11}
\drawhypervertex{5}{11}
\drawhyperedge{12}{6}
\drawhypervertex{1}{12}
\drawhypervertex{3}{12}
\drawhypervertex{4}{12}
\drawhyperedge{13}{6}
\drawhypervertex{1}{13}
\drawhypervertex{3}{13}
\drawhypervertex{5}{13}
\drawhyperedge{14}{6}
\drawhypervertex{2}{14}
\drawhypervertex{3}{14}
\drawhypervertex{4}{14}
\drawhyperedge{15}{6}
\drawhypervertex{2}{15}
\drawhypervertex{3}{15}
\drawhypervertex{5}{15}
\end{tikzpicture} } 
%
 \vc{ 
  \begin{tikzpicture}[flag_pic]
  \draw[opacity=0](0,0)--(0,-3.4);
\outercycle{6}{0}
\draw
(x0) node{\tiny$X_1$}
(x1) node{\tiny$X_2$}
(x2) node{\tiny$X_3$}
(x3) node{\tiny$X_4$}
(x4) node{\tiny$X_4$}
(x5) node{\tiny$X_4$}
;
\drawhyperedge{0}{6}
\drawhypervertex{0}{0}
\drawhypervertex{1}{0}
\drawhypervertex{2}{0}
\drawhyperedge{1}{6}
\drawhypervertex{0}{1}
\drawhypervertex{1}{1}
\drawhypervertex{3}{1}
\drawhyperedge{2}{6}
\drawhypervertex{0}{2}
\drawhypervertex{1}{2}
\drawhypervertex{4}{2}
\drawhyperedge{3}{6}
\drawhypervertex{0}{3}
\drawhypervertex{1}{3}
\drawhypervertex{5}{3}
\drawhyperedge{4}{6}
\drawhypervertex{0}{4}
\drawhypervertex{2}{4}
\drawhypervertex{3}{4}
\drawhyperedge{5}{6}
\drawhypervertex{0}{5}
\drawhypervertex{2}{5}
\drawhypervertex{4}{5}
\drawhyperedge{6}{6}
\drawhypervertex{0}{6}
\drawhypervertex{2}{6}
\drawhypervertex{5}{6}
\drawhyperedge{7}{6}
\drawhypervertex{1}{7}
\drawhypervertex{2}{7}
\drawhypervertex{3}{7}
\drawhyperedge{8}{6}
\drawhypervertex{1}{8}
\drawhypervertex{2}{8}
\drawhypervertex{4}{8}
\drawhyperedge{9}{6}
\drawhypervertex{1}{9}
\drawhypervertex{2}{9}
\drawhypervertex{5}{9}
\end{tikzpicture} } 
%
 \vc{ 
  \begin{tikzpicture}[flag_pic]
\draw[opacity=0](0,0)--(0,-3.4);
  \outercycle{6}{0}
\draw
(x0) node{\tiny$X_1$}
(x1) node{\tiny$X_2$}
(x2) node{\tiny$X_3$}
(x3) node{\tiny$X_3$}
(x4) node{\tiny$X_4$}
(x5) node{\tiny$X_4$}
;
\drawhyperedge{0}{6}
\drawhypervertex{0}{0}
\drawhypervertex{1}{0}
\drawhypervertex{2}{0}
\drawhyperedge{1}{6}
\drawhypervertex{0}{1}
\drawhypervertex{1}{1}
\drawhypervertex{3}{1}
\drawhyperedge{2}{6}
\drawhypervertex{0}{2}
\drawhypervertex{1}{2}
\drawhypervertex{4}{2}
\drawhyperedge{3}{6}
\drawhypervertex{0}{3}
\drawhypervertex{1}{3}
\drawhypervertex{5}{3}
\drawhyperedge{4}{6}
\drawhypervertex{0}{4}
\drawhypervertex{2}{4}
\drawhypervertex{3}{4}
\drawhyperedge{5}{6}
\drawhypervertex{0}{5}
\drawhypervertex{2}{5}
\drawhypervertex{4}{5}
\drawhyperedge{6}{6}
\drawhypervertex{0}{6}
\drawhypervertex{3}{6}
\drawhypervertex{5}{6}
\drawhyperedge{7}{6}
\drawhypervertex{0}{7}
\drawhypervertex{4}{7}
\drawhypervertex{5}{7}
\drawhyperedge{8}{6}
\drawhypervertex{1}{8}
\drawhypervertex{2}{8}
\drawhypervertex{3}{8}
\drawhyperedge{9}{6}
\drawhypervertex{1}{9}
\drawhypervertex{2}{9}
\drawhypervertex{4}{9}
\drawhyperedge{10}{6}
\drawhypervertex{1}{10}
\drawhypervertex{3}{10}
\drawhypervertex{5}{10}
\drawhyperedge{11}{6}
\drawhypervertex{1}{11}
\drawhypervertex{4}{11}
\drawhypervertex{5}{11}
\end{tikzpicture} } 
%
 \vc{ 
  \begin{tikzpicture}[flag_pic]
 \draw[opacity=0](0,0)--(0,-3.4);
 \outercycle{6}{0}
\draw
(x0) node{\tiny$X_1$}
(x1) node{\tiny$X_1$}
(x2) node{\tiny$X_2$}
(x3) node{\tiny$X_2$}
(x4) node{\tiny$X_3$}
(x5) node{\tiny$X_3$}
;
\drawhyperedge{0}{6}
\drawhypervertex{0}{0}
\drawhypervertex{1}{0}
\drawhypervertex{2}{0}
\drawhyperedge{1}{6}
\drawhypervertex{0}{1}
\drawhypervertex{1}{1}
\drawhypervertex{3}{1}
\drawhyperedge{2}{6}
\drawhypervertex{0}{2}
\drawhypervertex{2}{2}
\drawhypervertex{4}{2}
\drawhyperedge{3}{6}
\drawhypervertex{0}{3}
\drawhypervertex{3}{3}
\drawhypervertex{4}{3}
\drawhyperedge{4}{6}
\drawhypervertex{1}{4}
\drawhypervertex{2}{4}
\drawhypervertex{5}{4}
\drawhyperedge{5}{6}
\drawhypervertex{1}{5}
\drawhypervertex{3}{5}
\drawhypervertex{5}{5}
\drawhyperedge{6}{6}
\drawhypervertex{2}{6}
\drawhypervertex{4}{6}
\drawhypervertex{5}{6}
\drawhyperedge{7}{6}
\drawhypervertex{3}{7}
\drawhypervertex{4}{7}
\drawhypervertex{5}{7}
\end{tikzpicture} } 
%
 \vc{ 
  \begin{tikzpicture}[flag_pic]
  \draw[opacity=0](0,0)--(0,-1.4);
\outercycle{6}{0}
\draw
(x0) node{\tiny$X_1$}
(x1) node{\tiny$X_2$}
(x2) node{\tiny$X_3$}
(x3) node{\tiny$X_3$}
(x4) node{\tiny$X_3$}
(x5) node{\tiny$X_3$}
;
\drawhyperedge{0}{6}
\drawhypervertex{0}{0}
\drawhypervertex{1}{0}
\drawhypervertex{2}{0}
\drawhyperedge{1}{6}
\drawhypervertex{0}{1}
\drawhypervertex{1}{1}
\drawhypervertex{3}{1}
\drawhyperedge{2}{6}
\drawhypervertex{0}{2}
\drawhypervertex{1}{2}
\drawhypervertex{4}{2}
\drawhyperedge{3}{6}
\drawhypervertex{0}{3}
\drawhypervertex{1}{3}
\drawhypervertex{5}{3}
\end{tikzpicture} } 
 \vc{ 
  \begin{tikzpicture}[flag_pic]
  \draw[opacity=0](0,0)--(0,-1.4);
  \outercycle{6}{0}
\draw
(x0) node{\tiny$X_1$}
(x1) node{\tiny$X_2$}
(x2) node{\tiny$X_2$}
(x3) node{\tiny$X_3$}
(x4) node{\tiny$X_3$}
(x5) node{\tiny$X_3$}
;
\drawhyperedge{0}{6}
\drawhypervertex{0}{0}
\drawhypervertex{1}{0}
\drawhypervertex{2}{0}
\drawhyperedge{1}{6}
\drawhypervertex{0}{1}
\drawhypervertex{1}{1}
\drawhypervertex{3}{1}
\drawhyperedge{2}{6}
\drawhypervertex{0}{2}
\drawhypervertex{1}{2}
\drawhypervertex{4}{2}
\drawhyperedge{3}{6}
\drawhypervertex{0}{3}
\drawhypervertex{2}{3}
\drawhypervertex{5}{3}
\drawhyperedge{4}{6}
\drawhypervertex{0}{4}
\drawhypervertex{3}{4}
\drawhypervertex{5}{4}
\drawhyperedge{5}{6}
\drawhypervertex{0}{5}
\drawhypervertex{4}{5}
\drawhypervertex{5}{5}
\end{tikzpicture} } 
%
 \vc{ 
  \begin{tikzpicture}[flag_pic]
  \draw[opacity=0](0,0)--(0,-1.4);
  \outercycle{6}{0}
\draw
(x0) node{\tiny$X_1$}
(x1) node{\tiny$X_1$}
(x2) node{\tiny$X_1$}
(x3) node{\tiny$X_2$}
(x4) node{\tiny$X_2$}
(x5) node{\tiny$X_2$}
;
\drawhyperedge{0}{6}
\end{tikzpicture} 
} 
\end{center}
\caption{The family $\mathcal{F}$ of seven $6$-vertex $3$-graphs.}\label{seven-six-graphs}
\end{figure}

Observe that $\mathcal{F}$ consists of the empty $3$-graph and the $6$-vertex blow-ups of $K_3^3$, $K_3^4$, and $K_5^3$. Thus, $\mathcal{F}$ is exactly the set of $6$-vertex induced subgraphs of the balanced blow-up of $K_5^3$.

Using flag algebras, we show that for sufficiently large $n$, if $G_n$ is a $F_{4,2}$-free $n$-vertex 3-graph with $\delta_2^+(G_n) \geq \frac{3n}{5}$, then
the 6-vertex 3-graphs with positive density in $G_n$ are in $\mathcal{F}$.
The details of the calculations by computer are available at \oururl.
Using the induced removal lemma (Theorem~\ref{thm:removal}), $G_n$ has a small edit distance to $G'_n$, where every 6-vertex subgraph vertices belongs to $\mathcal{F}$.

Recall that $\gamma^+(K_4^3) \leq 0.543 < \frac{3}{5}$, so by Lemma~\ref{positive codegree cleanup}, $G'_n$ contains a subgraph with minimum positive co-degree larger than $(0.543 + \varepsilon)n$ for some $\varepsilon > 0$. Thus, we can assume $G'_n$ contains $K_4^3$, say on vertices $v_1,v_2,v_3,v_4$.
A search through $\mathcal{F}$ shows that there are only two possible subgraphs on $5$ vertices containing $K_4^3$.
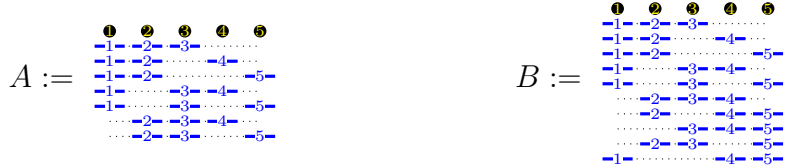
\begin{align*}
A&:=
 \vc{ 
  \begin{tikzpicture}[flag_pic]\outercycle{5}{0}
\draw (x0) node[unlabeled_vertex]{};\draw (x1) node[unlabeled_vertex]{};\draw (x2) node[unlabeled_vertex]{};\draw (x3) node[unlabeled_vertex]{};\draw (x4) node[unlabeled_vertex]{};
\labelvertex{0}\labelvertex{1}\labelvertex{2}\labelvertex{3}\labelvertex{4}\drawhyperedge{0}{5}
\drawhypervertex{0}{0}
\drawhypervertex{1}{0}
\drawhypervertex{2}{0}
\drawhyperedge{1}{5}
\drawhypervertex{0}{1}
\drawhypervertex{1}{1}
\drawhypervertex{3}{1}
\drawhyperedge{2}{5}
\drawhypervertex{0}{2}
\drawhypervertex{1}{2}
\drawhypervertex{4}{2}
\drawhyperedge{3}{5}
\drawhypervertex{0}{3}
\drawhypervertex{2}{3}
\drawhypervertex{3}{3}
\drawhyperedge{4}{5}
\drawhypervertex{0}{4}
\drawhypervertex{2}{4}
\drawhypervertex{4}{4}
\drawhyperedge{5}{5}
\drawhypervertex{1}{5}
\drawhypervertex{2}{5}
\drawhypervertex{3}{5}
\drawhyperedge{6}{5}
\drawhypervertex{1}{6}
\drawhypervertex{2}{6}
\drawhypervertex{4}{6}
\end{tikzpicture} } 
&
B&:=
 \vc{ 
  \begin{tikzpicture}[flag_pic]\outercycle{5}{0}
\draw (x0) node[unlabeled_vertex]{};\draw (x1) node[unlabeled_vertex]{};\draw (x2) node[unlabeled_vertex]{};\draw (x3) node[unlabeled_vertex]{};\draw (x4) node[unlabeled_vertex]{};
\labelvertex{0}\labelvertex{1}\labelvertex{2}\labelvertex{3}\labelvertex{4}
\drawhyperedge{0}{5}
\drawhypervertex{0}{0}
\drawhypervertex{1}{0}
\drawhypervertex{2}{0}
\drawhyperedge{1}{5}
\drawhypervertex{0}{1}
\drawhypervertex{1}{1}
\drawhypervertex{3}{1}
\drawhyperedge{2}{5}
\drawhypervertex{0}{2}
\drawhypervertex{1}{2}
\drawhypervertex{4}{2}
\drawhyperedge{3}{5}
\drawhypervertex{0}{3}
\drawhypervertex{2}{3}
\drawhypervertex{3}{3}
\drawhyperedge{4}{5}
\drawhypervertex{0}{4}
\drawhypervertex{2}{4}
\drawhypervertex{4}{4}
\drawhyperedge{5}{5}
\drawhypervertex{1}{5}
\drawhypervertex{2}{5}
\drawhypervertex{3}{5}
\drawhyperedge{6}{5}
\drawhypervertex{1}{6}
\drawhypervertex{3}{6}
\drawhypervertex{4}{6}
\drawhyperedge{7}{5}
\drawhypervertex{2}{7}
\drawhypervertex{3}{7}
\drawhypervertex{4}{7}
\drawhyperedge{8}{5}
\drawhypervertex{1}{8}
\drawhypervertex{2}{8}
\drawhypervertex{4}{8}
\drawhyperedge{9}{5}
\drawhypervertex{0}{9}
\drawhypervertex{3}{9}
\drawhypervertex{4}{9}
\end{tikzpicture} } 
\end{align*}
Note that $A$ is the (unique) $5$-vertex blow-up of $K_4^3$ and $B = K_5^3$. We now partition $V(G'_n)$ into five sets $X_1,\ldots,X_5$ as follows. 
Put $v_i \in X_i$ for $i \in [4]$.
For $v \in V(G'_n) \setminus \{v_1,v_2,v_3,v_4 \}$, define $G_v :=  G_n[v_1,v_2,v_3,v_4,v]$.
If $G_v$ is isomorphic to $A$, then put  $v \in X_i$, where $d_{G_v}(v,v_i) = 0$.
If $G_v$ is isomorphic to $B$, then put  $v \in X_5$.

An inspection of cases establishes the following claim.

\begin{claim}
$G'_n$ satisfies the following conditions.\\
(i) For every $i \in [5]$, no edge of $G'_n$ intersects $X_i$ in more than one vertex.\\
(ii) For every $a \in X_i, b \in X_j,$ and $c \in X_k$ with $i,j,k$ pairwise distinct, $abc \in E(G'_n)$.
\end{claim}

In particular, $G'_n$ is a blow-up of $K_5^3$. From here, the positive co-degree condition can be used to establish that the vertex-partition of $G_n'$ is essentially balanced.
\end{proof}

Finally, Theorem~\ref{thm:FlagDensity}  was proved using flag algebras. The certificates for the proofs are available at \oururl. These bounds are unlikely to be tight.

\section{Concluding remarks}\label{conclusion}

While we significantly expand the known sets of jumps and achievable values for $\gamma^+$, the general behavior of $\gamma^+$ remains mysterious, even for $r = 3$. It is unclear whether $\gamma^+$ has a jump everywhere, though we conjecture that more jumps exist than are characterized in Theorem~\ref{general jumps}.

\begin{ques}\label{jump question}
For $r \geq 3$, which values of $\alpha \in [\frac{2}{2r - 1}, 1]$ are $\gamma^+$-jumps? 
Does there exist an $\alpha \in [\frac{2}{2r - 1}, 1]$ which is \textit{not} a $\gamma^+$-jump?
\end{ques}

Many more concrete questions could be asked when $r=3$. 
For example, it is unclear, how far Theorem~\ref{densities} is from completely characterizing achievable densities in the range $[0,\frac{1}{2}]$ when $r=3$.

\begin{ques}\label{turan suspension}
Are there achievable values of $\gamma^+$ in $[\frac{2}{5}, \frac{1}{2}]$ that are not of the form $\frac{k-2}{2k - 3}$, when $r=3$? 
\end{ques}

A negative answer to Question~\ref{turan suspension} would suggest some similarity between $\gamma^+$ and $\pi$, since the extremal constructions in Theorem~\ref{densities} are a fairly natural analogue of Tur\'an graphs. However, note that there is no hope for an Erd\H{o}s-Stone-Simonovits-type result giving values of $\gamma^+$. 
Indeed, since a balanced blow-up of $K_4^3$ has positive co-degree density $\frac{1}{2}$, every $3$-graph $F$
with $0 < \gamma^+(F) < \frac{1}{2}$ is $4$-partite. 

To begin addressing either Question~\ref{jump question} or \ref{turan suspension}, it seems natural to start with the next interval between known achievable values. Even this next case seems difficult.

\begin{ques}
For $r = 3$, is every $\alpha \in [\frac{2}{5}, \frac{3}{7})$ a $\gamma^+$-jump? 
Is there a family $\mathcal{F}$ with $\gamma^+(\mathcal{F}) \in \left(\frac{2}{5}, \frac{3}{7}\right)$?
\end{ques}

The use of flag algebras introduces the potential for new approaches to positive co-degree questions, particularly when combined with proving the existence of $\gamma^+$-jumps. As illustrated by Theorems~\ref{J4} and \ref{F42density}, flag algebra calculations have the potential to directly determine values of $\gamma^+$, and even inexact bounds (e.g., $\gamma^+(K_4^3) \leq 0.543$) given by flag algebras are sometimes substantially better and more useful than what seems tractable by hand. When combined with known jumps of the function $\gamma^+$, flag algebra bounds also have the potential to produce exact results. For instance, any $3$-graph $F$ that can be shown via a flag calculation to have $\gamma^+(F) < \frac{2}{5}$ must have $\gamma^+(F) \in \left\{0, \frac{1}{3}\right\}$. Thus, obtaining estimates via flags and ``rounding down'' via known jumps is a very efficient way to determine the densities of many small $3$-graphs. Since we can directly characterize those $3$-graphs with positive co-degree density in $\left\{0, \frac{1}{3}\right\}$, it is now also possible to directly determine by inspection whether a fixed $3$-graph $F$ has $\gamma^+(F) \in \left\{0, \frac{1}{3}\right\}$; however, it seems unlikely that the set 
\[\mathcal{F}(r,d) := \{F \text{ an $r$-graph} : \gamma^+(F) = d\}\]
can be as simply characterized for other values of $d$, even when $r = 3$.  We are interested to see whether further understanding of jumps will include characterizations of this type. If they do not, the potential combination of estimated densities with the theory of jumps is an appealing approach to determining densities exactly.

It is also open whether every achievable value of $\gamma^+$ can be achieved as the density of a single $r$-graph. Theorem~\ref{single graph} shows that $\frac{2}{5}$ is achievable by a single $3$-graph, as is every achievable density known outside the interval $\left(\frac{2}{5}, \frac{1}{2}\right)$.

\begin{ques}
For a fixed $k \geq 5$, is there a $3$-graph $F_k$ such that $\gamma^+(F_k) = \frac{k-2}{2k-3}$? More generally, if $\mathcal{F}$ is a family of $r$-graphs with $\gamma^+(\mathcal{F}) = \alpha$, does there always exist a single $r$-graph $F$ with $\gamma^+(F) = \alpha$?
\end{ques}

Every $3$-graph that is not the subgraph of a (blow-up of a) suspension will have positive co-degree density at least $\frac{1}{2}$. The complete list of known achievable densities at least $\frac{1}{2}$ is as follows: $\frac{1}{2}$ (achieved by $F_{3,2}$), $\frac{4}{7}$ (achieved by $J_4$), $\frac{3}{5}$ (achieved by $F_{4,2}$), and $\frac{2}{3}$ (achieved by the Fano plane). 

\begin{ques}
For $r= 3$, find other achievable values for $\gamma^+$ larger than $\frac{1}{2}$. 
Is there an $\alpha \in \left[\frac{1}{2}, 1\right)$ which is a $\gamma^+$-jump? Is there an $\alpha \in \left[\frac{1}{2}, 1\right)$ for which we can characterize the $3$-graphs with $\gamma^+(F) = \alpha$?
\end{ques}

Very little is known about the $\gamma^+$ function for $r$-graphs when $r \geq 4$.
A natural starting point would be to study extensions of $3$-graphs whose positive co-degree densities are known. For example, the \textit{r-daisy} $\mathcal{D}_r$ is the $6$-edge $r$-graph on $r+2$ vertices whose all six edges contain the same $r-2$ vertices and each pair of the remaining 4 vertices is in one edge. 
There was a recent breakthrough on the Tur\'an density of $r$-daisies~\cite{ellis2024daisies}.
Note that $J_4$ is the $3$-daisy. 

\begin{ques}
What is $\gamma^+(\mathcal{D}_r)$ for $r \geq 4$?
\end{ques}

\section{Acknowledgments}
This work used computing resources at the Center for Computational Mathematics, University of Colorado Denver, including the Alderaan cluster, supported by the National Science Foundation award OAC-2019089.
The authors thank Jan Volec for a discussion on Theorem~\ref{J4}.

Part of the work was done while the first and third authors were at AIM.

\bibliographystyle{abbrvurl}
\bibliography{references.bib}

\begin{thebibliography}{10}

\bibitem{MR2299696}
C.~Avart, V.~R\"odl, and M.~Schacht.
\newblock Every monotone 3-graph property is testable.
\newblock {\em SIAM J. Discrete Math.}, 21(1):73--92, 2007.
\newblock \href {https://doi.org/10.1016/j.endm.2005.06.086}
  {\path{doi:10.1016/j.endm.2005.06.086}}.

\bibitem{BaberTuran}
R.~Baber.
\newblock Turán densities of hypercubes, 2012.
\newblock \href {https://arxiv.org/abs/1201.3587} {\path{arXiv:1201.3587}}.

\bibitem{do-jump}
R.~Baber and J.~Talbot.
\newblock Hypergraphs do jump.
\newblock {\em Combin. Probab. Comput.}, 20(2):161--171, 2011.
\newblock \href {https://doi.org/10.1017/S0963548310000222}
  {\path{doi:10.1017/S0963548310000222}}.

\bibitem{jozsituran}
J.~Balogh.
\newblock The {T}ur\'{a}n density of triple systems is not principal.
\newblock {\em J. Comb. Theory Ser. A}, 100(1):176–180, Oct. 2002.
\newblock \href {https://doi.org/10.1006/jcta.2002.3285}
  {\path{doi:10.1006/jcta.2002.3285}}.

\bibitem{l2norm}
J.~Balogh, F.~C. Clemen, and B.~Lidick\'y.
\newblock Hypergraph {T}ur\'an problems in {$\ell_2$}-norm.
\newblock In {\em Surveys in combinatorics 2022}, volume 481 of {\em London
  Math. Soc. Lecture Note Ser.}, pages 21--63. Cambridge Univ. Press,
  Cambridge, 2022.
\newblock \href {https://doi.org/10.1017/9781009093927.003}
  {\path{doi:10.1017/9781009093927.003}}.

\bibitem{BLP}
J.~Balogh, N.~Lemons, and C.~Palmer.
\newblock Maximum size intersecting families of bounded minimum positive
  co-degree.
\newblock {\em SIAM J. Discrete Math.}, 35(3):1525--1535, 2021.
\newblock \href {https://doi.org/10.1137/20M1336989}
  {\path{doi:10.1137/20M1336989}}.

\bibitem{balogh2024turandensitylongtight}
J.~Balogh and H.~Luo.
\newblock Turán density of long tight cycle minus one hyperedge.
\newblock {\em Combinatorica}, 44(5):949–976, Apr. 2024.
\newblock \href {https://doi.org/10.1007/s00493-024-00099-y}
  {\path{doi:10.1007/s00493-024-00099-y}}.

\bibitem{Bollobascancellative}
B.~Bollob\'as.
\newblock Three-graphs without two triples whose symmetric difference is
  contained in a third.
\newblock {\em Discrete Math.}, 8:21--24, 1974.
\newblock \href {https://doi.org/10.1016/0012-365X(74)90105-8}
  {\path{doi:10.1016/0012-365X(74)90105-8}}.

\bibitem{DaisyBollobas}
B.~Bollob\'as, I.~Leader, and C.~Malvenuto.
\newblock Daisies and other {T}ur\'an problems.
\newblock {\em Combin. Probab. Comput.}, 20(5):743--747, 2011.
\newblock \href {https://doi.org/10.1017/S0963548311000319}
  {\path{doi:10.1017/S0963548311000319}}.

\bibitem{ConlonFox}
D.~Conlon and J.~Fox.
\newblock Graph removal lemmas.
\newblock In {\em Surveys in combinatorics 2013}, volume 409 of {\em London
  Math. Soc. Lecture Note Ser.}, pages 1--49. Cambridge Univ. Press, Cambridge,
  2013.

\bibitem{MR1829685}
A.~Czygrinow and B.~Nagle.
\newblock A note on codegree problems for hypergraphs.
\newblock {\em Bull. Inst. Combin. Appl.}, 32:63--69, 2001.

\bibitem{FanoFuredi}
D.~De~Caen and Z.~F\"uredi.
\newblock The maximum size of 3-uniform hypergraphs not containing a {F}ano
  plane.
\newblock {\em J. Combin. Theory Ser. B}, 78(2):274--276, 2000.
\newblock \href {https://doi.org/10.1006/jctb.1999.1938}
  {\path{doi:10.1006/jctb.1999.1938}}.

\bibitem{ding2023vanishingcodegree}
L.~Ding, H.~Liu, S.~Wang, and H.~Yang.
\newblock Vanishing codegree {T}ur\'{a}n density implies vanishing uniform
  {T}ur\'{a}n density, 2023.
\newblock \href {https://arxiv.org/abs/2312.02879} {\path{arXiv:2312.02879}}.

\bibitem{ellis2024daisies}
D.~Ellis, M.-R. Ivan, and I.~Leader.
\newblock Tur\'an densities for daisies and hypercubes.
\newblock {\em Bull. Lond. Math. Soc.}, 56(12):3838--3853, 2024.
\newblock \href {https://doi.org/10.1112/blms.13171}
  {\path{doi:10.1112/blms.13171}}.

\bibitem{ErSi}
P.~Erd\H{o}s and M.~Simonovits.
\newblock A limit theorem in graph theory.
\newblock {\em Studia Sci. Math. Hungar.}, 1:51--57, 1966.

\bibitem{ErSt}
P.~Erd\H{o}s and A.~H. Stone.
\newblock On the structure of linear graphs.
\newblock {\em Bull. Amer. Math. Soc.}, 52:1087--1091, 1946.
\newblock \href {https://doi.org/10.1090/S0002-9904-1946-08715-7}
  {\path{doi:10.1090/S0002-9904-1946-08715-7}}.

\bibitem{codF32falgas}
V.~Falgas-Ravry, E.~Marchant, O.~Pikhurko, and E.~R. Vaughan.
\newblock The codegree threshold for 3-graphs with independent neighborhoods.
\newblock {\em SIAM J. Discrete Math.}, 29(3):1504--1539, 2015.
\newblock \href {https://doi.org/10.1137/130926997}
  {\path{doi:10.1137/130926997}}.

\bibitem{FalgasK4-}
V.~Falgas-Ravry, O.~Pikhurko, E.~Vaughan, and J.~Volec.
\newblock The codegree threshold of {${K^-_4}$}.
\newblock {\em J. Lond. Math. Soc. (2)}, 107(5):1660--1691, 2023.
\newblock \href {https://doi.org/10.1112/jlms.12722}
  {\path{doi:10.1112/jlms.12722}}.

\bibitem{FALGAS-RAVRY_VAUGHAN_2013}
V.~Falgas-Ravry and E.~R. Vaughan.
\newblock Applications of the semi-definite method to the {T}urán density
  problem for 3-graphs.
\newblock {\em Combinatorics, Probability and Computing}, 22(1):21–54, 2013.
\newblock \href {https://doi.org/10.1017/S0963548312000508}
  {\path{doi:10.1017/S0963548312000508}}.

\bibitem{F5Frankl}
P.~Frankl and Z.~F\"uredi.
\newblock A new generalization of the {E}rd{\H o}s-{K}o-{R}ado theorem.
\newblock {\em Combinatorica}, 3(3-4):341--349, 1983.
\newblock \href {https://doi.org/10.1007/BF02579190}
  {\path{doi:10.1007/BF02579190}}.

\bibitem{K43-extremalfrankl}
P.~Frankl and Z.~F\"uredi.
\newblock An exact result for {$3$}-graphs.
\newblock {\em Discrete Math.}, 50(2-3):323--328, 1984.
\newblock \href {https://doi.org/10.1016/0012-365X(84)90058-X}
  {\path{doi:10.1016/0012-365X(84)90058-X}}.

\bibitem{hg-no-jump}
P.~Frankl and V.~R\"odl.
\newblock Hypergraphs do not jump.
\newblock {\em Combinatorica}, 4(2-3):149--159, 1984.
\newblock \href {https://doi.org/10.1007/BF02579215}
  {\path{doi:10.1007/BF02579215}}.

\bibitem{F32furedi}
Z.~F\"uredi, O.~Pikhurko, and M.~Simonovits.
\newblock The {T}ur\'an density of the hypergraph {$\{abc,ade,bde,cde\}$}.
\newblock {\em Electron. J. Combin.}, 10:Research Paper 18, 7, 2003.
\newblock \href {https://doi.org/10.37236/1711} {\path{doi:10.37236/1711}}.

\bibitem{halfpap2025}
A.~Halfpap.
\newblock personal communication.

\bibitem{halfpap2024positive}
A.~Halfpap, N.~Lemons, and C.~Palmer.
\newblock Positive co-degree density of hypergraphs, 2024.
\newblock \href {https://arxiv.org/abs/2207.05639} {\path{arXiv:2207.05639}}.

\bibitem{halfpapmagnanspanning}
A.~Halfpap and V.~Magnan.
\newblock Positive co-degree thresholds for spanning structures, 2024.
\newblock \href {https://arxiv.org/abs/2409.09185} {\path{arXiv:2409.09185}}.

\bibitem{Kamcev2024}
N.~Kamčev, S.~Letzter, and A.~Pokrovskiy.
\newblock The {T}urán density of tight cycles in three-uniform hypergraphs.
\newblock {\em International Mathematics Research Notices}, 2024(6):4804--4841,
  08 2023.
\newblock \href {https://doi.org/10.1093/imrn/rnad177}
  {\path{doi:10.1093/imrn/rnad177}}.

\bibitem{Keevashsurvey}
P.~Keevash.
\newblock Hypergraph {T}ur\'an problems.
\newblock In {\em Surveys in combinatorics 2011}, volume 392 of {\em London
  Math. Soc. Lecture Note Ser.}, pages 83--139. Cambridge Univ. Press,
  Cambridge, 2011.

\bibitem{lidicky2024c5-}
B.~Lidick\'y, C.~Mattes, and F.~Pfender.
\newblock The hypergraph {T}ur\'{a}n densities of tight cycles minus an edge,
  2024.
\newblock \href {https://arxiv.org/abs/2409.14257} {\path{arXiv:2409.14257}}.

\bibitem{ma2024codegreeturandensity3uniform}
J.~Ma.
\newblock On codegree {T}ur\'an density of the 3-uniform tight cycle
  {$C_{11}$}, 2024.
\newblock \href {https://arxiv.org/abs/2409.02765} {\path{arXiv:2409.02765}}.

\bibitem{MubayiFano}
D.~Mubayi.
\newblock The co-degree density of the {F}ano plane.
\newblock {\em J. Combin. Theory Ser. B}, 95(2):333--337, 2005.
\newblock \href {https://doi.org/10.1016/j.jctb.2005.06.001}
  {\path{doi:10.1016/j.jctb.2005.06.001}}.

\bibitem{F32Mubayi}
D.~Mubayi and V.~R\"odl.
\newblock On the {T}ur\'an number of triple systems.
\newblock {\em J. Combin. Theory Ser. A}, 100(1):136--152, 2002.
\newblock \href {https://doi.org/10.1006/jcta.2002.3284}
  {\path{doi:10.1006/jcta.2002.3284}}.

\bibitem{Mubayicode}
D.~Mubayi and Y.~Zhao.
\newblock Co-degree density of hypergraphs.
\newblock {\em J. Combin. Theory Ser. A}, 114(6):1118--1132, 2007.
\newblock \href {https://doi.org/10.1016/j.jcta.2006.11.006}
  {\path{doi:10.1016/j.jcta.2006.11.006}}.

\bibitem{codegreeconj}
B.~Nagle.
\newblock Tur\'an related problems for hypergraphs.
\newblock In {\em Proceedings of the {T}hirtieth {S}outheastern {I}nternational
  {C}onference on {C}ombinatorics, {G}raph {T}heory, and {C}omputing ({B}oca
  {R}aton, {FL}, 1999)}, volume 136, pages 119--127, 1999.

\bibitem{piga}
S.~Piga, M.~Sales, and B.~Sch\"ulke.
\newblock The codegree {T}ur\'an density of tight cycles minus one edge.
\newblock {\em Combin. Probab. Comput.}, 32(6):881--884, 2023.
\newblock \href {https://doi.org/10.1017/s0963548323000196}
  {\path{doi:10.1017/s0963548323000196}}.

\bibitem{piga2024codegreeturandensity3uniform}
S.~Piga, N.~Sanhueza-Matamala, and M.~Schacht.
\newblock The codegree {T}ur\'an density of $3$-uniform tight cycles, 2024.
\newblock \href {https://arxiv.org/abs/2408.02588} {\path{arXiv:2408.02588}}.

\bibitem{piga2023smallcodegree}
S.~Piga and B.~Schülke.
\newblock Hypergraphs with arbitrarily small codegree {T}ur\'an density, 2023.
\newblock \href {https://arxiv.org/abs/2307.02876} {\path{arXiv:2307.02876}}.

\bibitem{pikhurko2023limitpositiveelldegreeturan}
O.~Pikhurko.
\newblock On the limit of the positive {$\ell$}-degree {T}ur\'an problem.
\newblock {\em Electron. J. Combin.}, 30(3):Paper No. 3.25, 15, 2023.
\newblock \href {https://doi.org/10.37236/11912} {\path{doi:10.37236/11912}}.

\bibitem{Raz07}
A.~A. Razborov.
\newblock Flag algebras.
\newblock {\em J. Symbolic Logic}, 72(4):1239--1282, 2007.
\newblock \href {https://doi.org/10.2178/jsl/1203350785}
  {\path{doi:10.2178/jsl/1203350785}}.

\bibitem{MR2604289}
V.~R\"odl and M.~Schacht.
\newblock Generalizations of the removal lemma.
\newblock {\em Combinatorica}, 29(4):467--501, 2009.
\newblock \href {https://doi.org/10.1007/s00493-009-2320-x}
  {\path{doi:10.1007/s00493-009-2320-x}}.

\bibitem{VeraSos}
V.~T. S\'os.
\newblock Remarks on the connection of graph theory, finite geometry and block
  designs.
\newblock In {\em Colloquio {I}nternazionale sulle {T}eorie {C}ombinatorie
  ({R}oma, 1973), {T}omo {II}}, pages 223--233. Accad. Naz. Lincei, Rome, 1976.

\bibitem{MR177847}
P.~Tur\'an.
\newblock Research problems.
\newblock {\em Magyar Tud. Akad. Mat. Kutat\'o{} Int. K\"ozl.}, 6:417--423,
  1961.
\newblock \href {https://doi.org/10.1007/BF02017934}
  {\path{doi:10.1007/BF02017934}}.

\bibitem{flagmatic}
E.~Vaughan.
\newblock Flagmatic.
\newblock URL: \url{http://lidicky.name/flagmatic/flagmatic.html}.

\bibitem{volec}
J.~Volec.
\newblock personal communication.

\bibitem{Wu}
Z.~Wu.
\newblock Positive co-degree {T}urán number for ${C}_5$ and ${C}_5^{-}$, 2022.
\newblock \href {https://arxiv.org/abs/2212.12815} {\path{arXiv:2212.12815}}.

\end{thebibliography}

\end{document}